\documentclass [12pt, reqno]{amsart}  %

\usepackage{amssymb,amsxtra,amsfonts}
\usepackage{epsf}              %
\usepackage{epsfig}             %
\usepackage{graphics}

\usepackage[colorlinks]{hyperref}

\usepackage{verbatim}  
\usepackage{bbding}   

\def\reE@DeclareMathSymbol#1#2#3#4{%
    \let#1=\undefined
    \DeclareMathSymbol{#1}{#2}{#3}{#4}}
\DeclareSymbolFont{symbolsC}{U}{txsyc}{m}{n}
\SetSymbolFont{symbolsC}{bold}{U}{txsyc}{bx}{n}
\DeclareFontSubstitution{U}{txsyc}{m}{n}
\reE@DeclareMathSymbol{\strictiff}{\mathrel}{symbolsC}{76}

\newcommand\beq{\begin{equation}}
\newcommand\eeq{\end{equation}}
\newcommand\bal{\begin{align*}}
\newcommand\eal{\end{align*}}   
\newcommand\bmx{\left(\begin{matrix}}
\newcommand\emx{\end{matrix}\right)}
\newcommand\bsmx{\left(\begin{smallmatrix}}
\newcommand\esmx{\end{smallmatrix}\right)}

\newcommand{\spq}{/\!\!/}

\providecommand{\spqa}[1]{\underset{#1}{/\!\!/}}
\newcommand{\st}{\ \bigl\vert\ }
\providecommand{\abs}[1]{\lvert#1\rvert}

\providecommand{\im}{\text{\rm Im}}

\def\part#1{\frac{\partial\phantom{q}}{\partial#1}}

\newcommand {\flb}{\lbrack\!\lbrack}
\newcommand {\frb}{\rbrack\!\rbrack}
\newcommand {\flp}{(\!(}
\newcommand {\frp}{)\!)}

\newcommand{\union}{\cup}

\newcommand{\id}{\text{\rm Id}} 
\newcommand{\Id}{\text{\rm Id}}

\newcommand{\Rat}{\mathop{\rm Rat}} 
 
\newcommand{\Lie}{{\mathop{\rm Lie}}}

\newcommand{\Ad}{{\mathop{\rm Ad}}}
\newcommand{\ad}{{\mathop{\rm ad}}}
\newcommand{\rank}{\mathop{\rm rank}}
\DeclareMathOperator{\pr}{pr}

\newcommand{\res}{{\mathop{\rm Res}}}
\newcommand{\Res}{{\mathop{\rm Res}}}
\newcommand{\Prod}{\prod}
\newcommand{\Lambdareplacement}{L}

\newcommand{\tr}{{\mathop{\rm Tr}}}
\newcommand{\Tr}{{\mathop{\rm Tr}}}
\DeclareMathOperator{\Hom}{Hom}         
\DeclareMathOperator{\Aut}{\mathop{\rm Aut}}

\newcommand{\SL}{{\mathop{\rm SL}}}
\newcommand{\PSL}{{\mathop{\rm PSL}}}
\newcommand{\GL}{{\mathop{\rm GL}}}

\newcommand{\SO}{{\mathop{\rm SO}}}

\DeclareMathOperator{\Rep}{\rm Rep}
\renewcommand{\Im}{\mathop{\rm Im}}

\renewcommand{\ker}{\mathop{\rm Ker}}

\newcommand{\End}{\mathop{\rm End}}
\newcommand{\diag}{{\mathop{\rm diag}}}
\newcommand{\diagonals}{{\mathop{\rm diagonals}}}
\newcommand{\mdiagonals}{{\setminus\mathop{\rm diagonals}}}

\newcommand{\hk}{{hyperk\"ahler }}   

\newcommand{\PVI}{{$\text{\rm P}_{\text{\rm VI}}$}}   

\newcommand{\hPVIn}{{$\text{\rm hP}^n_{\text{\rm VI}}$}}   
\newcommand{\hPVn}{{$\text{\rm hP}^n_{\text{\rm V}}$}}  
\newcommand{\hPIVn}{{$\text{\rm hP}^n_{\text{\rm IV}}$}} 
\newcommand{\hPVI}{{$\text{\rm hP}_{\text{\rm VI}}$}}   
\newcommand{\hPV}{{$\text{\rm hP}_{\text{\rm V}}$}}  
\newcommand{\hPIV}{{$\text{\rm hP}_{\text{\rm IV}}$}}

\newcommand{\bfa}{{\bf a}}
\newcommand{\ba}{{\bf a}}

\newcommand{\bd}{{\bf d}}

\newcommand{\bp}{{\bf p}}
\newcommand{\bq}{{\bf q}}

\newcommand{\IB}{\mathbb{B}}
\newcommand{\IC}{\mathbb{C}}

\newcommand{\IF}{\mathbb{F}}
\newcommand{\IG}{\mathbb{G}}

\newcommand{\IL}{\mathbb{L}}
\newcommand{\IM}{\mathbb{M}}
\newcommand{\IN}{\mathbb{N}}
\newcommand{\IP}{\mathbb{P}}

\newcommand{\IZ}{\mathbb{Z}}

\newcommand{\bP}{{\bf P}}

\newcommand{\cA}{\mathcal{A}}
\newcommand{\cB}{\mathcal{B}}

\newcommand{\Conn}{\mathcal{C}onn}

\newcommand{\cE}{\mathcal{E}}

\newcommand{\G}{\mathcal{G}}\newcommand{\cG}{\mathcal{G}}

\newcommand{\cI}{\mathcal{I}}

\newcommand{\cL}{\mathcal{L}}

\newcommand{\cM}{\mathcal{M}}

\newcommand{\cN}{\mathcal{N}}
\newcommand{\cO}{\mathcal{O}}

\newcommand{\cP}{\mathcal{P}}

\newcommand{\cQ}{\mathcal{Q}}       

\newcommand{\cR}{\mathcal{R}}

\newcommand{\cU}{\mathcal{U}}

\newcommand{\g}{       \mathfrak{g}     }

\newcommand{\lt}{\mathfrak{t}}
\newcommand{\lh}{\mathfrak{h}}

\newcommand{\gl}{       \mathfrak{gl}     } 

\newcommand{\h}{\mathfrak{h}}

\newcommand{\wt}{\widetilde}

\newcommand{\wh}{\widehat}

\newcommand{\al}{\alpha}

\newcommand{\be}{\beta}
\newcommand{\ga}{\gamma}
\newcommand{\de}{\delta}
\newcommand{\De}{\Delta}

\newcommand {\eps}{\varepsilon}

\newcommand{\Ga}{\Gamma}

\newcommand{\la}{\lambda}
\newcommand{\La}{\Lambda}

\newcommand{\Si}{\Sigma}
\renewcommand{\th}{\theta}

\renewcommand{\bar}{\overline}

\newcommand{\bla}{{\pmb{\lambda}}}    

\makeatletter
 \newlength{\typesize}
\makeatother

\newlength{\vvoff}
\newlength{\hhoff}

\def\mapright#1{\smash{
        \mathop{\longrightarrow}\limits^{#1}}}

\def\underset#1#2{\ \smash{\mathop{ #2 }\limits_{#1}}\ }

\newcommand{\pf}{\begin{bpf}}

\newcommand{\pfms}{\begin{bpfms}}
\newcommand{\epf}{\end{bpf}\hfill$\square$\\}           
\newcommand{\epfms}{\end{bpfms}\hfill$\square$\\}       

\newcommand{\idea}{\begin{bidea}}

\newcommand{\eidea}{\end{bidea}\hfill$\square$\\}           

\newcommand{\sk}{\begin{bsk}}    

\newcommand{\esk}{\end{bsk}\hfill$\square$\\}           
\newcommand{\sketch}{\begin{bsketch}}

\newcommand{\esketch}{\end{bsketch}\hfill$\square$\\}

\newtheorem {hypo}{\bf\hspace{-\parindent}Hypothesis}
\newtheorem {thm}[hypo]{Theorem}   
\newtheorem {prop}[hypo]{Proposition}

\newtheorem {cor}[hypo]{Corollary}
\newtheorem {lem}[hypo]{Lemma}

\newtheorem {defn}[hypo]{Definition}
\newtheorem{eg}[hypo]{Example}
\theoremstyle{remark}\newtheorem{rmk}[hypo]{Remark}

\numberwithin{hypo}{section}

\numberwithin{equation}{section}

\begin{document}
\title{Simply-laced isomonodromy systems}
\author{Philip Boalch}

\begin{abstract}
A new class of isomonodromy equations will be introduced and shown to admit Kac--Moody Weyl group symmetries. This puts into a general context some results of Okamoto on the 4th, 5th and 6th Painlev\'e equations, and shows where such Kac--Moody Weyl groups 
and root systems occur ``in nature". 
A key point is that one may go beyond the class of affine Kac--Moody root systems.
As examples, 
by considering certain hyperbolic Kac--Moody Dynkin diagrams,
we find there is a sequence of higher order Painlev\'e systems lying over each of the classical Painlev\'e equations.
This leads to a conjecture about the Hilbert scheme of points on some Hitchin systems. 
\end{abstract}

\maketitle

\tableofcontents

\section{Introduction}

The aim of this article is to introduce and study a new system of isomonodromic deformation equations.
The best known isomonodromy equations are the Schlesinger equations
\cite{schles-icm1908} 
controlling deformations of Fuchsian systems on the Riemann sphere.
Geometrically these equations constitute  a nonlinear flat connection 
on a bundle
$$\cM^*\times \IB\to \IB$$
over a space of parameters (the ``times'') $\IB\cong \IC^m\setminus\diagonals$, where 
$$\cM^*\cong \left(\cO_1\times\cdots \times \cO_m\right)\spq \GL_n(\IC)$$
is the symplectic quotient of a product of coadjoint orbits.
This nonlinear connection may be interpreted as a nonabelian analogue of the Gauss--Manin connection (cf. \cite{smid} \S7) and admits degenerations into Hitchin-type integrable systems 
(cf. \cite{garn1919}).
Thus in general one may view an ``isomonodromy system'' as a system of nonlinear differential equations obtained by deforming a Hitchin integrable system (and whose solutions will involve more complicated functions, such as the Painlev\'e transcendents, than the abelian or theta functions involved in solving such integrable systems).

The next simplest family of isomonodromy equations are due to Jimbo--Miwa--M\^ori--Sato (JMMS \cite{JMMS}) and arose as equations for correlation functions of the quantum nonlinear Schr\"odinger equation.
The JMMS equations are the isomonodromic deformation equations for linear differential systems of the form
\beq \label{eq: JMMS system1}
 \frac{d}{dz}- \left(T_0+\sum_1^m \frac{R_i}{z-t_i}\right)
\eeq
(where $T_0$ is a diagonal matrix) 
having an irregular singularity at $z=\infty$. 
A remarkable feature of the JMMS equations is that there are now two sets of times; one may deform the pole positions (the $t_i$, as in the Fuchsian case) as well as the eigenvalues of $T_0$ (the ``irregular times''). 
Further Harnad \cite{Harn94} has shown that the JMMS equations admit a symmetry under which one may {\em swap} the roles of the two sets of times.

In this article we will write down and study the ``next simplest'' class of isomonodromy equations, which are even more symmetric: in effect the two set of  times are extended to $k$ sets of times, all of which may be permuted, and Harnad's discrete duality is extended to an action of a continuous $\SL_2(\IC)$ symmetry group.
(Further, a Hamiltonian formulation will be given enabling the definition of new $\tau$ functions.)

One motivation to study such symmetric isomonodromy systems was to better understand and generalise the affine Weyl symmetry groups of the Painlev\'e equations.
In effect the Painlev\'e equations are the simplest examples of isomonodromy equations: they are the second order nonlinear differential equations which arise as the explicit form of
the isomonodromy connection when the fibres have dimension two, i.e. $\dim_\IC\cM^*=2$.
Okamoto has shown that the six Painlev\'e equations admit certain affine Weyl symmetry groups (\cite{OkaP24, OkaPVI, OkaPV, Okamoto-dynkin}), and 
 this is intriguing since the underlying root systems are not immediately apparent from the geometry. 
For example the simplest nontrivial case of the Schlesinger equations is equivalent to the Painlev\'e VI equation, which has symmetry group the affine Weyl group of type $D_4$.
However no link to the loop group of $\SO_8$ is manifest (one starts with a Fuchsian system with four poles on a rank two bundle). Similarly Painlev\'e V has affine $A_3$ symmetry group and  Painlev\'e IV has affine $A_2$ symmetry group (the other cases are not simply-laced and will be ignored here for simplicity).

\begin{figure}[h]
	\centering
	\input{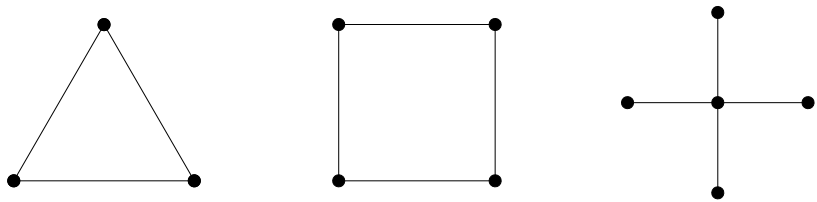}
	\caption{Affine Dynkin diagrams for Painlev\'e equations IV, V and VI.}\label{fig: p456}
\end{figure}

In the second part of this article  (\S\ref{sn: reductions}) we will single out a large class of graphs which contains these three graphs, and attach an isomonodromy system to each such graph (and some data on the graph), such that the corresponding Kac--Moody Weyl group acts by symmetries, relating the corresponding isomonodromy systems. 
Thus any such graph appears as a ``Dynkin diagram'' for an isomonodromy system.
The class of graphs for which this result holds contains all the complete $k$-partite graphs for any $k$ (and in particular all the complete graphs).
For example one may ask why there is no second order Painlev\'e equation attached to the pentagon (the affine $A_4$ Dynkin diagram): from our viewpoint this is because it is not a complete $k$-partite graph for any $k$, whereas the triangle and the square are (as is the four-pointed star).

A simple corollary of this way of thinking is that we are able to canonically attach an isomonodromy system of order $2n$  (i.e. $\dim_\IC\cM^*=2n$) to each of the six Painlev\'e equations for each integer $n=1,2,3,4, \ldots$ so that the $n=1$ case is isomorphic to the original Painlev\'e system. (These ``higher Painlev\'e systems'' look to be completely different to the well-known ``Painlev\'e hierarchies'').

\ 

\begin{figure}[h]
	\centering
	\input{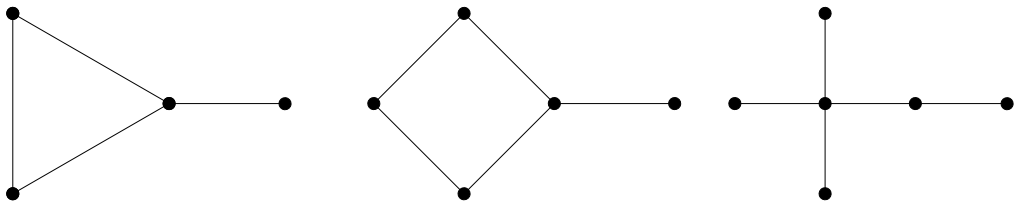}
	\caption{Dynkin diagrams for simply-laced higher Painlev\'e systems.}\label{fig: hp456}
\end{figure}

In the remainder of this introduction we will recall more background, summarise the main result, 
and describe the JMMS equations from our point of view, serving as the main prototype for the extension we have in mind.

\subsection{Further remarks and viewpoints}

1) In 1981 a quite general isomonodromy system was established by Jimbo--Miwa--Ueno \cite{JMU81}, controlling isomonodromic deformations of linear differential systems whose most singular coefficient has distinct eigenvalues at each pole.
This work was revisited from a moduli theoretic viewpoint in  \cite{smid} and the symplectic nature of these (JMU) equations was established.
More general moduli spaces were then constructed, without the distinct eigenvalue condition, in \cite{wnabh} (in fact in arbitrary genus and with compatible parabolic structures and stability conditions) and shown to be hyperk\"ahler manifolds. 
The distinct eigenvalue condition implies that in general the JMU equations are not symmetric, and so here we go back and generalise in a different direction the earlier viewpoint of \cite{JMMS}.  
This is guided by the descriptions in \cite{JMMS, Harn94, smid} of the moduli spaces involved as symplectic quotients, enabling us to see that some of them are isomorphic to (Nakajima) quiver varieties, and thus how the graphs arise from the geometry (cf. 
\cite{CB-additiveDS} in the Fuchsian case and \cite{quad} Exercise 3 for early examples in the irregular case).
A new viewpoint we have found very useful here is to identify certain quiver varieties (or twisted versions of them) as moduli spaces of presentations of modules for the first Weyl algebra.

2) Noumi--Yamada (\cite{noumi-yamada-cmp98} Problem 5.1) considered the problem of finding systems of nonlinear differential equations for each affine root system (or more generally for each generalised Cartan matrix) on which the corresponding Weyl group acts as Backlund transformations.
As an example in \cite{noumi-yamada-Al1} (see also Noumi's ICM talk \cite{noumi-icm})
they wrote down a sequence of isomonodromy equations for each of the type $A$ affine Weyl groups, generalising the  Painlev\'e IV and V equations. See also various articles of Y. Sasano such as \cite{sasano-4d}. 
(Note they have also found \cite{NY-so8} a direct link between \PVI\  and $\SO(8)$.)
Their work is in a sense orthogonal to ours: although we consider a much larger class of Kac--Moody root systems, whereas they only construct systems in certain affine Kac--Moody cases, the intersection of the set of our graphs with theirs is amongst the usual second order Painlev\'e equations. (I do not know if there is some common generalisation.)

3) The graphs we are using give a reasonably efficient way to start to classify some isomonodromy systems (or at least to tell when they might be isomorphic). Some four dimensional examples were found in this way in \cite{rsode} by looking amongst the hyperbolic Kac--Moody graphs (the next simplest class after the affine ones). 
Not every isomonodromy system will be simply-laced though of course. 
Some other possible approaches to such questions are as follows.
i) Recall that Cosgrove (e.g. in \cite{cosgrove-1}) has done much work classifying higher order equations with the Painlev\'e property although his approach seems difficult for higher rank equations (e.g. sixth order).
ii) On the other hand Malgrange \cite{malgrange--galois-feuil} and 
Umemura \cite{umemura--dgtinfd} have introduced a nonlinear differential Galois theory which one might hope would help with such classification, although it turns out that the strength of their theory is its vast generality and present results suggest it does not even distinguish amongst the six Painlev\'e equations themselves (when they have generic parameters).\footnote{One could see this attempt to classify such nonlinear algebraic differential equations as a basic step in ``differential algebraic geometry''---the extension of algebraic geometry obtained by allowing derivatives in the equations---the isomonodromy equations should form a basic class of objects to be studied, much as abelian varieties or Calabi--Yau manifolds are in classical algebraic geometry. On the other hand one can view this subject as a half-way step to ``noncommutative algebraic geometry''.}

4) In essence we are associating what might be called a ``wild nonabelian Hodge structure'' to a certain class of graphs with some extra data on them. This structure consists of a hyperk\"ahler manifold $\cM$ (as in \cite{wnabh}) which in one complex structure is a moduli space of meromorphic connections and in another is a space of meromorphic Higgs bundles. The isomonodromy system is naturally associated to this structure (it controls the isomonodromic deformations of the meromorphic connections). It seems that the work of S. Szabo \cite{szabo-nahm} (interpreting certain Fourier--Laplace transforms as Nahm transforms) can be extended to show that the full hyperk\"ahler metric is also preserved by all the Kac--Moody Weyl group symmetries.

\begin{figure}[h]
	\centering
	\input{partite3.pstex_t}
	\caption{Complete $k$-partite graphs from partitions of $N\le 6$}\label{fig: graph table}
\center{\ (omitting the stars $\cG(1,n)$ and the totally disconnected graphs $\cG(n)$)}\label{fig: partite-intro}
\end{figure}

\newpage
\subsection{Summary of main result}

The main result involves a class of graphs that we will call {\em supernova graphs} (cf.  Definition \ref{defn: supernova graph}), generalising the class of star-shaped graphs.
In brief 
$\wh\cG$ is a supernova graph with nodes $\wh I$ if for some $k$ there is a complete $k$-partite subgraph $\cG\subset  \wh \cG$ with nodes $I\subset \wh I$, and 
$\wh \cG$ is obtained from $\cG$ by gluing on some legs.
In turn recall a complete $k$-partite graph is a graph whose nodes 
$$I = \bigsqcup_{j\in J} I_j$$
are partitioned into parts parameterised by a set $J$ with $\abs{J}=k$, and two nodes of $\cG$ are joined by a single edge if and only if they are not in the same part.

Let $\wh \cG$ be a supernova graph and choose data $\bd, \bla, \ba$, consisting of:

1)  An integer $d_i \ge 0$ for each node  $i\in\wh I$,

2) A scalar $\la_i\in \IC$ for each $i\in\wh I$, 
such that $\sum \la_id_i=0$,

3) A distinct point $a_j$ of the Riemann sphere for each part $j\in J$.

\begin{thm}\label{thm: intro thm}
$\bullet$ 
There is an isomonodromy system 
\beq\label{eq: slimd}
\cM_{st}^*(\bla, \bd, \ba)\times \IB\to \IB
\eeq
where the space of times is 
$\IB\cong \Prod_{j\in J} (\IC^{\abs{I_j}}\setminus \diagonals).$
It controls isomonodromic deformations of certain linear differential systems on bundles of rank
\beq\label{eq: ranks}
\sum_{i\in I\setminus I_\infty} d_i,\eeq
where $I_\infty\subset I$ is the part with $a_j=\infty$. 

$\bullet$ If $i\in \wh I$ is a node of $\wh \cG$ and  $\la_i\ne 0$ then there is an isomorphism of isomonodromy systems
$$
\cM_{st}^*(\bla, \bd, \ba)\times \IB \cong 
\cM_{st}^*(r_i(\bla), s_i(\bd),  \ba)\times \IB $$
where $s_i,r_i$ are the simple reflections (and dual simple reflections) generating the Weyl group of the Kac--Moody root system attached to the graph $\wh \cG$.

$\bullet$ If $g\in \SL_2(\IC)$ then there is an isomorphism 
$$\cM_{st}^*(\bla, \bd, g(\ba))\times \IB \cong 
  \cM_{st}^*(\bla, \bd, \ba)\times \IB $$
of isomonodromy systems, where $g$ acts on $\ba$ by diagonal M\"obius transformations.
\end{thm}

Also, in \S\ref{sn: +iDSP}, precise criteria involving the Kac--Moody root system will be established for when the spaces $\cM_{st}^*$ are nonempty---this is an additive, irregular analogue of the Deligne--Simpson problem. 
Notice in particular that the action of $\SL_2(\IC)$ enables us to change the ranks appearing in \eqref{eq: ranks}, by moving different points $a_j$ to $\infty$.
Once we have this possibility, it is straightforward to obtain all the reflections geometrically.
(This extends the viewpoint of \cite{pecr, k2p} which derived the action of the Okamoto symmetries of Painlev\'e VI on linear monodromy data from Harnad duality/Fourier--Laplace.)

\subsection{Prototype: The lifted JMMS equations}

By writing the (lifted) JMMS equations in a slightly novel way
we will show how graphs appear naturally. 
In particular we will explain 
how the square (which Okamoto attached to 
Painlev\'e V) emerges from the viewpoint of \cite{JMMS} Appendix 5. 
Our main result is obtained by extending this to include the triangle.

Choose two finite dimensional complex vector spaces 
$W_0,W_\infty$ and consider the space $\IM = \Hom(W_0,W_\infty)\oplus \Hom(W_\infty,W_0)$.
The (lifted) JMMS equations govern how $(P,Q)\in \IM$ should vary with respect to some times $T_0,T_\infty$, and we will write these equations in the following form:

\beq \label{eq: altJMMSintro}
\begin{split} 
dQ = Q\wt{PQ} +\wt{QP}Q  + T_0 Q dT_\infty + dT_0 Q  T_\infty  \\
-dP = P\wt{QP} + \wt{PQ}P + T_\infty P dT_0 + dT_\infty P  T_0.
\end{split}
\eeq

Here $T_0, T_\infty$ are semisimple matrices 
$T_0\in \End(W_0), T_\infty\in \End(W_\infty)$
which may have repeated eigenvalues, but are restricted so that no further eigenvalues are allowed to coalesce, and the corresponding eigenspace decompositions of $W_0,W_\infty$ are held constant.
The tilde operation appearing in \eqref{eq: altJMMSintro} is defined as follows: if $R\in \End(W_i)$ then 
$$\wt R := \ad^{-1}_{T_i}[dT_i,R]$$
which is a one-form with values in the image in $\End(W_i)$ 
of $\ad_{T_i}$. (To see this makes sense note that $\ad_{T_i}$ is invertible when restricted to its image.)
These equations are known to have the Painlev\'e property and special cases include the fifth and sixth Painlev\'e equations. Indeed for example in a special case they imply 
the Schlesinger equations: if one restricts to the case when $T_0=0$ then \eqref {eq: altJMMSintro}
simplifies to
$$dQ = Q\wt{PQ},\qquad dP =  -\wt{PQ}P.$$

It follows immediately that if we write 
$T_\infty = \sum t_i\Id_i$ where $\Id_i$ is the idempotent for the $i$th eigenspace of $T_\infty$ and set
$R_i = Q\Id_i P \in\End(W_0)$ then
$$dR_i = -\sum_{j\ne i} [R_i,R_j]d\log(t_i-t_j)$$
which are the Schlesinger equations.
They govern the isomonodromic deformations of a logarithmic connection on the vector bundle
$W_0\times \IP^1\to\IP^1$.
In general the JMMS equations govern the isomonodromic deformations of a meromorphic connection (corresponding to the differential system \eqref{eq: JMMS system1}) 
on the vector bundle
$W_0\times \IP^1\to\IP^1$ having an irregular singularity of 
Poincar\'e rank one and arbitrarily many logarithmic singularities.
Although quite general nonlinear 
isomonodromy equations have been written down by Jimbo--Miwa--Ueno \cite{JMU81}, the JMMS equations are notable since they have the following symmetry
(now quite transparent in the way we are writing the equations):

\begin{thm}[Harnad duality \cite{Harn94}]
The permutation
$$
(W_0,W_\infty,P,Q,T_0,T_\infty)
\mapsto
(W_\infty,W_0,Q,-P,-T_\infty,T_0)$$
preserves the JMMS equations.
\end{thm}

This is remarkable since it implies that the {\em same equations} also control isomonodromic deformations of a connection on the vector bundle $W_\infty\times \IP^1\to\IP^1$, which of course will in general have different rank.

Our basic aim is to show that the JMMS equations have a natural generalisation that admits an enriched symmetry group.
To describe the picture let us first describe the combinatorics of Harnad's duality in terms of graphs.
First consider the graph with two nodes labelled by $0$ and $\infty$ connected by a single edge. 
We put the vector space $W_j$ at the node $j$ for $j=0,\infty$ and view the maps $P,Q$ as maps in both directions along the edge.
At this level (up to a sign) Harnad duality simply flips over the graph.

$$
	\input{bipartite1.pstex_t}
$$

The next step (which enables us to see the link to Dynkin graphs) is to refine the above graph by splaying each node according to the eigenspaces of the times $T_0,T_\infty$.
In other words suppose $W_j = \bigoplus_{i\in I_j}V_i$ is the eigenspace decomposition of $T_j$ for $j=0,\infty$. Then we break up the node corresponding to $W_0$ into $\abs{I_0}$ nodes and thereby splay the graph (and similarly at the other node), as in the diagram below (for  the case $\abs{I_0}=3, \abs{I_\infty} =2$):

\ 

\begin{figure}[h]
	\centering
	\input{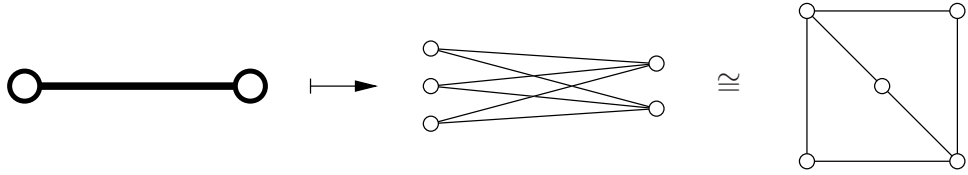}
	\caption{Splaying both nodes.}\label{fig: bipartite-splay}
\end{figure}

Note that the class of refined graphs which arise in this way are precisely the {\em complete bipartite graphs}.
Now the fifth Painlev\'e equation is equivalent to a basic case of the JMMS equations (appearing in the title of \cite{JMMS}): it occurs when both $W_0$ and $W_\infty$ have dimension two and each time $T_j$ has two distinct eigenvalues. 
The crucial observation then is that in this case the refined graph is a square, i.e. the affine Dynkin diagram $A_3^{(1)}$ that Okamoto associated to Painlev\'e V appears almost directly from the work of JMMS.

\ 

\begin{figure}[h]
	\centering
	\input{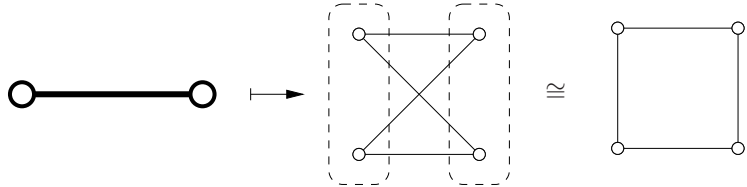}
	\caption{How $A_3^{(1)}$ appears in the graphical approach to the JMMS system.}\label{fig: square}
\end{figure}

This is more than a coincidence since, as we will confirm, Harnad's duality (and the well-known Schlesinger/Backlund transformations) yield the Okamoto symmetries.
Moreover we see how to put this in a more general context since {\em any} complete bipartite graph also appears in the same way. 

The final step ``reduction'' (see \S\ref{sn: reductions}) is to choose an adjoint orbit $\breve\cO_i\subset \End(V_i)$ for each node of the refined graph and quotient by the symmetry group $\prod\GL(V_i)$.
In terms of graphs this will correspond to gluing a leg (a type $A$ Dynkin graph) onto each node of the refined graph: the resulting graph is the Dynkin diagram of the Kac--Moody root system whose Weyl group naturally acts on the nonlinear equations.

\subsection{Generalisation}

The generalisation we have in mind 
involves replacing the initial graph in the JMMS story above (the interval) by 
an {\em arbitrary complete graph}.
(Recall that the complete graph with $k$ nodes  has exactly one edge between any two vertices.)

\ 

\begin{figure}[h]
	\centering
	\input{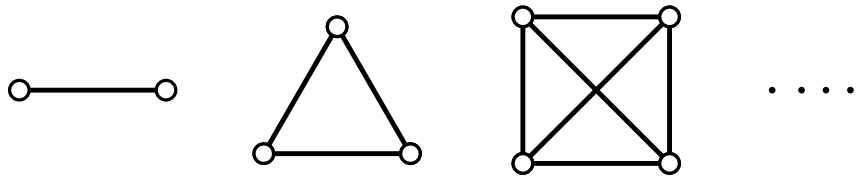}
	\caption{The sequence of complete graphs.}\label{fig: complete graphs}
\end{figure}

Each node of the graph  will be labelled by a distinct element 
$a_j\in \IC\cup\{\infty\}$ of the Riemann sphere.
Thus if $J$ is the set of nodes we may view $J$ as a subset of the Riemann sphere. 
(We will see below that $\infty$ plays a distinguished role.)
In the case of the JMMS equations above we had $J=\{0,\infty\}$.
As above we attach a vector space $W_j$ to each node $j\in J$ and consider a set of times which are semisimple elements $T_j\in \End(W_j)$ for each $j\in J$.
The set of unknown variables  in the nonlinear equations, 
generalising $(P,Q)$ above,
again consists of the linear maps in both directions along each edge:
$$B_{ij}\in \Hom(W_j,W_i),\qquad\text{for all $i\ne j \in J$}.$$
The integrable nonlinear equations we will find governing these then have the form:

\begin{align*} %
dB_{ij} =\sum_{k\in J} 
\wt{X_{ik}B_{ki}}B_{ij} 
+B_{ij}\wt{B_{jk}X_{kj}}+dT_iX_{ik}B_{kj}
+ B_{ik}X_{kj}dT_j
-X_{ik}dT_k X_{kj}/\phi_{ij}
\end{align*}

\noindent 
plus some terms linear in $B_{ij}$ that we will neglect in the introduction, where $\phi_{ij}$ is a complex number depending only on the embedding $J\hookrightarrow \IC\cup\{\infty\}$  and $X_{ij} = \phi_{ij}B_{ij}$.
(Observe that the three quadratic terms are absent in the original JMMS equations.)

The remainder of the story is then similar to above: we splay the nodes of the complete graph according to the eigenspaces of the times $T_j$ to obtain a refined graph.
(The class of graphs which arise in this way is exactly the class of complete $k$-partite graphs.)
Then we reduce as above, gluing on some legs, to obtain a graph whose Kac--Moody Weyl group acts. 
The simplest case, with $k=3$ and each $W_j$ of dimension one, is equivalent to the fourth Painlev\'e equation (and the refined graph, which in this case equals the unrefined graph, is the triangle, agreeing with Okamoto).

\

\noindent
{\bf Acknowledgments.}
\noindent
The basic picture (in particular the link between a large class of Kac--Moody root systems and irregular connections) appeared in the  preprint \cite{rsode} (arxiv June 2008). 
Here we go further by writing down the explicit isomonodromy system, giving its Hamiltonian formulation and showing it is symmetric\footnote{The present article basically subsumes \cite{rsode}, which is not planned for publication, although \cite{rsode} still contains some further motivation, examples and directions.}.
The author gave talks about these further results in Paris September 2010 and St. Petersburg June 2011.
Thanks are due to  H. Nakajima for some helpful correspondence.
Some of this work was done whilst visiting the IHES February 2010.
The author is partially supported by ANR grants 
08-BLAN-0317-01/02 (SEDIGA), 09-JCJC-0102-01 (RepRed).

\section{Modules over the Weyl algebra}

This section describes the class of modules over the first Weyl algebra that we wish to consider. By taking this viewpoint, rather than just that of meromorphic connections, more symmetries are apparent. (Appendix \ref{apx: Harnad}  explains how this generalises Harnad's duality.) 
Choose three $n\times n$ complex matrices $\al,\be,\ga$
and consider the matrix
$$M=\al \partial + \be z - \ga$$
with values in the first Weyl algebra $\cA_1=\IC\langle z,\partial\rangle$, where $\partial = d/dz$.

We wish to study (local) systems of differential equations of the form
$$M(\wh v) = 0$$
for a vector $\wh v$ of holomorphic functions. This may be rephrased as the local system of holomorphic solutions of the (left) $\cA_1$-module $\cN$ defined by the exact sequence

$$\cA^n_1\quad\mapright{(\cdot M)}\quad \cA^n_1\  \to \ \cN\ \to\  0.$$

We will restrict to the case where $\al$ and $\be$ are commuting diagonalisable matrices whose kernels intersect only at zero. This class of modules is clearly preserved under the Fourier--Laplace transform, and in fact by more general symplectic transformations
\beq\label{eq: spt}
(\partial,z) \mapsto (a\partial + b z, c\partial + d z)
\eeq
where $ad-bc=1$:

\begin{lem}
If $\al$ and $\be$ are two commuting semisimple matrices whose kernels intersect only at zero, then so are
$$
a\al +c \be\qquad\text{and}\qquad
b\al +d \be
$$ for any complex numbers $a,b,c,d$ with $ad-bc=1$.
\end{lem}
\pf
In each joint eigenspace the corresponding pair of eigenvalues $(\al_i,\be_i)$ is a nonzero point of $\IC^2$. We are acting on this element of $\IC^2$ by an invertible matrix so it will remain nonzero.
\epf

Thus such a module determines a finite number of points $a_i:=[-\be_i:\al_i]\in \bP:=\IP^1(\IC)$ of the projective line (the minus sign appearing here will be useful later). 
We label this sphere $\bP$ (and call it the ``Fourier sphere'') 
to avoid confusion later with the Riemann sphere $\IP^1_z$ on which $z$ is a local coordinate.
Let $J$ be the set of such points of $\bP$, and let $k=\abs{J}$ be the number of points so obtained (which is in general less than the number of joint eigenspaces, since two pairs of eigenvalues may differ by an overall scalar).
For each point $j\in J$ there is an associated subspace 
$W_j\subset V:=\IC^n$ 
(the joint eigenspaces having pairs of eigenvalues lying over $j$), so that
$V = \bigoplus_{j\in J} W_j$.

We will further assume for each $j$ that the component of $\gamma$ in $\End(W_j)$  is semisimple.%

Let $\infty = [1:0]$ be the point of $\bP$ where $\al_i=0$, corresponding to the kernel of $\al$ (we are not assuming $\infty \in J$; this kernel may be trivial).
Any other point of $\bP$ corresponds to a complex number 
($a_i=-\be_i/\al_i\in \IC$); we identify $\bP$ with $\IC\union\{\infty\}$ in this way.

By multiplying $M$ on the left by a constant invertible matrix, we may then normalize $M$ uniquely so it is of the form
\beq\label{eq: normn}
M= \left(\begin{matrix}
0 &  \\ 
 & 1 
\end{matrix}\right)\partial +
\left(\begin{matrix}
1 & 0 \\ 
0 & - A
\end{matrix}\right) z-\ga
\in \End(W_\infty\oplus U_\infty)\otimes \cA_1
\eeq
where
$U_j = V\ominus W_j=\bigoplus_{i\in J\setminus\{j\}} W_i$.
Here $A=\sum_{j\neq \infty} a_j\Id_j\in \End(U_\infty)$
where $\Id_j$ is the idempotent for $W_j\subset V$. 
The decomposition of $V$ allows us to decompose
$\End(V)$
as 
$$\End(V) = \bigoplus_{i\ne j\in J} \Hom(W_i,W_j) \oplus 
\bigoplus_{j\in J} \End(W_j)$$
and so we may decompose 
$$\ga = \ga^\circ +\delta(\ga)$$
with respect to this decomposition.
We will also write (after normalisation)
\beq\label{eq: Ga decomp}
\delta(\ga) = \wh T = \bmx C & \\ & T \emx,\qquad
\ga^\circ = \Ga =
\left(\begin{matrix}
0 & P \\ 
Q & B
\end{matrix}\right)
 \in \End(W_\infty\oplus U_\infty).
\eeq
Thus $\al,\be,\gamma$ determine linear maps 
$B,T\in\End(U_\infty), C\in\End(W_\infty)$ and $P:U_\infty\to W_\infty$ and $Q:W_\infty\to U_\infty$ and by assumption $C$ and $T$ are semisimple.

The corresponding system of differential equations then takes the form
\beq\label{eq: system}
\partial v= (Az+B+T + Q(z-C)^{-1}P) v
\eeq
for a holomorphic function $v$ (with values in $U_{\infty}$).
This system corresponds to the meromorphic connection on the trivial holomorphic bundle on $\IP^1$ with fibre $U_\infty$ determined by the matrix of meromorphic one-forms:
\beq\label{eq: conn}
\cA= (Az+B+T + Q(z-C)^{-1}P)dz.
\eeq
In general this connection
will have a pole of order three (irregular singularity of Poincar\'e rank two) at $z=\infty$ and a simple pole (Fuchsian singularity) at each eigenvalue of $C$.

If instead we first perform a symplectic transform \eqref{eq: spt} before passing to the connection \eqref{eq: conn} then a different connection will be obtained, usually on a different rank bundle.
Indeed if the transform moves $j \in J$ to $\infty$ then the resulting bundle will have fibre $U_j= V\ominus W_j$.
On the other hand (after a generic symplectic transform) it will have fibre $V$---if no points of $J$ are at $\infty$.
Thus in general connections will be obtained on $k+1$ bundles of different ranks (and in fact on each of these bundles a whole family of connections is obtained corresponding to the affine subgroup of $\SL_2(\IC)$ stabilizing $\infty\in \bP$).

Our basic aim is to show that the isomonodromic deformations of such connections (on different rank bundles) are governed by the {\em same} system of nonlinear differential equations when we vary  $T$ and $C$.

\section{Symplectic vector spaces}\label{sn: sympl}

Suppose (as above) we have a finite dimensional 
complex vector space $V$ graded by a finite set $J$, so that
$V = \bigoplus_{j\in J} W_j$ for subspaces $W_j$.
Consider the complex vector space
$$\IM:= \End(V)^\circ =  \bigoplus_{i\ne j\in J}\Hom(W_i,W_j)$$
of linear maps in both directions between each pair of  distinct vector spaces $W_j$.

Now suppose we have an injective map $\ba:J\hookrightarrow \bP=\IC\cup\{\infty\}$ (as above) so we may identify $J$ with its image in $\bP$.
Thus in effect we have a vector space $W_p$ for each $p\in \bP$
(taking $W_p=\{0\}$ if $p\notin J$). 
We write $a_i=\ba(i)$ for the complex number corresponding to $i\in J\setminus \{\infty\}$.
Then we may define a complex symplectic structure on $\IM$ by the formula
$$\omega = \omega_\ba=
\sum_{i\ne j\in J\setminus\{\infty\}} 
\frac{\tr(dB_{ij}\wedge dB_{ji})}{2(a_i-a_j)}
+
\sum_{i\in J\setminus\{\infty\}} \tr(dB_{i\infty}\wedge dB_{\infty i})
$$
where $B_{ij}\in \Hom(W_j,W_i)$.
We will sometimes write $\IM=\IM_\ba$ when we think of it as a complex symplectic manifold with the symplectic form $\omega_\ba$.
If we define constants $\phi_{ij}\in \IC$ for $i,j\in J$ so that 
$\phi_{ii}=0$, 
$\phi_{ij} = (a_i-a_j)^{-1}$ if $i,j\neq \infty$ and 
$\phi_{i\infty} = 1 = -\phi_{\infty i}$
then
$$\omega = 
\frac{1}{2}\sum_{i,j\in J} \phi_{ij}
\tr(dB_{ij}\wedge dB_{ji}) = 
\frac{1}{2}\sum_{i,j\in J}
\tr(dX_{ij}\wedge dB_{ji}).$$
\noindent  where  we have defined $X_{ij} = \phi_{ij} B_{ij}\in \Hom(W_j,W_i)$ 
for all $i,j\in J$.

Now suppose we have matrices $\al,\be,\ga$ as above, and we normalise them as in \eqref{eq: normn}.
Then for $i\ne j\in J$ we will identify  $B_{ij}$ above with the component of $\ga$ in 
$\Hom(W_j,W_i)$, 
so that 
$P$ has components $B_{\infty i}$,  
$Q$ has components $B_{i \infty}$ and the $\Hom(W_j,W_i)$ component of $B$ is $B_{ij}$ (for $i\ne j\in J\setminus \{\infty\}$). 
In these terms
\begin{align}
\omega &= 
\tr( dQ\wedge dP ) + \frac{1}{2}\tr \left(dX\wedge dB\right)
\label{eq: omega QPXB}
\end{align}
where $X=\ad_A^{-1}(B)\in \Im(\ad_A)\subset \End(U_\infty)$.

Finally if we define $\Xi = \phi(\Ga)$ where $\phi:\End(V)\to \End(V)$ is the linear map given in components by $(B_{ij})\mapsto(\phi_{ij}B_{ij})$ (so that $X_{ij}$ is always the $\Hom(W_j,W_i)$ component of 
$\Xi$) 
or in other words
$$
\Xi = \bmx 0 & -P \\ Q & X \emx,\qquad
\Ga = \bmx 0 & P \\ Q & B \emx \in \End(W_\infty\oplus U_\infty),
$$
then the symplectic form is
\begin{align} \label{eq: XiGa sp form} 
\omega &= \frac{1}{2}\tr (d\Xi\wedge d\Gamma). 
\end{align}

Note the basic property of $\phi$, that it is skew-adjoint: 
$\Tr(\phi(E) F) = - \Tr(E\phi(F))$.
In particular this implies $\tr(\Xi\Ga)=0$.
The key fact about these symplectic forms is the following:
\begin{prop}\label{prop: invce of omega}
The symplectic form $\omega$ on $\IM$ is equivariant under the symplectic transformations \eqref{eq: spt} of $M$.
\end{prop}
Said differently, 
given an injective map $\bfa: J\to \bP$ and a symplectic matrix $g\in \SL_2(\IC)$ we obtain a map 
$\varphi:\IM_\ba\to\IM_{g\cdot \ba}$ (by acting on $M$ with $g$ and renormalising).
We are claiming that $\varphi$ is symplectic.

\pf
The group $\SL_2(\IC)$ is generated by 
i) the scalings $(\partial,z)\mapsto (\partial/c, cz)$ for $c\in \IC^*$, 
ii) the shearings $(\partial,z)\mapsto (\partial + cz, z)$ for $c\in \IC$,
and iii) the Fourier--Laplace transform 
$(\partial,z)\mapsto (-z, \partial)$ (beware this convention is sometimes taken to be the inverse Fourier--Laplace transform).
Note that i) and ii) correspond (on the level of connections \eqref{eq: conn}) respectively to scaling the coordinate $z$ and to tensoring by the connection $czdz$ on the trivial line bundle. 

For i), the resulting action on $\Ga$ is 
$$
\bmx 0 & P \\ Q & B \emx
\mapsto 
\bmx 0 & P/c \\ cQ & cB \emx
$$
and $a_j\mapsto c^2 a_j$ for all $j$. Thus $\Ga$ becomes $\eps \Ga$ where $\eps= \diag(1/c,c)$ and in turn it follows that $\Xi$ becomes $\Xi \eps^{-1}$. Thus from \eqref{eq: XiGa sp form}, this transformation leaves $\omega$ unchanged.

For ii), $\ga$ is unchanged and each $a_j$ is replaced by $a_j-c$, which again fixes $\omega$.

For iii) recall $W_0$ is the kernel of $A$ (or equivalently of $\beta$), so we may write 
\beq\label{eq: 3x3 Ga}
\be = \bmx 1 & & \\ & 0 & \\ & & -A_r \emx,\quad 
\Ga=
\bmx 
0  & P_0 & P_r \\
Q_0 & 0 & B_{0r} \\
Q_r & B_{r0} & B_{rr}
\emx \in \End(W_\infty\oplus W_0 \oplus V^r)
\eeq
where $V^r= \bigoplus_{j\ne 0,\infty} W_j$ is the ``rest of $V$''. 
After Fourier--Laplace and renormalising we then find $\be,\Ga$ become:
\beq\label{eq: after FL}
\be' = \bmx 0 & & \\ & 1 & \\ & & A^{-1}_r \emx,\quad 
\Ga'= \eps\Ga
\eeq
where $\eps = \diag(1,-1,-A_r^{-1})$.
Noting that the roles of $W_0$ and $W_\infty$ have been swapped (and $a_i\mapsto -1/a_i$), we find directly that 
\begin{lem}\label{lem: multipliers}
Under the Fourier--Laplace transform $\Xi$ becomes 
$\Xi'=\Xi\eps^{-1}$.
\end{lem}
Thus $\omega$ in  \eqref{eq: XiGa sp form} is clearly preserved.
\epf

In other words the above formulae define 
a symplectic action of $\SL_2(\IC)$ on  the symplectic vector bundle  
$\IM\times \bP^J\setminus \text{(diagonals)}$ (with fibres $\IM_\ba$) covering the standard action (diagonal M\"obius transformations) on the base 
$\bP^J\setminus \text{(diagonals)}$.

\section{Invariance of spectral invariants}\label{sn: isospectral}

As is by now well known, isomonodromy systems have degenerations (dating back at least to \cite{garn1919}) which are isospectral integrable systems (solvable by spectral curve methods i.e. in terms of abelian functions). Here we will briefly touch on this in our context, focusing on the results we will need elsewhere in this article (our main interest being the isomonodromic system, requiring in general more complicated functions to solve them).
This extends some work of Adams--Harnad--Hurtubise \cite{AHH-dual}.

Suppose we fix $\wh T$ and $A$ as above (except in this section we do  not need to assume $\wh T$ is semisimple).
Let $\cI = \cI(\wh T, A)$ denote the ring of functions on $\IM$ generated by the 
coefficients of the polynomial 
$$\cP(\la,z) = 
\det(z-C)\det\left(\la - ( Az+B+T + Q(z-C)^{-1}P)\right)\in \IC[\la,z].$$

This is the ring of  spectral invariants of the matrix \eqref{eq: conn} of meromorphic
one-forms.

Now if we perform a symplectic transformation \eqref{eq: spt} changing the matrix $M$ then the resulting matrix \eqref{eq: conn} (and $\wh T, A$) will in general change.
However it turns out that the ring
$\cI$ is unchanged (as a ring of functions on $\IM$).

\begin{lem}
The equality
$$\cP(\la,z) = \det(\al \la + \be z - \ga)$$
holds, and
so the symplectic transformation \eqref{eq: spt} (and more generally the corresponding $\GL_2(\IC)$ action) just changes the coordinates $\la,z$ on $\IC^2$, and thus preserves $\cI$.
\end{lem}
\pf
Let $N=\al \la + \be z - \ga\in \End(V)\otimes \IC[\la,z]$ as in 
\eqref{eq: normn},\eqref{eq: Ga decomp}. 
Now consider the equation 
$$
N =
\left(\begin{matrix}
1 & 0 \\ 
s &  1
\end{matrix}\right) 
\left(\begin{matrix}
x & t \\ 
0 & y
\end{matrix}\right)
\in 
\End(W_\infty\oplus U_\infty)\otimes \IC(\la,z). 
$$
This has unique solution: 
$$t=-P,\  x= z-C,\  s=Q(C-z)^{-1},\  y= \la-(Az+B+T + Q(z-C)^{-1}P)$$
so we see that $\det(N) = \det(x)\det(y) = \cP(\la,z)$ as required.
Thus under the transformation (and renormalising $N$) 
$\cP$ just undergoes a coordinate change (and multiplication by an invertible constant due to the renormalisation) so $\cI$ is unchanged. 
\epf

In particular the corresponding spectral curves (cut out in $\IC^2$ by $\cP(\la,z)$) are isomorphic (and the corresponding curves obtained by first clearing the factors of $\det(z-C)$ are birational).

On the other hand we know for general reasons that $\cI$ is a {\em Poisson commutative} subring of the ring of functions on the symplectic manifold $\IM$ (with the symplectic form $\omega$ determined by $A$).
This follows for example from the Adler--Kostant--Symes theorem, once we relate $\IM$  via a moment map to an appropriate loop algebra.
See appendix \ref{apx: AKS}.
Moreover this approach immediately gives formulae for the (isospectral, time independent) Hamiltonian vector fields determined by functions $H\in \cI$.

The case $A=0$ (so $B=0$ too) was studied in a sequence of papers by Adams--Harnad--Hurtubise--Previato (starting with \cite{AHP})
extending and answering questions in the work of Moser \cite{moser79}.
In particular the isomorphisms above generalize the duality of 
\cite{AHH-dual}
which  was used to explain a number of examples of integrable systems admitting different Lax representations.
Here we get more alternative Lax representations, even in the cases considered in \cite{AHH-dual} (since we may always, using a more general symplectic transform, pass to a system having just a pole of order three, on the trivial bundle with fibre $V$).
Indeed we now see the symmetric diagram in \cite{AHH-dual} p.303, \cite{previato-70} p.438 (which looks a little like two pages of an open book, if we draw the spine vertically down the middle) has a natural generalization where the book has $k+1$-pages, each page of which looks as follows (with the spine on the left):

\begin{figure}[ht]
	\centering
	\input{leaf2.pstex_t}
\end{figure}

Here $\Rat(U_\infty,C)$ is the Lie algebra of rational maps $\IP^1\to \End(U_\infty)$ having poles only at $\infty$ and at the eigenvalues of $C=T_\infty$, the top horizontal arrow is the map 
$\Ga\mapsto \cB=Az+B+T+Q(z-C)^{-1}P$ (and is shown to be a  moment map 
in Lemma \ref{lem: mmap}),
$G_\infty$ is the subgroup of $\GL(U_\infty)$ centralizing $A$ and $T$ (the subgroup of the group of global gauge transformations of the bundle $\IP^1\times U_\infty\to U_\infty$ preserving the shape of $\cB$), and $H_\infty$ is the subgroup  of $\GL(W_\infty)$ centralizing $C$. The other `pages' are obtained by performing a transformation to move the points $J\subset \bP$, changing which point is at $\infty$ and thus the spaces $U_\infty, W_\infty$.
Note that the product $G_\infty\times H_\infty$ is the centralizer in $\prod\GL(W_j)$ of $\wh T$, so does not change from page to page. 
As in the case of \cite{AHH-dual} the bottom horizontal map is an injective Poisson map and (at least when $\wh T$ is semisimple) all the quotients here are manifolds when restricted to dense open subsets. 

In other words our  space $\IM$ thus generalizes the ``generalized Moser space''  $M$ of Adams et al.

\section{Time-dependent Hamiltonians}

We will first review precisely what we mean by a time-dependent Hamiltonian system (with multiple times), since some strictly weaker  notions have been used in the context of isomonodromy recently.
This formulation (which we view as standard---as in \cite{JMMS}) 
enables us to see the link between 
the isomonodromy $\tau$-function and the time-dependent Hamiltonians.

Let $\IM$ be a complex symplectic manifold with symplectic form $\omega$, and let $\IB$ be another manifold (which will play the role of the space of times).
Consider the trivial symplectic fibre bundle
$$\pi:\IF:=\IM\times\IB\to \IB$$
with base $\IB$ and fibre $\IM$.
In particular the tangent space to $\IF$ at any point has a splitting into horizontal and vertical subspaces (beware that $\IF$ stands for `fibration' and not for `fibre' here).

Then a vector field 
$X$ on $\IB$ and a function $H$ on $\IF$ determine a vector field $\wt X=X+v_h$ on $\IF$, by taking the horizontal component to equal $X$ and the vertical component to equal  the Hamiltonian vector field 
$v_{H}$ of $H$, defined by $dH = \omega(\cdot, v_H)$.
In the case when $\IB$ has dimension one, this is just the usual notion of a time-dependent Hamiltonian, and the flows correspond to the differential equation $dm/dt = -v_H(m)$ for $m\in \IM$.
In other words if $\{m_i\}$ are local coordinates on $\IM$ then 
the coordinates of a solution $m(t)$ evolve according to
$$\frac{d m_i}{dt} = \{ m_i, H \}$$
where the Poisson bracket of $\omega$ is defined as usual via
$\{f,g\} = \omega(v_f,v_g)$, so that $v_H = \{H,\cdot\}$.
Thus in general we would like a Hamiltonian function on $\IF$ for each vector field on the base. 

\begin{defn}
A system of time-dependent Hamiltonians on $\pi:\IF=\IM\times\IB\to \IB$ is a global section $\varpi$ of the vector bundle $\pi^*(T^*\IB)$ over $\IF$, i.e. it is a one-form on $\IF$ whose vertical component is zero.
\end{defn}
Thus if we have a vector field $X$ on $\IB$ then we obtain a time-dependent Hamiltonian function $H=\langle \varpi, X \rangle$ on $\IF$.  
If we have some coordinates $\{t_i\}$ on $\IB$ we may trivialise the cotangent bundle and write 
$$\varpi = \sum_i H_i dt_i$$
and then $H_i$ is the Hamiltonian function for the vector field 
$\partial/\partial t_i$ on $\IB$.
The corresponding differential equations may be written as 
$d m_i = \{m_i , \varpi\}$, 
where $d$ is the exterior derivative on $\IB$.

Said differently $\varpi$ tells us how to modify the trivial connection on $\IF$ to obtain a new connection. 
This may be rephrased in terms of symplectic connections as follows.
First let $\wh \omega = \pi_1^*\omega$ be the two-form on $\IF$ obtained by pulling back the form $\omega$ on $\IM$ along the projection $\pi_1:\IF\to \IM$
(coming from the fact that $\IF$ is trivial as a fibration over $\IB$).
In general any two-form $\Omega$ on $\IF$ which restricts to $\omega$ on each fibre determines a connection on the bundle $\IF\to \IB$: the field of horizontal subspaces is given by the orthocomplement of the vertical subspace:
\begin{align*}
H_p &= \{ u\in T_p\IF \st \Omega(u,v) = 0 \text{ for all $v\in T_p\IM$ }\}\\
&= %
\ker\Bigl(T_p\IF\to T^*_p\IM; u\mapsto [v\mapsto \Omega(u,v)]\Bigr)
\end{align*} 
for all $p\in \IF$, which is easily seen to be complementary to the vertical subspace.
Taking $\Omega = \wh \omega$ defines the trivial connection on $\IF$.
The connection given by some time-dependent Hamiltonians $\varpi$
is given by the two-form 
$$\wh \omega - d\varpi$$ on $\IF$.
Note that by assumption 
this will again restrict to $\omega$ on each fibre. 
Thus from this point of view the time-dependent Hamiltonians give the difference between the original trivial connection and the new (interesting) one.
It is a general fact about symplectic connections (see \cite{GLSW} Theorem 4) that if $\Omega$ is closed then the local isomorphisms between open subsets of the fibres obtained by integrating the connection, will be symplectic (and in our situation $\wh \omega-d\varpi$ is clearly closed). \footnote{It is straightforward to check the two procedures to obtain a connection from $\varpi$ agree: e.g. in coordinates if $H=\langle\varpi,\partial/\partial t\rangle$ and $v_H$ is the corresponding (vertical) Hamiltonian vector field, so that $dH = \omega(\cdot,v_H)$ on each fibre. We should check that the vector field $u=\partial/\partial t + v_H$ on $\IF$ (defined by the first recipe) does indeed satisfy $\Omega(u,v)=0$ for all vertical $v$, where $\Omega = \wh\omega - dH\wedge dt$. 
But this is immediate as $\wh\omega(u,v) = \omega(v_H,v)$.}

We are mainly interested in the case where the resulting symplectic connection is integrable, i.e.
that the vector fields $X_i := \partial/\partial t_i + \{H_i,\cdot\}$ on $\IF$ commute (for local coordinates $t_i$ on $\IB$).
One may readily verify that 
\begin{lem}(cf. \cite{GLS96} (1.12))\label{lem: ham curvature}
The Lie bracket $[X_i,X_j]$ is the (vertical) Hamiltonian vector field
associated to the function
$$f_{ij}=\frac{\partial H_j}{\partial t_i} 
-\frac{\partial H_i}{\partial t_j} + \{H_i,H_j\}=-\Omega(X_i,X_j)$$
where $\Omega=\wh \omega -d\varpi$.
\end{lem}
Thus if the connection is integrable, each function $f_{ij}$ is constant on each fibre $\IM$ of $\IF$, so is the pullback of a function on $\IB$.
One may also demand the stronger condition (``strong integrability''), 
that $f_{ij}=\{H_i,H_j\}=0$.
In that case, since $\wh \omega(X_i,X_j) = \{H_i,H_j\}=0$, 
the restriction of $d\varpi$ to each solution leaf is zero, i.e.
$\varpi$ restricts to a closed one-form on each solution leaf.
Then, pulling back to the base, we may may regard $\varpi$ as a flat connection on the trivial line bundle on $\IB$, and locally define a holomorphic function $\tau$ on $\IB$ as its horizontal section (well defined upto a scalar multiple), so that
\beq \label{eq: tau function defn}
d\log \tau = \varpi.
\eeq
This will be the case in our situation, so we will get new $\tau$ functions, analogous to those of  \cite{JMMS, JMU81}.

For our purposes (with $\IM=\End(V)^\circ$ and $\omega = \tr(d\Xi \wedge d \Ga)/2$) the following description of the Hamiltonian equations will be useful:

\begin{lem}  \label{lem: ham to eq}
Suppose the Hamiltonian one-form $\varpi$ on $\IF=\IM\times \IB\to \IB$ 
satisfies $$d_\IM\varpi  = \tr(\cE\wedge d_\IM\Xi)$$
for some $\IM$ valued one-form $\cE$ on $\IF$ with vertical component zero.
Then the corresponding nonlinear differential equations (for local sections $\Ga$ of $\IF$) are given by
$$d_\IB \Ga = \cE.$$
\end{lem}
\pf Suppose there is just one time $t$ and $\varpi = H dt$, so 
\begin{align*}
\tr(\cE\wedge d_\IM\Xi ) &= d_\IM H\wedge dt = \omega(\cdot, v_H)\wedge dt \\
&= -\frac{1}{2}\iota_{v_H}\tr(d_\IM \Xi \wedge d_\IM \Ga)\wedge dt\\
&= \cdots = \Tr(d_\IM \Xi\langle d_\IM \Ga,v_H\rangle)\wedge dt
\end{align*}
where in the last line we use the skew-symmetry of $\phi$.
Thus since the trace pairing is nondegenerate we deduce
$\langle d_\IM \Ga,v_H\rangle dt = -\cE$, which yields the differential equation $d\Ga/dt = \langle \cE,d/dt\rangle$, i.e. $d_\IB\Ga = \cE$.
(For multiple times just repeat the above argument.) 
\epf

Now we will describe the time-dependent Hamiltonian system
that is the central focus of this article.

\subsection{The space of times}
Suppose we are in the situation of Section \ref{sn: sympl}
with a  vector space 
$V$ and normalised matrices $\al,\be,\ga\in\End(V)$  determining $J\subset \bP=\IC\union\{\infty\}$ and a $J$-grading $\bigoplus_{j\in J} W_j$ of $V$ etc.
Define $$\wh T=\delta(\gamma)\in\End(V),$$
to be the block-diagonal part of $\ga$ and let $T_j\in \End(W_j)$
be the component of $\wh T$ in  $\End(W_j)$.
(Thus $C=T_\infty$ and $T$ is the component of $\wh T$ in $\End(U_\infty)$.) 

By hypothesis $\wh T$  is semisimple and so
also each $T_j$ is semisimple.
(We do not make any assumption of distinct eigenvalues---any multiplicities are permitted.)
Thus each vector space $W_j$ has a finer decomposition into the eigenspaces of $T_j$. 
Let $I_j$ denote the set of eigenspaces of $T_j$
and let $V_i\subset W_j$ be the corresponding eigenspace, for 
$i\in I_j$. Thus 
\beq \label{eq: fine decomp}
W_j = \bigoplus_{i\in I_j}V_i\qquad\text{so that}\qquad
V=\bigoplus_{i\in I}V_i
\eeq
where $I=\bigsqcup_{j\in J}I_j$.
The space of times $\IB$ is the space of $\wh T$ such that the decomposition \eqref{eq: fine decomp} does not change. 
Explicitly we may write
$$\wh T = \sum_{i\in I}  t_i\Id_i$$
for some complex numbers $t_i$, where $\Id_i\in \End(V)$ is the idempotent for $V_i$. Thus $\wh T$ is identified with a point $\{t_i\}$ of $\IC^I$ and the space of times with an open subset of $\IC^I$: 
$$\IB = \left\{ \wh T = \sum_{i\in I}  t_i\Id_i\ \bigl\vert\ 
t_i\in \IC \text{ and if $i,i'\in I_j$ for some $j\in J$
then $t_{i}\ne  t_{i'}$} \right\}.$$
Thus %
 $\IB \cong \Pi_{j\in J} \left(\IC^{\abs{I_j}}\setminus\{\text{diagonals}\}\right).$ 
In particular the fundamental group of $\IB$ is a product of Artin braid groups.

\subsection{Hamiltonians}

The Hamiltonian one-form on $\IF=\IM\times\IB\to \IB$ we are interested in is given by the expression:
\beq \label{eq: hamns}
\varpi = \varpi_0 + \varpi_1
\eeq
where

\begin{center}
\fbox{
$\begin{array}{l}
  \varpi_0=\frac{1}{2}\tr\left(\wt{\Xi\Ga}\delta(\Xi\Ga)\right)
              -\tr\left(\Xi\gamma\Xi d\wh T\right),\\
\\
  \varpi_1 = \tr(X^2TdT) + \tr(PAQT_\infty dT_\infty),
 \end{array}$
}
\end{center}
and $\wt{\Xi\Ga} = \ad_{\wh T}^{-1}[d\wh T,\Xi\Ga]$.

This determines a symplectic connection $\wh\omega - d\varpi$ and thus 
some nonlinear differential equations for $\Ga$ 
as a function of the times $\wh T$.
These equations will be written explicitly in Section \ref{sn: differential equations} below.
Our next aim is to establish the following,  which is one of the main results of this article.

\begin{thm}\label{thm: Ham invariance} 
The one-form $\varpi$ is invariant under all the symplectic transformations \eqref{eq: spt} up to some simple global gauge transforms. 
These gauge transformations are tangent to the symmetries of $\IM$ and so the reduced equations are completely invariant.
\end{thm}

Before proving this we will first discuss such gauge transformations.
Suppose $g: \IB \to \GL(V)$ is a holomorphic map such that for each $\wh T\in \IB$ the map $\Ga\mapsto g\Ga g^{-1}$ is a well-defined symplectic automorphism of $\IM$. 
Then there is a new connection on $\IF$ whose horizontal sections are 
$\Ga=g\Ga_0 g^{-1}$
with $\Ga_0$ constant (and $g$ varying over $\IB$).
The gauge transformations we will need are of this form
where
\beq\label{eq: gen gauget}
g(\wh T) = \exp(\la \wh T^2/2) \in \GL(V)
\eeq
for some constant $\la\in \End(V)$ of the form 
$\la = \sum_{j\in J}\la_j \Id_j$ where $\la_j\in \IC$ and $\Id_j$ is the idempotent for $W_j\subset V$.
The horizontal sections of the new connection then satisfy
$d_\IB\Ga = [\th,\Ga]$ where $\th = g^{-1}dg = \la \wh T d\wh T$.
In turn if follows from 
Lemma \ref{lem: ham to eq} that a
 two-form on $\IF$ %
for the new connection is
\begin{align*}
\wh \omega ' 
&=  \wh \omega + \Tr (d\Xi [\th, \Ga]) \\
&=  \wh \omega + d\Tr (\Ga\Xi \th) 
\end{align*}
so the gauge transformation corresponds to a Hamiltonian  term 
$\Tr (\Ga\Xi \th)$. %
Two examples will cover the cases we need.

\begin{eg}\label{eg: gauget1}
 For example if $\la_j = 0$ for all $j\in J\setminus\{\infty\}$ 
then  the Hamiltonian term is 
$\Tr (\Ga\Xi \th) =  \la_\infty\Tr(PQ T_\infty d T_\infty).$
\end{eg}

\begin{eg}\label{eg: gauget2}
As a second example suppose 
$\la_j = a_j^{-1}$ if $j\in J\setminus \{0,\infty\}$ 
and is zero otherwise.
Then using the notation of \eqref{eq: 3x3 Ga} 
the Hamiltonian term is 
$\Tr (\Ga\Xi \th) =
 \tr \bigl((B_{r0}X_{0r}+B_{rr}X_{rr}-Q_rP_r)\th_r\bigr)$ 
where 
$\th_r= \sum_{j\in J\setminus\{0,\infty\}}  T_j d T_j/a_j = 
A_r^{-1} T_rdT_r$.
\end{eg}

\pfms (of Theorem \ref{thm: Ham invariance}).\  
From the discussion in Section \ref{sn: sympl} it is immediate that $\varpi_0$ is invariant under all the symplectic transformations 
(for example under the Fourier--Laplace transform $\ga, \Gamma, \Xi, \wh T$ become $\eps\ga, \eps \Gamma, \Xi\eps^{-1}, \eps\wh T$ respectively and the factors of $\eps$ cancel each other,
noting also that $\delta(\wt {\Xi\Ga})$ does not change when $\wh T$ is replaced by $\eps\wh T$).
We should check that the change in $\varpi_1$ may be compensated for by gauge transformations of the symplectic fibration $\IF\to \IB$. 
As before, we go through the subgroups generating the symplectic group.

i) for the scalings, each of $Q,B$ and $T$ is scaled by $c\in \IC^*$, $A $ by $c^2$ and $P, T_\infty, X$ by $c^{-1}$ and so both terms of $\varpi_1$
are invariant. 

ii) For the shearings, $\ga$ is unchanged and $A$ becomes $A-c\Id$ for a constant $c\in \IC$. Thus  $X=\ad_A^{-1}B$, and therefore the term $\tr(X^2TdT)$, 
is unchanged.  However from the final term we pick up an extra term 
$-c\tr(PQT_\infty dT_\infty)$.
This may be removed by the gauge transformations of Example \ref{eg: gauget1} with $\la_\infty = c$.
Note that in fact Example \ref{eg: gauget1} computes 
the gauge transform of the trivial connection $\wh \omega$ whereas we are interested in the transform of $\wh \omega - d\varpi$.
However the difference is the same in both cases since 
$\varpi$ is invariant under the action of $g$ (and indeed this holds for all transformations of the form \eqref{eq: gen gauget}). 

iii) For the Fourier--Laplace transform, using the formulae \eqref{eq: after FL} we find $\varpi_1$ 
minus its transformed version is 
$$\tr\Bigl( (X_{r0}X_{0r} +X_{rr}^2 - Q_rP_rA_r^{-1}  - X_{rr}A_rX_{rr}A_r^{-1})T_rdT_r\Bigr)$$
which simplifies further to 
equal the term $\Tr(\Ga\Xi\th)$ in Example \ref{eg: gauget2}, appearing from the gauge transformation by $g=\exp(A_r^{-1}T_r^2/2)$.
\epfms

\subsection{Further properties of the Hamiltonians}\label{ssn: further props}

\begin{thm}
Let $\varpi = \sum _{i\in I} H_i dt_i$ be the Hamiltonian one-form \eqref{eq: hamns} so $H_i$ is a function on the total space $\IF$. Then for all $i,j\in I$:

\noindent i) $\{H_i, H_j\} = 0$ as functions on any fibre $\IM$ of $\IF$,

\noindent ii) ${\partial H_i}/{\partial t_j} = {\partial H_j}/{\partial t_i}$ 

\noindent
and consequently the nonlinear connection on $\IF$ determined by $\varpi$ is integrable and the restriction of $\varpi$ to each solution leaf is closed (so new $\tau$ functions may be defined by \eqref{eq: tau function defn})
\end{thm}

\pf
For the first property we will show that each $H_i$ is in the subring 
$\cI$ of the functions on $\IM$ defined in section 
\ref{sn: isospectral}. The result then follows as $\cI$ is Poisson commutative.
If $i\in I_\infty$ we claim that
$$H_i = \frac{1}{2}\res_i\tr (\cA\cB)$$
where $\cB$ is the matrix valued rational function  $Az + B+T +Q(z-C)^{-1}P$ on $\IP^1$, $\cA= \cB dz$ 
and $\res_i$ denotes the residue at $t_i$ (an eigenvalue of $C=T_\infty$).
This shows that $H_i$ is the restriction of an invariant 
function on the loop algebra, 
and so $H_i$ is in $\cI$ by the Adler--Kostant--Symes theorem (see Appendix 
\ref{apx: AKS}).
The claim itself will be established in proof of Theorem \ref{thm: relate to JMU} below.
If $i\in I_j$ and $j\ne \infty$ then we can do a symplectic transform to move $j\in \bP$ to $\infty$ and repeat as above. However the Hamiltonians are not quite invariant under the symplectic transformations: we need to check that when we twist by the gauge transformations the result stays in $\cI$. It is sufficient to check that all the terms that modify the Hamiltonians, induced by gauge transformations, are themselves in $\cI$. These terms are of the form $\Tr(\Ga\Xi\th)$ (where $\th = \la\wh Td\wh T$). Expanding this out we see the Hamiltonian terms  are of the form 
$$-\sum_{i\in I} \la_i\tr(Q_iP_i)t_idt_i$$
for constants $\la_i$. Here $Q_i= \Ga\circ \iota_i:V_i\to V$ and 
$P_i = -\pi_i\circ \Xi: V\to V_i$ where $\iota_i$ and $\pi_i$ are the inclusion and projection for $V_i\subset V=\bigoplus V_i$. 
Thus (since $t_i$ is held constant in i)) it is sufficient to verify the following.

\begin{lem}
For any $i\in I$ the function $\tr(Q_iP_i)$ on $\IM$ is in $\cI$.
\end{lem}
\pf
For $i\in I_\infty$ this is clear since $\tr Q_iP_i$ is the residue of $\tr\cA$ at $z=t_i$ and this is the restriction of an $\Ad$-invariant function on the loop algebra, so the result follows from the Adler--Kostant--Symes theorem. For $i\in I_j, j\ne \infty$ 
we may do a symplectic transform to move $j\in \bP$ to $\infty$: we may write $\tr(Q_iP_i) = -\tr (\Ga \Id_i \Xi) =-\tr (\Id_i \Xi\Ga)$ 
where $\Id_i$ is the idempotent for $V_i\subset V$.
From the results of section \ref{sn: sympl}, this expression $\tr (\Id_i \Xi\Ga)$ is invariant under the symplectic transformations, so again $\tr(Q_iP_i)$ is the trace of a residue and thus in $\cI$. 
\epf

Thus $H_i\in \cI$ for all $i$ and so they all Poisson commute.
The second property is a straightforward verification.
The integrability now follows from Lemma 
\ref{lem: ham curvature}.
\epf

Note that we will see below in Section \ref{sn: reductions} that $-P_iQ_i$ is a moment map for the natural action of $\GL(V_i)$ on $\IM$, and so each of the functions $\tr(Q_iP_i)$ is constant on the reduced manifolds (the symplectic quotients of $\IM$)---i.e. the gauge transformations only change the reduced Hamiltonians by constants and so will not change the reduced differential equations.

\begin{thm}\label{thm: relate to JMU}
The Hamiltonian one-form $\varpi$ in \eqref{eq: hamns}
equals $\varpi_\infty+ \sum_{i\in I_\infty} \varpi_i$
where
$$\varpi_i =  \frac{1}{2} \res_{z=t_i}(\tr \cA\cB) dt_i$$
and
$$\varpi_\infty = 
\res_\infty \tr \left(
(d_z\wh g) \wh g^{-1} z dT
\right)$$
for any formal series $\wh g=1+g_1/z+\cdots$ putting $\cA$ into formal normal form at $z=\infty$.
\end{thm}

Here $\cA=\cB dz$, and $\wh g$ is discussed fully in Appendix \ref{apx: leading term}.

\pf
Computing the residues at finite distance yields 
\beq\label{eq: infty terms}
\sum_{i\in I_\infty} \varpi_i =\tr(PAQCdC) + \Tr (P(B+T)QdC) +\tr(PQ\wt{PQ})/2
\eeq
where $C=T_\infty$.
At $z=\infty$, by writing 
$\wh g = 1 +g_1/ z +g_2/z^2 +\cdots$ we find 
$$\varpi_\infty 
=-\Res_\infty\tr (g_1 dT dz/z) = \tr(g_1 dT).$$
This is more difficult to compute explicitly. The result is
$$
\tr(QT_\infty PdT) 
-\tr(XTXdT) 
+ \tr(X^2TdT) 
$$
\beq\label{eq: varpi-infty}
+\tr ([X,QP]dT)- \tr(XBXdT)
\eeq
$$
+ \tr(\wt{QP}\delta(QP))/2
+ \tr(\wt{QP}\delta(XB))
+ \tr(\wt{XB}\delta(XB))/2.
$$
See appendix \ref{apx: leading term} for the details (where we also give the expression that arises if we were working with other complex reductive groups).
Finally we have the pleasant task of showing that the sum of 
 \eqref{eq: infty terms} and \eqref{eq: varpi-infty}
equals the shorter expression \eqref{eq: hamns}
in terms of $\Ga,\Xi\in \End(V)$.
Indeed the quartic terms equal
$\frac{1}{2}\tr(\wt{\Xi\Ga}\delta(\Xi\Ga))$,
the cubic terms equal
$-\tr(\Xi\Ga\Xi d\wh T)$,
and the quadratic terms equal
$\varpi_1 -\tr(\Xi\wh T\Xi d\wh T)$.
\epf

If all the simple poles are nonresonant ($Q_iP_i$ has no eigenvalues differing by a positive integer for $i\in I_\infty$) then  
there is a unique formal isomorphism $\wh g_i\in G\flb z-t_i\frb$ with  constant term $1$, such that  
$\wh g_i[\cA] = Q_iP_idz/(z-t_i)=d_z\xi_i$, 
where $\xi_i = Q_iP_i\log(z-t_i)$.
Then since $d_\IB\xi_i$ is minus the  principal part of $\cB dt_i$ at $z=t_i$  
it follows that 
$$\varpi_i =
\res_{z=t_i} \tr \left(
(d_z\wh g_i) \wh g_i^{-1} d_\IB\xi_i
\right)
$$
for all $i\in I_\infty$ (and similarly at $z=\infty$), and thus that $\varpi$ is an extension to our context of the one-form of Jimbo--Miwa--Ueno \cite{JMU81} p.311, that they used to define $\tau$ functions.\footnote{
In \cite{JMU81} the connections are assumed to have regular semisimple leading coefficients at each irregular singularity: this immediately implies that the connections only have one level, and that the blocks $\La_i$ of the exponent of formal monodromy are nonresonant; we are considering a multilevel case and allow $\La_i$ to be resonant.}

\section{Isomonodromy}

So far we have written down a family of nonlinear connections and shown they are integrable, and that they are invariant under the full $\SL_2(\IC)$ group of symplectic transformations, up to some simple gauge transformations.
In this section we will show that local solutions to these nonlinear equations yield isomonodromic families of linear connections on the Riemann sphere.
Performing symplectic transforms then shows that the same equations control many different isomonodromic deformations (e.g. on different rank bundles with different numbers of poles).
Here we use the De\! Rham approach to isomonodromy, in terms of integrable absolute connections (cf. \cite{smid} Section 7).

\subsection{Full connections}

Suppose that $\Ga\in \IM$ is a holomorphic function of $\wh T\in \IB'$ defined on some open subset $\IB'\subset\IB$.
Write $V = W_\infty\oplus U_\infty$ as usual.
Consider the following $\End(U_\infty)$-valued 
meromorphic one-form on $\IP^1\times \IB'$:
\beq\label{eq: full conn}
\Omega = 
(Az+B)dz + d(zT)  + Qd\log(z-T_\infty)P +  
[dT,X] +\wt{\delta(XB)} + \wt{\delta (QP)}
\eeq
where $d$ is the exterior derivative on $\IP^1\times \IB'$,
$T\in \End(U_\infty)$ denotes the restriction of $\wh T$,
$\delta:\End(U_\infty)\to\End(U_\infty)$ denotes the restriction of the map $\delta$, where
$A,B,P,Q$ are determined from $\al,\be,\gamma$ by normalizing and writing
$$
\al \partial + \be z - \ga = 
\left(\begin{matrix}
z-T_\infty & -P \\ 
-Q & \partial - Az -B - T
\end{matrix}\right)\in \End(W_\infty\oplus U_\infty)\otimes \cA_1
$$
as before, with $T_\infty = C$.
Here $X = \ad_A^{-1}(B)\in \End(U_\infty)^\circ$ and, for any $R\in \End(U_\infty)$
$$\wt R := \ad^{-1}_T([dT,R]).$$
We will view $\Omega$ as a linear connection on the trivial vector bundle with fibre $U_\infty$ over the product $\IP^1\times\IB'$ (so local horizontal sections are maps $v:\cU\to U_\infty$ satisfying $dv = \Omega v$, with $\cU\subset\IP^1\times\IB'$). 
Note that the vertical  component of $\Omega$ is
\beq\label{eq: vconn}
\cA:=\langle \Omega,\partial\rangle dz = 
(Az+B+T)dz  + Q(z-T_\infty)^{-1}Pdz 
\eeq
$$
= (Az+B+T)dz  + \sum_{i\in I_\infty} \frac{Q_iP_i}{z-t_i}dz
$$
as in \eqref{eq: conn} where $U_\infty \underset{P_i}{\overset{Q_i}{\leftrightarrows}} V_i\subset W_\infty$ are the components of $P,Q$.

The main result we will establish in this section is the following.

\begin{thm}
If the local section $\Ga$ of $\IF$ is horizontal for the connection 
$\wh \omega - d\varpi$ determined by the Hamiltonians \eqref{eq: hamns}
then $\Omega$ is flat. 
\end{thm}

In principle this is possible by direct algebraic computation, which we will leave to the reader. Instead we will give a more conceptual approach. %

\pf
We need to see that $d\Omega = \Omega^2.$ Write $\Omega = \cB dz + \sum_{i\in I} \cB_i dt_i$ so $\cA = \cB dz$.
We should show that
$$\frac{\partial \cB}{\partial t_i} 
-\frac{\partial \cB_i}{\partial z} + [\cB,\cB_i] = 0$$
and
 $$\frac{\partial \cB_j}{\partial t_i} 
-\frac{\partial \cB_i}{\partial t_j} + [\cB_j,\cB_i] = 0$$
for all $i,j\in I$. The first set of equations give the `isomonodromic' evolution of $\cB$ with respect to  the times $t_i$: $\frac{\partial \cB}{\partial t_i} = \frac{\partial \cB_i}{\partial z} + [\cB_i,\cB]$.
We should check that the right-hand side coincides with the evolution determined by the nonlinear connection on $\IF$.

If $i\in I_\infty$ then as we saw in \S \ref{ssn: further props}  the corresponding Hamiltonian is 
$H_i = \res_i\tr(\cA\cB)/2$. The derivative of this is 
$\res_i\tr(\cB d\cA)$ so by Lemma \ref{lem: hamvfs} the corresponding Hamiltonian vector field at $\cB$ is $[\cB,\cB_-]$ where $\cB_-$ is the singular part of $\cB$ at $z=t_i$, i.e. $\cB_- = Q_iP_i/(z-t_i)$. 
Now we observe, from the expression for $\Omega$, 
that $\cB_i = -Q_iP_i/(z-t_i) = -\cB_-$. 
Thus the tangent to $\cB$ corresponding (under the map $\ga\mapsto \cB$) to the vector  
$\partial/\partial t_i + \{H_i,\cdot\}$ on $\IF$  is
$$\frac{Q_iP_i}{(z-t_i)^2} + [\cB,\cB_-]$$
where the first term is obtained by differentiating $\cB$ with respect to its explicit $t_i$-dependence. In turn this is 
$$\frac{\partial \cB_i}{\partial z} + [\cB_i,\cB]$$
as required.

If $i\in I_j, j\ne \infty$ then from Theorem \ref{thm: relate to JMU}  $H_i$ is the coefficient of $dt_i$ in 
$$\varpi_\infty=\Res_\infty \tr \left((d_z\wh g)\wh g^{-1} z dT \right)$$
where $\wh g$ is a family of formal isomorphisms to normal forms, as in Appendix \ref{apx: leading term}.
It is straightforward to check that $\wh g$ may be chosen so each coefficient depends holomorphically on the parameters (in the nonresonant case this is clear since $\wh g$ is uniquely determined).
This may be rewritten as
$\Res_\infty \tr \left(\wh g^{-1} (d_z\wh g) \cR \right)$
where $\cR = z \wh g^{-1} dT \wh g$.
Now since $\wh g[\cA] = \cA^0:= (Az+T +\wh\La/z)dz$ we have
$\wh g^{-1} (d_z\wh g)  = \wh g^{-1}\cA^0 \wh g - \cA$
so that $H_i$ is the coefficient of $dt_i$ in 
$$-\Res_\infty \tr \left(\cA \cR \right)$$
since $\Res_\infty \tr \left(\cA^0 zdT\right)=0$.
Next we claim that the derivative of $H_i$ is the coefficient of $dt_i$ in $-\Res_\infty \tr \left((d\cA) \cR \right)$.
To see this note that
\begin{align*}
   \Res_\infty \tr \left(\cA d\cR \right) 
&= \Res_\infty \tr \left(\cA [\cR,\wh g^{-1}d\wh g] \right) \\
&= -\Res_\infty \tr \left(\wh g^{-1}d_z\wh g [\cR,\wh g^{-1}d\wh g] \right)\quad\text{as $\cA^0$ and $dT$ commute} \\
&= \Res_\infty \tr \left([d_z\wh g \wh g^{-1}, d\wh g\wh g^{-1}]zdT\right)\\
&=0\qquad\text{by looking at the possible degrees in $z$.}
\end{align*}
Thus by Lemma \ref{lem: hamvfs} the  Hamiltonian vector field of $H_i$ at $\cB$ is 
$[\cR^i_-,\cB]$ where  $\cR^i$ is the coefficient of $dt_i$ in $\cR$
and $\cR^i_-$ is the polynomial 
part of $\cR^i \in \End(U_\infty)(\!(z^{-1})\!)$.
Thus to conclude we should check i) that $\cR_-^i=\cB_i$ (the $dt_i$ component of $\Omega$) and ii) that $\partial \cB_i/\partial z$
is the derivative of $\cB$
with respect its explicit $t_i$ dependence.

Part ii) is immediate (it only involves the term $d(zT)$). 
For i) write $\wh g = 1 + g_1/z +\cdots$ so that 
$$\cR_- = (z\wh g^{-1} dT \wh g)_- = zdT + [dT,g_1].$$
Then from Lemma \ref{lem: g1 calcn} 
$$[dT,g_1] = [dT,X] + \wt{\de (QP)} +\wt{\de(XB)}$$
so i) follows (comparing with the expression for $\Omega$).
 
Finally it is straightforward check the $t_i$-$t_j$ components commute. See e.g. \cite{Dickey} Proposition 11.2.12 p.198 for a streamlined direct algebraic approach, which extends immediately to our context. 
\epf

\begin{rmk}
Note that it is possible (similarly to \cite{JMU81})
to view the eigenvalues $a_j\in \IC$ of $A$ as times as well (higher irregular times), although  we will not do this here, since it is the lowest irregular times $T$ which are related under the $\SL_2(\IC)$ action to the pole positions $t_i$. Also 
 the Hamiltonian story is then more complicated (one runs into the symplectic trivialisation problem of \cite{smid} Remark 7.1; the symplectic structure on $\cM^*$ depends on these times so one needs to choose an a priori symplectic trivialization of $\IF$).  
In general  this problem does not arise if attention is restricted to the lowest irregular times (i.e. the coefficients of the irregular type closest to the residue). 
(Note also that due to the $\SL_2(\IC)$ action, these higher times will only give something new if $k=\abs{J}\ge 4$.)
\end{rmk}
\section{Nonlinear differential equations}\label{sn: differential equations}

In this section we will
write down the nonlinear Hamiltonian equations in various ways.
Note that all the equations here are equivalent to the (simpler looking) ``generic equations'' which appear when $J\subset \IC$ (see \eqref{eq: master} of \S\ref{ssn: generic equations}).

\begin{prop}\label{prop: main nlin eqns}
Horizontal sections of the connection $\wh \omega - d\varpi$ on 
$\IF\to \IB$ are determined by the differential equation:

$$%
d\Ga = \left[\wt{\delta(\Xi\Ga)},\Ga\right] 
+\left(\ga\Xi d\wh T + d\wh T\Xi \ga\right)^\circ
-\phi^{-1}\left(\Xi d\wh T\Xi\right) 
+%
\left(\begin{smallmatrix}
0 & -T_\infty dT_\infty PA\\
AQT_\infty dT_\infty & -TdTX - XTdT
\end{smallmatrix}\right)
$$

\noindent
where $\phi^{-1}:\End(V)\to \End(V)$
is defined in components by $\phi^{-1}(B_{ij})= B_{ij}/\phi_{ij}$ if $i\ne j$ and $0$ if $i=j$, and in general $R^\circ = R-\delta(R)$.
\end{prop}
\pf
This follows from Lemma \ref{lem: ham to eq}, and its analogue for the expression \eqref{eq: omega QPXB} for $\omega$. 
The calculation is simplified by noting the following, 
which is straightforward.
\begin{lem} \label{lem: useful}
If $R, F$ are a matrix valued function and one form respectively, then
(recalling that the commutator of matrices of one forms involves a plus sign):

i) $d(\wt R) = - \wt{dR}$, and ii) $\tr( \wt F R) = - \tr (F\wt R)$.
\end{lem}
Then it is easy, using the fact that $\de(\Xi\Ga) = -\de(\Ga\Xi)$, to verify that 
$d\tr(\wt{\Xi \Ga}\de(\Xi\Ga))/2 = \tr([\wt{\de(\Xi\Ga)}, \Ga ] d\Xi)$
yielding the first term of the equation.
Similarly the next two terms arise from $-\tr(\Xi\ga\Xi d\wh T)$,
and the final terms from the two remaining terms of $\varpi$.
\epf

Expanding in terms of $Q,P,B$ we obtain the equivalent equations:
\begin{align}\label{eq: QPB eqs}
dQ&=   Q\wt{PQ}  +  \wt{R}Q 
+ [dT,X]Q +(B+T)QdT_\infty 
+ dT Q T_\infty
+AQT_\infty dT_\infty,\notag\\
-dP&=  \wt{PQ}P  + P\wt R
+ P[dT, X]+dT_\infty P(B+T) 
+ T_\infty PdT 
+T_\infty dT_\infty PA,\\
dB &= 
[\wt{R},B] 
+[dT,QP] + BXdT + dTXB 
+[A,QdT_\infty P-XdTX]
+[T,[X, dT]],\notag
\end{align}
where $R= \de{(QP+XB)}$ and we tacitly apply $(\cdot)^\circ$ to the right-hand side of the third equation, and as usual $X=\ad_A^{-1}B$ and 
$\wt{(\cdot)} = \ad_t^{-1}[dt,\cdot]$ (with $t=T_i,T$ etc. as appropriate).

Alternatively one may rewrite the equations in terms of 
 the components $B_{ij}\in \Hom(W_j,W_i)$ of $\Gamma$ 
(for $i\ne j\in J$) as follows:
\begin{align}\label{eq: dBij}
\begin{split}
dB_{ij} =\quad
&\sum_{k} 
\wt{X_{ik}B_{ki}}B_{ij} 
+B_{ij}\wt{B_{jk}X_{kj}}\\
+&\sum_{k} dT_iX_{ik}B_{kj}
+ B_{ik}X_{kj}dT_j
-X_{ik}dT_k X_{kj}/\phi_{ij}\\
+&\ dT_iX_{ij}T_j+T_iX_{ij}dT_j -
\begin{cases}
T_idT_iX_{ij}+X_{ij}T_jdT_j &\text{if $i,j\ne\infty$, or}\\
T_\infty dT_\infty B_{ij}a_j&\text{if $i=\infty$, or}\\
-a_i B_{ij} T_\infty dT_\infty&\text{if $j=\infty$.}
\end{cases}
\end{split}
\end{align}
where  $X_{ij}= \phi_{ij}B_{ij}$ are the components of $\Xi$, and we set $B_{ii}=0$.

Now suppose we choose $j\in J$ and $i\in I_j$, 
so we have vector spaces 
$V_i\subset W_j$ 
and  $U_j=V\ominus W_j$.
If $i\in I_j$ we will also write $U_i:=U_j$.
The above equations may also be rewritten in terms of maps between these vector spaces.
To this end define
\beq\label{eq: PiQi def}
Q_i = \Gamma\circ \iota_i:V_i\to U_i,\qquad 
P_i = -\pi_i\circ\Xi:U_i\to V_i.
\eeq
Here $\pi_i$ and $\iota_i$ are the projection and the inclusion between $V$ and its summand $V_i$. (Note that the image of $Q_i$ is indeed in $U_i$, and for $P_i$ we tacitly restrict $\Xi$ to $U_i$.)
Thus all the data in $\Gamma$ is contained in the set of maps $\{Q_i\}$
and also in the set of maps $\{P_i\}$ as $i$ ranges over $I$.
Note also that if $j=\infty$ and $i\in I_\infty$ then
$Q_i,P_i$ are components of $Q,P$ (i.e.
$Q_i = Q\circ \iota_i$, and  
$P_i = \pi_i\circ P$)
so the notation \eqref{eq: PiQi def} is consistent with that used earlier, and now extends to all $i\in I$.

\begin{prop} \label{prop: local conns}
For all $j\in J$ and  $i\in I_j$
there are one-forms $\Omega_i$ on $\IB$ 
(depending explicitly on $\ga$) 
with values in $\End(U_i)$
such that
the equations \eqref{eq: dBij} are equivalent to the following overdetermined system: 
$$d Q_i = \Omega_i Q_i,\qquad dP_i = -P_i\Omega_i$$
for all $i\in I$.
\end{prop}
Note that just one of these two sets of equations is equivalent to 
the system \eqref{eq: dBij}. The point to note is that 
the same $\Omega_i$ appears in both equations.

\pf
If $i\in I_\infty$ then $\Omega_i$ is the restriction of the full connection $\Omega$ to the divisor $z=t_i$ (along which it has a logarithmic singularity): $\Omega_i = \Omega\bigl\vert_{z={t_i}}$. 
Indeed from the definition of the tilde operation it follows that
$$\sum_{j\in I_\infty\setminus\{i\}} Q_jP_jQ_i\frac{dt_i-dt_j}{t_i-t_j}=Q\wt{PQ}\circ\iota_i$$
and  similarly for $\pi_i\circ \wt{PQ}P$.
Thus from the first two equations in \eqref{eq: QPB eqs}, 
we see that
\beq\label{eq: Omi}
\Omega_i = 
\sum_{j\in I_\infty\setminus\{i\}} Q_jP_j\frac{dt_i-dt_j}{t_i-t_j}
+\wt R + [dT,X]  +d(T t_i) + (At_i +B)dt_i= \Omega\bigl\vert_{z={t_i}}\eeq
where $R=\de(QP+XB)$.
For the other possible $i\in I$ (say $i\in I_k, k\ne \infty\in J$) this may be shown by direct computation, from \eqref{eq: dBij}.
One finds that 
$$\Omega_i = 
\sum_{j\in I_k\setminus\{ i\}} Q_jP_j \frac{dt_i-dt_j}{t_i-t_j}+
{\pr}_k\left(\wt{\delta(\Xi\Ga)}
+(d\wh T) \Xi +
\left(\Ga dt_i
- \varphi_k^{-1}\Xi d\wh T
+d(t_i\wh T) +  C_i\right)\varphi_k\right)
$$
where $\pr_k:\End(V)\to \End(U_k)$ is the projection,
$\varphi_k = \id_k+\sum_{j\in J} \phi_{jk}\id_j\in \Aut(V)$ (where $\Id_j$ is the idempotent for $W_j$) and 
$C_i = \bsmx
 a_i T_\infty dT_\infty & \\ & -t_idt_i - TdT
\esmx\in \End(W_\infty\oplus U_\infty)$.
In fact the same expression works also for $k=\infty$ provided we replace $C_i$ by 
$\bsmx
 0 & \\ & At_idt_i
\esmx$.
\epf

Thus there is a direct geometric interpretation 
in the case when $i\in I_\infty$; 
the full connection $\Omega$ has a logarithmic singularity along the divisor $z=t_i$, with residue $R_i=Q_iP_i\in \End(U_\infty)$.
The one-form  $\Omega_i$ should be interpreted as a connection and as such is just the restriction of $\Omega$ to this divisor. 
Things have been arranged so that $Q_i$ is then a (collection of) horizontal sections of $\Omega_i$, and $P_i$ is then a (collection of) horizontal sections of the dual connection, so the residue is a horizontal section of the adjoint connection.
The geometry behind $\Omega_i$ for $i\in I_j, j\ne \infty$ is not so immediately transparent, but it may be obtained by performing a 
symplectic  transform to move $a_j\in \bP$ to $\infty$ and then using the above interpretation (taking care to do the explicit gauge transformations so the normalizations match up---this only involves the ``constants'' $C_i$).

\begin{cor}
Let $R_i:=Q_iP_i\in \End(U_i),\La_i:=-P_iQ_i\in \End(V_i)$.
If $\Gamma$ is a local horizontal section  then
$$dR_i = [\Omega_i,R_i],\qquad d\La_i=0$$
for all $i\in I$.
\end{cor}

Thus the adjoint orbit of each residue $R_i$ is preserved under the flow, for all $i\in I$. 
Note these orbits are also preserved under the symplectic transformations:

\begin{prop}
Under the action of $\SL_2(\IC)$ the orbit of $R_i\in \End(U_i)$ is preserved, and the value of $\La_i\in \End(V_i)$ is preserved.
\end{prop}
\pf
This follows from the definition \eqref{eq: PiQi def} of $P_i, Q_i$ together with the formulae for the $\SL_2(\IC)$ action on $\Ga,\Xi$ given in the proof of Proposition \ref{prop: invce of omega}. For example if $(\Ga,\Xi)\mapsto (\eps \Ga, \Xi \eps^{-1})$ then $\eps$ cancels in the definition of $\La_i$, and acts to conjugate $R_i$.
\epf

Finally note that the orbits of the elements $\La_i\in \End(V_i)$ cannot be assigned arbitrarily since $\La_i = -P_iQ_i =\pi_i\Xi\Ga\iota_i$ so that 
\beq\label{eq: fuchs reln}
\sum_{i\in I}\tr(\La_i) = \Tr(\Xi\Ga)=0
\eeq
due to the skew-adjointness of the map $\phi$.
(If we project to meromorphic connections this corresponds to the sum of the traces of the residues being zero.)

\subsection{Projected equations}\label{ssn: projected eqns}

Considering the projection $\Ga\mapsto \cA$ to the space of meromorphic connections it is easy to see that the nonlinear equations descend (to equations on the coefficients of $\cA$), as follows.
Recall that
$$\cA=\left(Az+B+T+\sum_{i\in I_\infty} \frac{R_i}{z-t_i}\right)dz$$
on the trivial bundle $U_\infty\times \IP^1$, where $R_i=Q_iP_i\in \End(U_\infty)$.
The corresponding nonlinear equations are
\beq\label{eq: inter eqs}
\begin{cases}
\text{the third equation in \eqref{eq: QPB eqs} for $dB$, and} \\
\text{the equations 
$dR_i = [\Omega_i,R_i]$ with $i\in I_\infty$}
\end{cases}\eeq
where 
$\Omega_i$ is as in \eqref{eq: Omi}.
(It is easy to see that these equations only depend on $\{R_i\}$  and not on their lifts $Q,P$; for example $QdT_\infty P = \sum_{I_\infty}R_idt_i$.)
These are the equations which arise as the condition for the vanishing of the curvature of the full connection
\beq\label{eq: full conn en bas}
\Omega = 
(Az+B)dz + d(zT)  + \sum_{I_\infty}R_i \frac{dz-dt_i}{z-t_i} +  
[dT,X] +\wt{R}
\eeq
with $R = \de(XB+\sum_{I_\infty} R_i)$.
The point is that by lifting up to $\IM$ we see the symplectic $\SL_2(\IC)$ symmetry group which is not apparent at this (intermediate) reduced level.
The Hamiltonians also descend to this level as follows.
Fix adjoint orbits $\cO_i\subset \End(U_\infty)$ for $i\in I_\infty$ and define
$$\wt\cM^* = \cO_B\times \prod_{I_\infty}\cO_i$$
where $\cO_B = \{(Az+T+B)dz\st B\in \Im(\ad_A)\subset \End(U_\infty)\}$
is the coadjoint orbit through $(Az+T)dz$ under the group of jets at $\infty$ of gauge transformations tangent to the identity (cf. Appendix \ref{apx: AKS} and \cite{smid} section 2).
Then observe that the Hamiltonians on $\IF$ 
descend to the total space of the trivial symplectic fibration
$$\wt\cM^*\times \IB \to \IB.$$
Indeed this follows from the expressions  \eqref{eq: infty terms} and \eqref{eq: varpi-infty}
noting e.g. that $QT_\infty P = \sum_{I_\infty} R_i t_i$ and
$$\tr(PQ\wt{PQ}) = 
\sum_{i\ne j\in I_\infty}\tr(R_iR_j)\frac{dt_i-dt_j}{t_i-t_j}.$$
The space
$\wt\cM^*$  is a symplectic leaf of the quotient $\IM/H_\infty$
where $H_\infty=\prod_{I_\infty} \GL(V_i)$
(cf. Lemma \ref{lem: mmap}).
Basic example of these projected equations
\eqref{eq: inter eqs}
are the Schlesinger equations and the JMMS equations, as will be made explicit in the following section.

Note that we view the projected equations \eqref{eq: inter eqs} as essentially equivalent to the full equations of Proposition \ref{prop: main nlin eqns}, since one may lift any solution of \eqref{eq: inter eqs} by only solving linear differential equations (as in \cite{k2p} Proposition 15).

\section{Examples}

\subsection{JMMS}
Consider the case where  $J=\{0,\infty\}\subset \bP$ so $A=B=0$, 
and
$$
\Ga = \bmx 0 & P \\ Q & 0 \emx,\qquad
\Xi = \bmx 0 & -P \\ Q & 0 \emx,\qquad
\wh T = \bmx T_\infty & 0 \\ 0 & T_0 \emx\in \End(W_\infty\oplus W_0).
$$
The phase space $\IM$ is 
$\{(P,Q)\}=\Hom(W_0,W_\infty)\oplus\Hom(W_\infty,W_0)$
with the symplectic form 
$\tr(d\Xi\wedge d\Ga)/2 = \tr(dQ\wedge dP).$
The Hamiltonian one-form then is:
$$\varpi = 
\frac{1}{2}\tr (Q\wt{PQ}P) + \frac{1}{2}\tr (P\wt{QP}Q)
+\tr (PT_0QdT_\infty) +\tr (QT_\infty PdT_0).$$
In this situation the full connection specializes to:
$$\Omega = 
 d(zT_0)  + Qd\log(z-T_\infty)P + \wt{QP}
$$
on the vector bundle $W_0\times\IP^1\times \IB \to \IP^1\times \IB$, where the space of times $\IB$ is 
$\{\wh T\} \cong  
(\IC^{\abs{I_0}}\setminus\diagonals)\times(\IC^{\abs{I_\infty}}\mdiagonals).$
The nonlinear differential equations \eqref{eq: QPB eqs}  are:
\beq \label{eq: altJMMS}
\begin{split} 
dQ = Q\wt{PQ} +\wt{QP}Q  + T_0 Q dT_\infty + dT_0 Q  T_\infty  \\
-dP = P\wt{QP} + \wt{PQ}P + T_\infty P dT_0 + dT_\infty P  T_0.
\end{split}
\eeq
These are equivalent to the JMMS equations \cite{JMMS};
They may be rewritten as follows.
Let $\{V_i\}$ be the eigenspaces of 
$T_\infty\in\End(W_\infty)$ (labelled by $i\in I_\infty$).
Let $\iota_i:V_i\to W_\infty$, $\pi_i:W_\infty\to V_i$ be the corresponding inclusions and projections.
Write 
$$Q_i = Q\circ \iota_i:V_i\to W_0,\qquad 
P_i = \pi_i\circ P:W_0\to V_i$$ 
for the corresponding components of $P$ and $Q$ respectively.
Proposition \ref{prop: local conns} then implies:
\begin{cor}
The equations \eqref{eq: altJMMS} may be rewritten as:
$$
dQ_i = \Omega_i Q_i,\qquad
dP_i = -P_i\Omega_i$$
for all $i\in I_\infty$,
where
$\Omega_i = \Omega\bigl\vert_{z=t_i} = 
d(t_iT_0)  + \wt{QP} + \sum_{j\ne i\in I_\infty} Q_jP_jd\log(t_i-t_j)$.
\end{cor}
These equations are the lifted JMMS equations which appear in  \cite{JMMS} (A.5.9). The JMMS equations themselves (\cite{JMMS} 4.44 or A.5.1) 
correspond to our projected equations \eqref{eq: inter eqs} obtained by setting $R_i=Q_iP_i$, which in this case, since $B=0$,  are
$$
dR_i = [\Omega_i, R_i]\qquad\text{where}\qquad
\Omega_i = 
 d(t_iT)+ \wt R + \sum_{j\in I_\infty\setminus\{i\}} R_j\frac{dt_i-dt_j}{t_i-t_j}
$$
and $R=\sum_{I_\infty}  R_i$.

\begin{rmk}
Strictly speaking  \cite{JMMS} only considers Hamiltonians in the case where $T_0$ has distinct eigenvalues, although as they write ``the general case is treated with minor modifications''. Under this restriction Harnad's symmetry only appears upon restricting $T_\infty$ to also have distinct eigenvalues, so that the residues $Q_iP_i$ are all rank one matrices.  This special case is highlighted on p.155 of \cite{JMMS} and the symmetry of the Hamiltonians may be seen (for this case) in \cite{JMMS}  equation (A.5.16).
\end{rmk}

\subsection{Schlesinger}
Now suppose we are in the special case of the JMMS equations where $T_0=0$.
Then the full connection specialises to:
$$\Omega =  Q\Delta P$$
where $\Delta=d\log(z-T_\infty)$, and the Hamiltonian one-form is
$\varpi= \tr(PQ\wt{PQ})/2$. 
The nonlinear differential equations are:
\beq\label{eq: lifted Schles}
dQ = Q\wt{PQ},\qquad dP =  -\wt{PQ}P.
\eeq
These equations are equivalent to the Schlesinger equations---more precisely they are the lifted Schlesinger equations, and imply the Schlesinger equations by projection as above.
Namely
the equations \eqref{eq: lifted Schles}
are equivalent to the equations
$$dQ_i = \Omega_i Q_i,\qquad dP_i = - P_i\Omega_i$$
where $\Omega_i=\Omega\bigl\vert_{z=t_i} = 
\sum_{j\ne i\in I_\infty} Q_jP_jd\log(t_i-t_j)$, and so
if we write 
$R_i = Q_iP_i\in\End(W_0)$ then
$$dR_i = [\Omega_i,R_i] = -\sum_{j\ne i} [R_i,R_j]d\log(t_i-t_j)$$
which are the Schlesinger equations \cite{schles-icm1908} p.67.

\subsection{Dual Schlesinger}
Suppose again we are in the situation of JMMS but 
instead that $T_\infty=0$. 
Then the nonlinear differential equations are:
$$dQ = \wt{QP}Q\qquad dP = -P\wt{QP}.$$
In this situation the full connection specializes to:
$$\Omega = 
 d(zT_0)  + QP\frac{dz}{z} +  \wt{QP},$$
and the Hamiltonians are $\varpi= \tr(QP\wt{QP})/2$. 
Setting $R=R_0= QP\in \End(W_0)$ the projected equations are
$$dR = [\wt R, R]$$
i.e. $dR = [\ad_T^{-1}[dT,R],R]$, with $T=T_0$.
Some special cases of these equations control semisimple Frobenius manifolds \cite{Dub95long} (3.74).
(By Harnad duality they are equivalent to certain cases of Schlesinger's equations---in fact they are within the special case of rank one residues considered in \cite{JMMS} p.155). 
A generalisation of these dual Schlesinger equations to arbitrary complex reductive groups $G$ was studied in \cite{bafi}---then the space of times is the regular part of a Cartan subalgebra (whose fundamental group is the pure  $G$-braid group).

\subsection{Generic equations}\label{ssn: generic equations}

Since the symplectic transformations enable us to move $J$ around in the Fourier sphere $\bP$ by M\"obius transformations, 
we see generically no point of $J$ will lie at $\infty\in\bP$.
In this case the full connection is
$$\Omega = 
(Az+B)dz + d(zT)   +  
[dT,X] +\wt{\delta(XB)}
$$
on the trivial bundle with fibre $V$  (since $W_\infty=0, U_\infty=V$---in effect
$P=Q=0$ and $T=\wh T$).
Upon restriction to $\IP^1$ this has just one pole of order three at $z=\infty$ and no others.
The Hamiltonians are 
\beq
\varpi = \frac{1}{2}\tr\left(\wt{XB}\delta(XB)\right)
-\tr\left(XBX dT\right) + \tr(X[X,T]dT)
\eeq
 and
the nonlinear equations are
\beq\label{eq: master}
dB = 
[\wt{\de(XB)},B] 
+[[dT,X],B+T]^\circ.
\eeq
These are the ``master equations'' in the sense that any of the other nonlinear equations considered here are equivalent to equations of this form (by moving $J\subset \bP$ so that $\infty\not\in J$).
Note that if $A$ has distinct eigenvalues 
(as in the work \cite{JMU81} of Jimbo--Miwa--Ueno) then $\wt{\de(XB)}=0$ and the equations are simpler.

For example if we consider the bipartite case with $J=\{0,1\} \subset \bP\setminus\infty$ so that $V=W_0\oplus W_1$, then upon writing
$B = \bsmx 0 & R\\ S & 0 \esmx$
the equations \eqref{eq: master} become
\begin{align*} 
dS &= S\wt{RS} +\wt{SR}S  \,+ T_1 S dT_0 + dT_1 S  T_0 -(ST_0dT_0 + T_1dT_1S)  \\
-dR &= R\wt{SR} + \wt{RS}R + T_0 R dT_1 + dT_0 R  T_1 - 
(RT_1dT_1 + T_0dT_0 R).
\end{align*}
The terms in parentheses can be gauged away, and we obtain the lifted JMMS equations \eqref{eq: altJMMS}.
In other words we have an alternative Lax pair for the JMMS equations, as controlling isomonodromic deformations of $\cA = (Az+B+T)dz$ on $V\times \IP^1$ with $V=W_0\oplus W_1$.
In turn by specialising (e.g. to $T_1=0$) this gives a new Lax pair for the Schlesinger equations, and specialising further even for Painlev\'e VI (cf. \S\ref{ssn: PVI}).

\section{Reductions, moduli spaces and relation to graphs}\label{sn: reductions}

In summary we have defined and studied a flat nonlinear connection on the bundle
$\IF = \IM\times \IB \to \IB$ and shown that it is invariant under an action of $\SL_2(\IC)$, up to some simple explicit gauge transformations. It is also invariant under the natural action of the group $\wh H= \Prod \GL(V_i)$.
The aim of this section is to consider the reductions of the nonlinear connection under this group of automorphisms $\wh H$ (rather than just the subgroup $H_\infty\subset \wh H$ considered in \S\ref{ssn: projected eqns}).
In brief this amounts to replacing $\IM$ by its symplectic quotient
$$\cM^* = \IM\spqa{\breve \cO} \wh H$$

\noindent
by $\wh H$ at a coadjoint orbit $\breve\cO$.
To make life simpler we will restrict to the subset of {\em stable} points $\cM^*_{st}\subset \cM^*$ throughout this section.

The main results are that 1) the resulting nonlinear connection on $\cM^*\times \IB\to \IB$ is (completely) invariant under the symplectic group of transformations, and 2) that after reduction extra symmetries become apparent that, when  combined with the symplectic transforms, immediately give the desired action of a Kac--Moody Weyl group.

This Kac--Moody Weyl group action simultaneously generalises that of 1) Okamoto \cite{Okamoto-dynkin} (in the theory of Painlev\'e equations, when $\dim_\IC(\cM^*)=2$), and 
2) Crawley--Boevey \cite{CB-additiveDS} (in the case of Fuchsian systems, corresponding to star-shaped Kac--Moody Dynkin graphs).

As we will explain, from the point of view of irregular connections this action is not at all mysterious: it basically amounts to changing the possible choices of orderings of the eigenvalues of the residues $R_i$ for all $i\in I$. Of course at any given moment, only the $R_i$ with $i\in I_\infty$ will appear as residues at simple poles, and we should use the symplectic transformations to realise the other $R_i$ as residues.

First of all some basic definitions related to graphs will be given.

\subsection{Representations of graphs}\label{ssn: reps of graphs}

Suppose $\cG$  is a graph with nodes $I$ (and edges $\cG$).
Let $\bar\cG$ be the set of oriented edges of $\cG$, i.e. the set of pairs $(e,o)$ such that  $e\in \cG$ is an edge of $\cG$ and $o$ is a choice of one of the two possible orientations of $e$.
Thus if $a\in \bar \cG$ is an oriented edge, the head $h(a)\in I$ and tail $t(a)\in I$ nodes of $a$  are well defined. 
For our purposes it is convenient to define a {\em representation} of the graph $\cG$ to be the following data

1) an $I$-graded vector space $V= \bigoplus_{i\in I} V_i$, and

2) for each oriented edge $a\in \bar\cG$,
a linear map $v_a : V_{t(a)}\to V_{h(a)}$ between the vector spaces at the tail and the head of $a$.

Thus the data in 2) amounts to choosing a linear map in both directions along each edge of $\cG$.
A {\em subrepresentation} of a representation $V$ of $\cG$ consists of an $I$-graded subspace $V'\subset V$ which is preserved by the linear maps, i.e. such that 
$v_a(V'_{t(a)})\subset V'_{h(a)}$ for each oriented edge 
$a\in \bar \cG$.
A representation $V$ is {\em irreducible} if it has no proper nontrivial subrepresentations.
Given an $I$-graded vector space $V$ we may consider the set 
$\Rep(\cG,V)$ of all
representations of $\cG$ on $V$. This is just the vector space
$$\Rep(\cG,V)=\bigoplus_{a\in \bar\cG} \Hom(V_{t(a)},V_{h(a)})$$  
of all possible maps.

\subsection{The Kac--Moody root system and Weyl group}\label{ssn: km}

Let $\cG$ be a graph with no edge loops.
Then one can define a (symmetric) Kac--Moody root system and Weyl group as follows.
Let $I$ be the set of nodes and let $n=\abs{I}$ be the number of nodes.
Define the $n\times n$ (symmetric) Cartan matrix to be 
$$C=2\ \id - A$$
where $A$ is the adjacency matrix of $\cG$; 
the $i,j$ entry of $A$ is the number of edges connecting the nodes $i$ and $j$.
The {\em root lattice} $\IZ^I=\bigoplus_{i\in I} \IZ\eps_i$ inherits a bilinear form defined by 
\beq\label{eq: KM bil form}
(\eps_i,\eps_j) = C_{ij}.
\eeq

The simple reflections $s_i$, acting on the root lattice, are defined  by the 
formula
$$s_i(\be) := \be-(\be,\eps_i)\eps_i$$
for any $i\in I$.
They satisfy (cf. \cite{Kac-book} p.41) the relations
$$s_i^2=1,\qquad 
s_is_j=s_js_i \text{\ if $A_{ij}=0$,}\qquad
s_is_js_i=s_js_is_j \text{\ if $A_{ij}=1$}.
$$
By definition the Weyl group is the group generated by these simple reflections.
There are also dual reflections $r_i$ acting on the vector space 
$\IC^I$ by the formula 
$$r_i(\la) = \la - \la_i \al_i$$
where $\la = \sum_{i\in I} \la_i\eps_i\in \IC^I$ with $\la_i\in\IC$ and 
$\al_i := \sum_{j} (\eps_i,\eps_j)\eps_j\in \IC^I.$ 
By construction one has that $s_i(\be)\cdot r_i(\la) = \be \cdot \la$,
where the dot denotes the pairing given by 
$\eps_i\cdot\eps_j = \delta_{ij}$.

The corresponding Kac--Moody root system is a subset of the {root lattice} 
$\IZ^I$.
It may be defined as the union of the set of real roots  and the set of imaginary roots, where

1) The simple roots are $\eps_i$ for $i\in I$,

2) The set of real roots is the Weyl group orbit of the set of simple roots,

3) Define the {\em fundamental region} to be the set of nonzero 
$\be\in\IN^I$ whose support is a connected subgraph of $\cG$ and such that $(\eps_i,\be)\le 0$ for all $i\in I$.
The set of imaginary roots is the union of the Weyl group orbit 
of the fundamental 
region and the orbit of minus the fundamental region.

This defines the root system (see \cite{Kac-book} Chapter 5 for the fact that this description does indeed give the roots of the corresponding Kac--Moody algebra).
By definition a root is positive if all its coefficients are $\ge 0$.
For example if $\cG$ is an $ADE$ Dynkin diagram this gives the root system of the corresponding finite dimensional simple Lie algebra, or if $\cG$ is an extended/affine $ADE$  Dynkin diagram then this is the root system of the corresponding affine Kac--Moody Lie algebra (closely related to the corresponding loop algebra), but of course there are many examples beyond these cases.

\subsection{Complete $k$-partite graphs}\label{ssn: k-partite graphs}

Let $\cG$ be a graph with nodes $I$ (and edges $\cG$).
Recall that by definition $\cG$ is a {\em complete $k$-partite graph}
if there is a partition $I=I_1\sqcup \cdots \sqcup I_k$ of its nodes (into $k$ nonempty parts $I_j$) such that two nodes are connected by a single edge if and only if they are not in the same part.
Thus there is a bijection between the set of partitions with $k$ parts, and the set of complete $k$-partite graphs.
Let $\cG(P)$ denote the complete $k$-partite graph corresponding to a partition $P$ (thought of equivalently either as a Young diagram, or as a partition of a finite set, or as a partition of an integer, or as a surjective map $\phi:I\to J$ onto the set $J$ of parts, so that $I_j=\phi^{-1}(j)$ for all $j\in J$).
For example the graph $\cG(1,1)$ corresponding to the partition $1+1$ is just a single edge connecting two nodes, and similarly $\cG(1,1,1)$ is the triangle and $\cG(2,2)$ is the square (a complete bipartite graph).
The star-shaped graph with $n$ legs of length one, is the bipartite graph $\cG(1,n)$. 
The graphs $\cG(n)$ have $n$ nodes and no edges, and
the graphs $\cG(1,1,\ldots,1)$ are the complete graphs (with every pair of nodes connected by a single edge).
See Figure \ref{fig: partite-intro} of the introduction.

\begin{defn}\label{defn: supernova graph}
A {\em (simply-laced) supernova graph} is a graph obtained by gluing a single leg (of arbitrary length $\ge 0$) onto  each node of a complete $k$-partite graph. 
\end{defn}

Here a ``leg'' of length $l$ is just a  Dynkin graph of type 
$A_{l+1}$, with $l$ edges. 
For example any star-shaped graph (with arbitrary length legs) may be viewed as a supernova graph with central subgraph of the form $\cG(1,n)$.
This motivated the name ``supernova'', as a star with more going on in the middle. 
(The not-necessarily simply-laced symmetric 
supernova graphs are the graphs described in appendix C of \cite{rsode}.)

In the next three subsections we will give different viewpoints on the symplectic reduction of $\IM$. Each viewpoint is useful for different reasons and
going between these different viewpoints yields the Kac--Moody reflections.

\subsection{Moduli of Weyl algebra module presentations}\label{ssn: wamps}

The first viewpoint is as the symplectic reductions of the space $\IM$ of presentations of modules for the first Weyl algebra.
Here the data we need is as follows.
Choose a finite set $J$ and an embedding $\ba:J\hookrightarrow \bP = \IC\union \{\infty\}$. Write $a_j = \ba(j)$ for $j\in J$. For each $j\in J$ choose a finite set $I_j$ and 
write $I = \bigsqcup I_j$.
Choose a finite dimensional complex vector space $V_i$ for each $i\in I$, and write $W_j = \bigoplus _{i\in I_j} V_i, 
V = \bigoplus _{i\in I} V_i.$
Finally choose an adjoint orbit $\breve \cO_i\subset \End(V_i)$ for each $i\in I$ and write $\breve \cO=\Prod \breve\cO_i$.
Note that we place no restriction on these orbits (for example the eigenvalues may be integral or  differ by integers).
In particular this data determines the symplectic manifold $\IM$ as in \S\ref{sn: sympl}.

\begin{prop}\label{prop: mmap La_i}
The group $\wh H=\prod \GL(V_i)$ acts on $\IM$ in a Hamiltonian fashion with moment map
$$\mu:\IM\to \prod\End(V_i);\qquad \Ga\mapsto (\La_i).$$
The induced action on $\IF=\IM\times \IB$ preserves the isomonodromy connection.  
\end{prop}
\pf
Here $\End(V_i)$ is identified with the dual of the Lie algebra of $\GL(V_i)$ by the trace pairing. 
The action is defined by $g(\Ga) = g\Ga g^{-1}$ for $g\in \wh H\subset \GL(V)$.
The nonlinear 
equations are invariant since this action commutes with $\phi:\End(V)\to \End(V)$ and clearly 
$g\wh T g^{-1} = \wh T$.
The moment map computation is straightforward: 
for example for any fixed $i\in I_j$ we may 
split $\IM$ symplectically as a product $\IM'\times T^*\Hom(V_i,U_j)$
(with symplectic form $\tr(dQ_i\wedge d P_i)$ on the cotangent bundle---cf. \eqref{eq: omega QPXB} in the case $j=\infty$ and the general case is similar). Then $g\in \GL(V_i)$ acts trivially on $\IM'$ and the action on the  cotangent bundle is easily seen to have moment map $\La_i = -P_iQ_i$.
\epf

Recall that the stable points of $\IM$ for the action of $\wh H$ (in the sense of Mumford's geometric invariant theory) are those points
whose $\wh H$ orbit is closed and of maximal possible dimension. 
These points may be described in terms of graph representations as follows. 
Let $\cG$ be the complete $k$-partite graph with nodes $I$ corresponding to the partition $I=\bigsqcup I_j$, so that $k=\abs{J}$.

\begin{prop}\label{prop: M as graph reps}
1) The space $\IM$ is isomorphic to the space $\Rep(\cG,V)$ of representations of the graph $\cG$ on the $I$-graded vector space $V$.

2) The stable points of $\IM$ for the action of $\wh H$ are the irreducible representations. 
\end{prop}
\pf
1) is straightforward: specifying a point of $\IM$ is the same as  the choice of linear maps $b_{ij}:V_j\to V_i$
for all $i,j\in I$ not in the same part.
Then 2) is a special case of a result of King \cite{king-quivers} on quiver representations.
\epf

Thus we may now perform the symplectic quotient of the stable part of $\IM$ by $\wh H$, at the (co)adjoint orbit $\breve\cO$ of $\wh H$.
Namely we define the moduli space of stable Weyl algebra module presentations to be
$$\cM_{st}^* = \cM_{st}^*(\cG,\breve\cO) = \{ \text{ stable points $\Ga\in \IM$ such that $\mu(\Ga)\in \breve\cO$ } \}/ \wh H$$
where $\mu$ is the moment map for the $\wh H$ action.
We will see below in Theorem \ref{thm: cM and Q} 
that this is a smooth (possibly empty) symplectic algebraic variety.

\begin{cor}
The isomonodromy connection on $\IF\to \IB$ 
descends to define a nonlinear connection on 
$\cM_{st}^*\times \IB\to \IB$, and this reduced connection is completely invariant under the symplectic action of $\SL_2(\IC)$.
\end{cor}
\pf
We have already shown the connection is $\wh H$ invariant, and so descends. For the symplectic invariance, we have already shown the connection on $\IF$ is invariant up to some explicit gauge transformation. But these gauge transformations act within the orbits of the $\wh H$ action, and so are trivial after reduction.
Said differently once the orbits of the elements $\La_i$ are fixed then the gauge terms changing the isomonodromy Hamiltonians are {\em constant} so that the (reduced) isomonodromy connections are the same (with their Hamiltonians differing by constants).
\epf

One consequence of stability that will be useful is the following.
Recall that given $\Ga\in \IM$ we have defined, for any $i\in I_j$, two linear maps $Q_i:V_i \to U_j$ and $P_i:U_j\to V_i$
(see \eqref{eq: PiQi def}).
 
\begin{lem}\label{lem: inj-surj}
If $\Ga\in \IM$ is stable then $Q_i$ is injective and $P_i$ is surjective, for all $i\in I$.
\end{lem}
\pf
Thinking in terms of representations of the complete $k$-partite graph $\G$, $P_i$ encodes all of the maps to $V_i$ and $Q_i$ encodes all the maps from $V_i$.
Thus if $P_i$ was not surjective, then we could replace $V_i$ by the image $P_i(U_j)\subset V_i$ to obtain a nontrivial proper subrepresentation, contradicting stability. 
If $Q_i$ was not injective, then we could choose a nonzero subspace 
$K\subset V_i$ of its kernel to define a subrepresentation (taking the zero vector space at all other nodes), again contradicting stability.
\epf

\begin{cor}
If $\Ga\in \IM$ is stable then fixing the adjoint orbit of $\La_i\in \End(V_i)$ is equivalent to fixing the adjoint orbit of $R_i\in \End(U_i)$. 
\end{cor}
\pf
Since $R_i=Q_iP_i$ and $\La_i=-P_iQ_i$, this follows from the injectivity/surjectivity conditions in Lemma \ref{lem: inj-surj}. The exact relation between the orbits is summarised in the Appendix \ref{apx: relating orbits}.
\epf

\begin{rmk}
One can also consider the naive symplectic quotient
$$\cM^* = \{ \text{ $\Ga\in \IM$ such that $\mu(\Ga)\in \breve\cO$ } \}/ \wh H.$$
Note that for sufficiently generic orbits $\breve\cO$ this coincides with $\cM_{st}^*$, since for any subrepresentation $V'\subset V$ one has 
\beq\label{subrep reln}
\sum\tr(\La'_i)= 0
\eeq
(where $\La_i'\in \End(V'_i)$ is the analogue of $\La_i$ determined by the subrepresentation $V'$) and the eigenvalues of $\La'_i$ will be a subset of the eigenvalues of $\La_i$, so that if $\breve\cO$ is sufficiently generic there will be no relation of the form \eqref{subrep reln} and thus no proper nontrivial subrepresentations. 
Thus in such cases $\cM^*$ itself is a smooth algebraic variety.
\end{rmk}

\subsection{Representations of supernova graphs}
Now we will describe the above moduli spaces in terms of representations of supernova graphs.
In brief the choice of the orbits $\breve \cO_i$ is replaced by the choice of a scalar on each node of the legs;  describing the choices in this way enables us to see the underlying Kac--Moody root system.
In turn this enables us to attach isomonodromy equations to representations
of supernova graphs.
(If the graph is star-shaped the equations will be equivalent to the Schlesinger equations, and more generally if the central part of the graph is bipartite, i.e. $\abs{J}=2$, then  the equations will be equivalent to the JMMS equations.)

Let $\wh \cG$ be a supernova graph with nodes $\wh I$, as defined in 
\S\ref{ssn: k-partite graphs}.
Thus $\wh \cG$ consists of a complete $k$-partite subgraph 
$\cG\subset \wh \cG$ with nodes $I\subset \wh I$  and 
a `leg' glued on to the $i$th node for each $i\in I$.
We will call $\cG$ the {\em core} of $\wh \cG$; it is uniquely determined except in the star-shaped case.
Let $I = \bigsqcup_{j\in J} I_j$ be the partition of $I$ into parts $I_j$
labelled by the set $J$, so that  $k=\abs{J}$.
Let $l_i \in \IZ_{\ge 0}$ be the length of the $i$th leg for each $i\in I$.

The further data required to construct the desired spaces are as follows:

$\bullet$ a complex number $\la_i\in \IC$ for each $i\in \wh I$, 

$\bullet$ a distinct point of the Fourier sphere $a_j\in\bP$ for each part $j\in J$,

$\bullet$ an integer $d_i\ge 0$ for each $i\in \wh I$.

Given such integers $d_i$ define $V_i=\IC^{d_i}$ and thus a vector space  $\wh V = \bigoplus_{i\in \wh I} V_i$ graded by $\wh I$.
The space of times for the corresponding isomonodromy equations 
is  
$\IB = \Prod_{j\in J} \left(\IC^{\abs{I_j}}\setminus\diagonals\right)\subset \IC^I$.
In other words $\IB$ consists of 
the sequences $\{t_i\st i\in I\}\in \IC^I$ such that 
$t_i\ne t_{i'}$ whenever $i$ and $i'$ are in the same part of $I$, i.e. there is a time variable $t_i\in \IC$ for each node of the core, and the times in the same part should be pairwise distinct.
Given such data let
$$\wh \IM = \Rep(\wh \cG, \wh V)$$ 
be the space of representations of $\wh \cG$ on $\wh V$.
By dividing up the edges of $\wh \G$ into those in the core or in a leg there is a product decomposition:
$$\wh \IM = \IM \times \IL$$
where $\IM$ is as above (identified with the space of representations of $\cG$ on $V=\bigoplus_{i\in I} V_i$) and
$\IL = \Prod_{i\in I} \IL_i$ where $\IL_i$ is the space of representations of the $i$th leg:
$$\IL_i=\prod_{j=1}^{l_i}T^*\Hom(V_{i,j},V_{i,j+1})$$
where $V_{i,j}\subset \wh V$ 
denotes the vector space on the $j$th node of the $i$th leg (labelled going down the leg so that $V_{i,1}=V_i$ for all $i\in I$).
Here we have identified $\Hom(V,W)\oplus\Hom(W,V)$ with the cotangent bundle $T^*\Hom(V,W)$, thus giving each $\IL_i$ and therefore $\IL$ a complex symplectic structure (using the natural symplectic structure on the cotangent bundle).
Thus $\wh \IM$ inherits a product  symplectic structure 
from that on  $\IL$ and the (nonstandard) symplectic form of \S\ref{sn: sympl} on $\IM$, depending on the chosen points $a_j\in \bP$ (and defining $W_j = \bigoplus_{i\in I_j} V_i$ to relate the definition of $\IM$ in \S\ref{sn: sympl} to the present definition).

A point of $\IL_i$ will be denoted $(\bp_i,\bq_i)$ with
$\bp_i = (p_{i1},p_{i2},\ldots), \bq_i = (q_{i1},q_{i2},\ldots)$
with $p_{ij} : V_{i,j}\to V_{i,j+1}$ and
 $q_{ij} : V_{i,j+1}\to V_{i,j}$
for $j=1,2,\ldots,l_i$.
The group
$$\wh \IG = \prod_{i\in \wh I}\GL(V_i)$$
acts naturally on $\wh \IM$ via its action on the spaces $V_i$, and preserves the symplectic structure.
Moreover it is a Hamiltonian action with moment map as follows:

$$\wh\mu_i = \La_i - q_{i1}\circ p_{i1}\in \End(V_i)$$
if $i\in I$, and
$$\wh\mu_{ij} = p_{i(j-1)}\circ q_{i(j-1)} - q_{ij}\circ p_{ij} \in \End(V_{i,j})$$
for the $j$th node of the $i$th leg, for $j>1$.
Write $$\wh \mu : \wh \IM \to \Prod_{i\in \wh I} \End(V_i)$$
for the resulting moment map, with components $\wh \mu_i$ and $\wh \mu_{ij}$.  (The right-hand side is identified with the dual of the Lie algebra of $\wh \IG$ using the trace pairing on each factor.)

Now identify the chosen complex number $\la_i\in \IC$ 
with the scalar matrix $\la_i \Id_{V_i}\in \End(V_i)$ for each $i\in \wh I$
and denote this $\wh I$-tuple by $\bla\in \Prod \End(V_i)$. 
Since $\bla$ is a central element of the Lie algebra of $\wh \IG$ we can equally well view the point $\{ \bla \}$ as a (co)adjoint orbit of $\wh \IG$.
Then we may define the ``twisted quiver variety'' associated to this data  to be
$$\cQ(\wh \cG, \bla, \bd) = \wh \IM^{st} \spqa{\bla} \wh \IG = \{ \text{ stable } \rho\in \wh\IM \st \wh\mu(\rho)=\bla\}/\wh \IG $$
as the complex symplectic quotient of the subset of stable points of $\wh \IM$ by the group $\wh \IG$, at the value $\bla$ of the moment map.
Since $\wh \IM$ is the space of representations of a graph (and $\wh \IG$ is the full automorphism group),  
\cite{king-quivers} again implies the stable points of $\wh \IM$ are the irreducible representations.

\begin{prop}\label{prop: qv constr}
The twisted quiver variety $\cQ(\wh \cG, \bla, \bd)$
is a smooth complex symplectic algebraic variety, which is either empty or of dimension $2-(\bd,\bd)$, where the inner product $(\ ,\ )$ on the root lattice of $\wh \cG$ is defined in \S\ref{ssn: km}. 
\end{prop}
\pf
The group $\IG=\wh\IG/\IC^*$ acts freely (since
the stability condition implies the stabiliser is finite, and it then follows in the quiver setting that the stabiliser is trivial, cf. \cite{cas-slod} p.283). Smoothness follows as in \cite{cas-slod} p.270. The dimension computation is now straightforward, since $\dim\wh\IM=\bd\cdot A\bd$ and 
$\dim(\wh \IG) = \bd\cdot \bd$, where $A$ is the adjacency matrix of $\wh \cG$, so that $\dim(\cQ) = \bd\cdot A\bd - 2 (\bd\cdot \bd-1) = 2-(\bd,\bd)$.
\epf

Now we will relate these twisted quiver varieties to the moduli spaces 
$\cM_{st}^*(\cG,\breve\cO)$ of Weyl module presentations of \S\ref{ssn: wamps}.
The extra data needed to define an isomorphism from $\cM_{st}^*(\cG,\breve\cO)$ to a twisted quiver variety is a marking of each orbit $\breve \cO_i$.

\begin{defn}
Suppose $\cO\subset \gl_n(\IC)$ is an adjoint orbit.
A {\em `marking'} of $\cO$ is
a finite ordered set $(\xi_1,\xi_2,\ldots,\xi_w)$ of complex numbers such that $\Prod_1^w(A-\xi_i) = 0$ for any $A\in \cO$.
\end{defn}

\begin{figure}[h]
	\centering
	\input{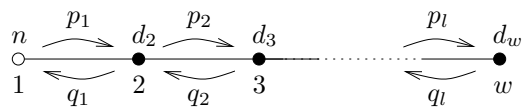}
	\caption{Representation of a type $A_w$ Dynkin graph.}\label{fig: leg}
\end{figure}

Equivalently a marking is the choice of a monic annihilating polynomial $f\in \IC[x]$, such that $f(A)=0$ for all $A\in \cO$, together with a choice of ordering of the roots of $f$.
A marking will be said to be {\em minimal} if $w=\deg(f)$ is minimal (so that $f$ is the minimal polynomial of $A\in \cO$).
A marking is {\em special} if the first root is zero ($\xi_1=0$).
Given a marking of $\cO\subset \gl_n(\IC)$, 
define  complex numbers 
\beq\label{eq: la defn}
\la_i = \xi_{i}-\xi_{i-1}
\eeq
(including $\la_1=\xi_1$) %
and integers 
$$d_i = \rank (A-\xi_1)\cdots (A-\xi_{i-1})$$ 
$i=1,2,\ldots$ (for any $A\in \cO$) so that $d_1=n$.
Then consider the type $A_w$ Dynkin graph (a leg) with $w$ nodes and $l:=w-1$ edges, as in Figure \ref{fig: leg}.

\begin{lem} (cf. \cite{CB-additiveDS}).\label{lem: orbits}  
If $\{(p_i,q_i)\}$ is a representation of this leg 
(type $A_w$ Dynkin graph) on 
the vector space $V=\bigoplus_1^w \IC^{d_i}$ such that each 
$p_i$ is surjective and each $q_i$ is injective and the moment map conditions
$$\La = q_1p_1+\la_1,\ p_{1}q_{1}=q_2p_2 + \la_2,\ \ldots, \ 
p_{l-1}q_{l-1}=q_lp_l + \la_l,\  p_lq_l=\la_w$$
hold, then $\La\in \cO$.
\end{lem}
\pf
It is clear that the orbit of $\La$ is uniquely determined by these conditions (cf. Proposition \ref{prop: orbit relations}).
Thus we just need to check that the orbit determined in this way is indeed the orbit 
$\cO$ we started with.
But if $\La$ is any element of $\cO$ we may define a representation satisfying these conditions by setting $V_1=\IC^n$ and then inductively 
$p_i=(\La - \xi_i)\bigl\vert_{V_i}, V_{i+1}=\Im(p_i)$ and taking 
$q_i$ to be the inclusion 
$V_{i+1}=
\Prod_1^i(\La-\xi_j)(V)
\hookrightarrow \Prod_1^{i-1}(\La-\xi_j)(V)= V_i$.
\epf

Now suppose we have data $\cG, V, \breve \cO, \ba$ used to define $\cM_{st}^*$ as in \S\ref{ssn: wamps}.
If we choose a marking of each orbit $\breve\cO_i$ then we can define a 
supernova graph $\wh \cG$ by gluing the first node of the leg corresponding to $\breve \cO_i$ onto the node $i$ of $\cG$.
Let $\wh I$ be the set of nodes of $\wh \cG$, so the markings determine a parameter $\la_i\in \IC$ and a dimension $d_i\in\IZ_{\ge 0}$ 
for all $i\in \wh I$.
This determines the variety $\cQ(\wh \cG, \bla, \bd)$.

\begin{thm} \label{thm: cM and Q}
The spaces $\cM^*_{st}(\cG,\breve\cO)$ and $\cQ(\wh \cG, \bla,\bd)$ are isomorphic.
\end{thm}
\pf
The first step is to see there is a well-defined map
$\wh\mu^{-1}(\bla)^{st}\to \mu^{-1}(\breve\cO)^{st}$ 
from the stable points of 
$\wh\mu^{-1}(\bla)\subset\wh\IM$ to the stable points of 
$\mu^{-1}(\breve\cO)\subset \IM$,
given by restricting to the subgraph $\cG$ a representation $\rho$ of $\wh\cG$ (on the vector space $\wh V = \bigoplus_{\wh I} V_i$).
Indeed if $\rho$ is stable then all the maps $p_{ij}$ down the legs are surjective and all the maps $q_{ij}$ up the legs are injective. (For example if $p_{i1}$ was not surjective, replace $V_{ij}$ by $V'_{ij}=p_{i(j-1)}\cdots p_{i1}(V_{i1})$ for $j>1$ (and $i$ fixed) to obtain a  non-zero subrepresentation, contradicting stability.)
Then using the condition $\wh \mu(\rho)=\bla$ and Proposition 
\ref{prop: orbit relations} repeatedly, implies 
$\La_i\in \breve\cO_i$ as in Lemma \ref{lem: orbits}.
The resulting point of  $\mu^{-1}(\breve\cO)$ is stable, since if it had a proper subrepresentation $V'$, then we could extend it to a subrepresentation of $\wh V$ as above by setting $V'_{ij}=p_{i(j-1)}\cdots p_{i1}(V'_{i1})$ for $j>1$ and all $i\in I$.
Thus the restriction map is well-defined. 
It is surjective since (up to choosing bases) we may define 
$p_{i1} = \La_i - \xi_{i1}$ (and identify $V_{i2}$ with its image, and 
$q_{i2}$ with the inclusion of this image in $V_{i1}$), and then repeat down each leg so $p_{i2} = (\La_i - \xi_{i2})\bigl|_{V_{i2}}$ etc.
The resulting point is in $\wh\mu^{-1}(\bla)$, and is clearly stable.
Further, since all the maps $p_{ij}$ are surjective and the $q_{ij}$ are injective, the fibres of the restriction map are precisely the orbits of the action of the group $\Prod_{\wh I\setminus I} \GL(V_i)$.
Finally the restriction map is equivariant under the group $\wh H = \Prod_{I} \GL(V_i)$ so the two quotients may be identified.
\epf

In particular we may denote this common space as 
$\cM^*_{st}(\wh \cG, \bla, \bd)$ or $\cM^*_{st}(\bla, \bd,\ba)$.

\begin{cor}
If $i\in \wh I\setminus I$ is a node which is not in the core, and $\la_i\neq 0$ then the spaces 
$\cM^*_{st}(\wh \cG, \bla, \bd)$ and $\cM^*_{st}(\wh \cG, r_i(\bla), s_i(\bd))$
are isomorphic (where $r_i,s_i$ are the reflections of \S\ref{ssn: km}).
\end{cor}
\pf
These reflections arise by changing the choice of marking.
Thus both spaces are isomorphic to the same space 
$\cM_{st}^*(\cG,\breve\cO)$ (which does not depend on the markings).
Explicitly if $i+1$ and $i-1$ denote the nodes adjacent to $i$ then under the reflections
$\bd$ is unchanged except for the component 
$$d_i\mapsto -d_i+d_{i-1}+d_{i+1}$$
and $\bla$ is unchanged except the components
$$(\la_{i-1},\la_i,\la_{i+1})\mapsto (\la_{i-1}+\la_{i},- \la_i,\la_{i+1}+\la_{i})$$
which, via the definition \eqref{eq: la defn} of $\bla$, indeed corresponds to swapping the two $\xi$'s occurring in $\la_i$.
\epf

To obtain the reflections corresponding to the core  nodes in this way we need to take a third viewpoint, that of meromorphic connections, and use the symplectic $\SL_2(\IC)$ transformations. This will be done in the next section.

\begin{rmk}\label{rmk: relation with quivers}
Note that by definition a representation of a graph is the same thing as a representation of the double of any quiver obtained by choosing an orientation of the graph, as appears for example in Nakajima's theory of quiver varieties \cite{nakaj-duke94}. In that theory the orientation of the graph is used to determine the symplectic structure on the space of representations, but for us an orientation is unneeded (since we use a nonstandard symplectic structure determined by the choice of an embedding of $J$ in $\bP$) and indeed from our viewpoint choosing an orientation is unnatural and breaks the $\SL_2(\IC)$ symmetry.
\end{rmk}

\subsection{Meromorphic connections on trivial bundles}\label{ssn: mero conns}

The third point of view is as moduli spaces of meromorphic connections on a trivial vector bundle on $\IP^1$. 
Since the setup is similar to \cite{smid} \S2 we will be brief.
Fix $G=\GL_{r}(\IC)$, let $\lt\subset \g=\Lie(G)$ be the diagonal matrices and fix a compact Riemann surface $\Si$.
In general specifying a (symplectic) moduli space of meromorphic connections on 
vector bundles of rank $r$ on $\Si$ (with unramified normal forms) involves specifying an effective divisor $D = \sum k_i(t_i)$ on $\Si$ and at each point 
$t_i$ specifying the corresponding local data.
The local data consists of an ``irregular type'' and a residue orbit.
Here an irregular type is an element
$Q\in \lt\flp z \frp/\lt\flb z \frb$, 
where  $z$ is a local coordinate on $\Si$ vanishing at $t_i$, 
i.e. it is an element
$$Q = \sum_1^{k_i-1} \frac{A_j}{z^j}$$ 
with $A_j\in\lt$, so that $dQ$ has a pole of order at most $k_i$.
Given an irregular type $Q$ consider the group $H\subset G$ centralising $Q$, i.e. consisting of elements $g\in G$ such that 
$gA_jg^{-1}=A_j$ for all $j$.
Then the remaining choice (of residue orbit) is the choice of an adjoint orbit $\cO\subset \lh = \Lie(H)$.
Given such local data $(Q, \cO)$ (at $t_i$) we consider connections on vector bundles over $\Si$ which are locally holomorphically isomorphic at $t_i$ to $$dQ+\La\frac{dz}{z} + holomorphic $$
for some $\La\in \cO$  
(and similarly at the other points of $D$, for other choices of local data).
Note that if $Q=0$ we are just fixing the adjoint orbit of the residue of a logarithmic connection.
We will say such a connection is {\em nonresonant} (at $t_i$) if 
$\ad_\La\in \End(\lh)$ has no eigenvalues which are nonzero integers.
(In the nonresonant case fixing the local data is equivalent to fixing the formal isomorphism class of the connection, but not in general.)

Given the divisor $D$ and local data at each point $t_i$ of $D$ there are several moduli spaces one can consider
(cf. the spaces $\cM$ and their approximations $\cM^*$ in \cite{smid}, and their generalisations, the hyperk\"ahler spaces of stable parabolic meromorphic connections in \cite{wnabh}).
Here we fix $\Si=\IP^1$ and will consider some simple examples 
of moduli spaces of connections with such fixed local data:
we will consider meromorphic connections on trivial vector 
bundles $U\times\IP^1\to \IP^1$ which take the form

\beq \label{eq: conn expr}
\cA= \left(Az+B+T + \sum_{I_\infty} \frac{R_i}{z-t_i}\right)dz
\eeq
as considered in \eqref{eq: conn} (so in particular $A,T$ are semisimple).
Such a connection will be said to be {\em stable} if it admits no proper nontrivial subconnections. (Here we work in the category of connections on trivial bundles so any subconnection should also be on a trivial bundle.)

\begin{lem}
Specifying such a connection
is equivalent to specifying the following data:

A finite dimensional complex vector space $U$,

A finite set $I_\infty$ and distinct $t_i\in \IC$ for each $i\in I_\infty$,

A grading $U = \bigoplus_{I'} V_i$ of $U$ by a finite set $I'$,

A partition $I' = \bigsqcup_{j\in J'} I_j$ of $I'$, 

Distinct complex numbers $a_j\in \IC$ for each $j\in J'$,

Complex numbers $t_i\in \IC$ for each $i\in I'$ such that 
$t_i \ne t_{i'}$ if $i,i'\in I_j$ for some $j$,

A linear map $B_{ij}\in \Hom(W_j,W_i)$ for each $i\ne j\in J'$ where 
$W_j = \bigoplus_{l\in I_j} V_l$,

Elements $R_i\in \End(U)$ for all $i\in I_\infty$. 
\end{lem}
\pf
Given such data set $A = \sum_{J'} a_j \Id_j$ where $\Id_j$ is the idempotent for $W_j$, and set $T = \sum_{I'} t_i \Id_i$ 
where $\Id_i$ is the idempotent for $V_i$, and set $B = \sum B_{ij}$, to obtain all the coefficients of $\cA$. The converse is straightforward.
\epf

Given such data (or the corresponding connection $\cA$) let 
$\cG'$ be the complete $\abs{J'}$-partite graph with nodes $I'$, corresponding to the partition $I' = \bigsqcup_{j\in J'} I_j$.
Specifying $B$ is then equivalent to specifying a representation $\rho$ of the graph $\cG'$ on  the $I'$-graded vector space $U$.
The following is an elementary exercise.

\begin{lem}
The connection $\cA$ is stable if and only if there are no nontrivial proper  subrepresentations
$U'\subset U$ of $\rho$ such that $R_i(U')\subset U'$ for all $i\in I_\infty$.  
\end{lem}

We wish to consider such connections with fixed local data.
Thus, at the simple poles, we fix adjoint orbits  $\cO_i\subset \End(U)$ and require $R_i\in \cO_i$ for each $i\in I_\infty$.
At $z=\infty$ we fix an adjoint orbit $\breve \cO_H\subset \lh$, where
$\lh$ is the Lie algebra of $H= \prod_{i\in I'} \GL(V_i)$ and restrict to connections which are locally holomorphically isomorphic to connections of the form
$$
\left(Az+ T +  \frac{\La}{z} + \cdots \right)dz
$$
near $z=\infty$, for some $\La\in \breve \cO_H$ (i.e. we fix the orbit of the residue of the normal form).

Now, with the data $A,T,\{t_i\}, \breve \cO_H, \{ \cO_i \}$ fixed, consider the space
$$\Conn_{st}(\breve \cO_H, \{ \cO_i \})$$ 
of isomorphism classes of such connections  which are stable.

The aim is to identify
this with a space $\cM_{st}^*$ of stable Weyl algebra module presentations.
First we  set up the required data to define $\cM_{st}^*$.
For each $i\in I_\infty$ set $d_i = \rank(R_i)$ for any $R_i\in \cO_i$ and set $V_i = \IC^{d_i}$. Then define 
$$\breve \cO_i \subset \End(V_i)$$
to be the unique orbit with the following property 
(cf. Appendix \ref{apx: relating orbits}):
if $R_i\in \cO_i$ and  $R_i=Q_iP_i$ for a surjective map 
$P_i:U\to V_i$ and an injective map  $Q_i:V_i\to U$, then 
$-P_iQ_i\in \breve\cO_i$.  
Define $I=I'\sqcup I_\infty$ so that, 
since specifying $\breve \cO_H$ amounts to specifying an orbit 
$\breve \cO_i\subset \End(V_i)$ for each $i\in I'$, we now have 
$\breve \cO_i\subset \End(V_i)$ for all $i\in I$.

Further define
$J = J'\union \{\infty \}$ 
(unless $I_\infty$ is empty, in which case  we set  $J=J'$).
(Here we have identified $J'$ with $\{a_j\st j\in J'\}\subset \IC$.)
Then let $\cG$ be the complete $k$-partite graph with nodes $I$, corresponding to the partition $I=\bigsqcup_J I_j$, where $k=\abs{J}$.
In particular we have a subgraph $\cG'\subset\cG$. 

Thus we now have all the data $\cG, V, \breve \cO, \ba$ necessary  to define 
 $\cM_{st}^*(\cG,\breve\cO)$ as in \S\ref{ssn: wamps}.

\begin{thm}\label{thm: conn isom}
The space $\Conn_{st}(\breve \cO_H, \{ \cO_i \})$ of isomorphism classes of stable connections is isomorphic to 
$\cM_{st}^*(\cG,\breve\cO)$.
\end{thm}
\pf
Recall that $\cM_{st}^* = \mu^{-1}(\breve\cO)^{st}/\wh H$ where 
$\mu^{-1}(\breve\cO)\subset \IM$. 
Thus suppose we have a stable representation $\rho$ of $\cG$ 
in $\mu^{-1}(\breve\cO)^{st}$.
By restricting to $\cG'\subset \cG$ we obtain the coefficient $B$ of $\cA$, and $\rho$ determines the coefficients $R_i=Q_iP_i$ as usual.
By Lemma \ref{lem: inj-surj} and Proposition \ref{prop: orbit relations} 
the stability of $\rho$ and the fact that $\La_i\in \breve \cO_i$
implies $R_i\in \cO_i$ for all $i\in I_\infty$.
Next we need to check that fixing $\La_i\in \breve \cO_i$ for $i\in I'$ corresponds to fixing the residue $\La$ of the normal form of $\cA$ to be in $\breve \cO_H$. In other words (as the notation suggests) we need to check that $\La_i$ is indeed the $i$th component of $\La$.
But $\La$ is computed in appendix \ref{apx: leading term} to be  
$\pi_\lh(QP+[X,B]/2)= \pi_\lh(QP+XB)$ where $\pi_\lh$ is the projection onto  the $\lh$ component, and by definition the $i$th component of this is 
$$ \pi_i(QP+XB) \iota_i   = \pi_i\Xi\circ\Ga \iota_i = -P_i Q_i = \La_i$$
as desired, for $i\in I'$.
Finally one may check the resulting connection is stable, that all connections with the given local data are obtained this way and that $\wh H$ orbits correspond to isomorphism classes, all of which is now straightforward.
(The natural symplectic structures also match up, cf. 
\cite{smid} \S2, and  Lemma \ref{lem: mmap} below.)
\epf

Thus we may also denote this space of stable connections as $\cM^*_{st}(\breve \cO_H, \{ \cO_i \})$.
Combined with Theorem \ref{thm: cM and  Q} it is thus also isomorphic to a twisted quiver variety for a supernova graph.

Note that giving a marking of the orbit $\breve \cO_i$ is the same as giving a {\em special} marking of $\cO_i$. 
Explicitly suppose  $\La_i = -P_iQ_i\in \breve \cO_i$ and 
$R_i=Q_iP_i\in \cO_i$.
Then if $(\xi_{i1},\ldots,\xi_{iw})$ is a 
marking of $\breve\cO_i$ then $(0,-\xi_{i1},\ldots,-\xi_{iw})$ is the corresponding special marking of $\cO_i$. Indeed
$$ 
R_i\Prod_{l=1}^w(R_i+\xi_{il}) = 
(-1)^wQ_i\Prod_{l=1}^w(\La_i-\xi_{il})P_i = 0.$$
Thus the choices involved in identifying $\Conn_{st}(\breve \cO_H, \{ \cO_i \})$ as a twisted quiver variety amount to choosing a marking of $\breve \cO_i$ for $i\in I'=I\setminus I_\infty$ and a 
special marking of each orbit $\cO_i$ (for $i\in I_\infty$).
Given such choices let $\wh \cG$ be the corresponding supernova graph and let 
$\bd,\bla$ be the corresponding data, so that 
\beq\label{eq: conn to qv isom}
\Conn_{st}(\breve \cO_H, \{ \cO_i \})\cong \cQ(\wh\cG,\bla,\bd) \cong 
\cM^*_{st}(\wh\cG,\bla,\bd).
\eeq

One useful input of the viewpoint of meromorphic connections is that we may perform the following scalar shifts.
Choose constants $c_i\in \IC$ for each $i\in I_\infty$ and set 
$c=\sum c_i$. Given orbits $(\breve \cO_H, \{ \cO_i \})$ as above, consider the shifted orbits:
$$
\breve \cO'_H = \breve\cO_H + c\Id_U,\qquad 
\cO'_i = \cO_i + c_i\Id_{U}.$$

\begin{lem}\label{lem: conn isoms}
The moduli spaces $\Conn_{st}(\breve \cO_H, \{ \cO_i \})$ and 
$\Conn_{st}(\breve \cO'_H, \{ \cO'_i \})$ are isomorphic.
\end{lem}
\pf
This is almost immediate from the expression \eqref{eq: conn expr} for the connections $\cA$; the map is given by replacing each $R_i$ by $R_i+c_i$ and leaving $B$ unchanged. The orbit $\breve\cO_H$ is shifted as stated since  $\La=\pi_\lh(QP+XB)$ and $QP = \sum R_i$.
\epf

It is clear from the projected equations \eqref{eq: inter eqs} that this operation relates the corresponding isomonodromy connections.

Finally we can deduce isomorphisms corresponding to the 
reflections at the nodes $i\in I$ of the core. Suppose that the support of $\bd$ intersects at least two parts $I_j\subset I$ of the core nodes (i.e. that we are not in a trivial case with just one part).

\begin{cor}
If $i\in I$ and $\la_i\neq 0$ then the space
$\cM^*_{st}(\wh \cG, \bla, \bd)$ is isomorphic to the space
 $\cM^*_{st}(\wh \cG, r_i(\bla), s_i(\bd))$
(where $r_i,s_i$ are the reflections of \S\ref{ssn: km}).
\end{cor}
\pf
We may suppose $i\in I_\infty$, since if $i\in I_j$ and $a_j\ne \infty$ then we may conjugate by a symplectic $\SL_2(\IC)$ transformation moving $a_j$ to $\infty$.
Then in terms of connections the idea is to first
change the marking of the orbit $\cO_i$, swapping the 
order of the first two eigenvalues.
The resulting 
marking will not be special, so we then 
perform a scalar shift to return to a special marking. 
This gives the desired reflection, as we will now verify in detail. 
Suppose $(\xi_{j1},\xi_{j2},\ldots)$ is the marking of 
$\breve \cO_j$ for any $j\in I$, so that 
$$(0,-\xi_{i1},-\xi_{i2},\ldots)$$ 
is the marking of $\cO_i$. Let $\bla $ denote the corresponding set of parameters (so for example 
$\la_{j1} = \xi_{j1}, \la_{j2} = \xi_{j2} - \xi_{j1}$ on the first two nodes of the $j$th leg).
After reordering, the marking of $\cO_i$ is changed to  
$(-\xi_{i1},0, -\xi_{i2},\ldots)$.
Since this is not special we perform the scalar shift by 
$c_i=c= \xi_{i1}$, so $\cO_i$ is replaced by $\cO'_i = \cO_i+c$ which has the special marking $(0,c, c -\xi_{i2},c-\xi_{i3},\ldots)$,
and $\breve \cO_j$ is replaced by 
$\breve\cO'_j=\breve \cO_j + c$, which has marking 
$(c+\xi_{j1},c+\xi_{j2},\ldots)$, for all $j\in I'$.
The isomorphism  $\Conn_{st}(\breve \cO_H, \{ \cO_j \})\cong\Conn_{st}(\breve \cO'_H, \{ \cO'_j \})$ of Lemma \ref{lem: conn isoms}
then yields the desired isomorphism $\cQ(\wh \cG, \bla, \bd)\cong\cQ(\wh \cG, r_i(\bla), s_i(\bd))$, once we pass from connections to 
twisted quiver varieties using the chosen markings.
Indeed computing the parameters corresponding to the new markings yields 
$$\la_{i1}\mapsto -c = -\xi_{i1} = -\la_{i1},$$
$$\la_{i2}\mapsto (\xi_{i2}-c)-(-c) = \xi_{i2} = \la_{i2}+\la_{i1},$$
$$\la_{j1}\mapsto \la_{j1} + c = \la_{j1}+\la_{i1}$$
for all $j\in I'$, with all the other components unchanged; These are the components of $r_i(\bla)$.
Considering the dimensions, $d_{i1}$ is the only component of $\bd$ which is changed, and the new value $d'_{i1}$ may be computed in terms of $\bd$ as follows:
$$d_{i1} = \rank(R_i) = \dim(U) -\dim\ker(R_i)$$
$$d'_{i1} = \rank(R_i+c) = \dim(U) -\dim\ker(R_i+c)$$
$$d_{i2} =  \dim(U) - \dim\ker(R_i) - \dim\ker(R_i+c)$$
so that $d'_{i1} = \dim(U) + d_{i2} - d_{i1}$, which is the
corresponding component of $s_i(\bd)$, given that 
$\dim(U) = \sum_{j\in I'} d_{j1}$.
\epf

Note that these scalar shifts generalise those in \cite{k2p} (e.g. bottom of p.185, in the case with just one simple pole, $\abs{I_\infty}=1$), 
whose origins lie in the twisted Fourier--Laplace transform of \cite{BJL81}.  
They were used in \cite{k2p} to better understand the Okamoto symmetries of Painlev\'e VI (in particular the action on linear monodromy data was deduced using this viewpoint
in \cite{k2p} Corollary 35---see also \cite{pecr} Remark 4).
Another approach is possible using the middle convolution operation 
\cite{Katz-rls}, but this may be derived from Fourier--Laplace (cf. \cite{Katz-rls} \S2.10, \cite{nlsl} Diagram 1, \cite{yamakawa-mc+hd}). 
The reflection isomorphisms constructed here 
are analogous to those for Nakajima quiver varieties 
(cf. \cite{quad} Theorem 1 and references therein), but it is not clear if they are actually equivalent;
the moduli theoretic approach here enables us to see the isomonodromy systems are preserved, and shows that they will extend to the full \hk wild nonabelian Hodge moduli spaces (after extending \cite{szabo-nahm}).

\section{Additive irregular Deligne--Simpson problems} \label{sn: +iDSP}

Suppose we fix $\Si=\IP^1$ to be the Riemann sphere and have an effective divisor $D=\sum k_i(t_i)$ on $\Si$.
If we fix local data (as in \S\ref{ssn: mero conns}) consisting of an irregular type and residue orbit at each point $t_i$, then we may consider the moduli space $\Conn_{st}$ of 
stable meromorphic connections on the trivial bundle over $\Si$ with the given local data at each point of $D$.
The (unramified) {\em additive irregular Deligne--Simpson problem}
is  then to characterise the local data for which $\Conn_{st}$
is nonempty (i.e. for which there exists such stable connections on the trivial bundle).
This is the natural extension of the usual additive Deligne--Simpson problem to the irregular case.
We will solve some cases of this here using \cite{CB-mmap} (these results appeared in the preprint \cite{rsode}). 
In the Fuchsian case (all irregular types zero) Crawley--Boevey 
\cite{CB-additiveDS} established a precise criterion in terms of roots of an associated  Kac--Moody root system for a star-shaped graph.
This was proved by identifying the space of stable Fuchsian systems 
with the stable points of a quiver variety and then 
using earlier results \cite{CB-mmap} characterising exactly when certain quiver varieties have stable points.
In \S\ref{sn: reductions}  we have identified some 
more general spaces of meromorphic connections 
with (twisted) quiver varieties and so can again use 
\cite{CB-mmap} to give a precise criterion for the existence of stable points, as follows.
The setup is the same as in \S\ref{ssn: mero conns}:

Choose a complex vector space $U=\IC^n$, 
distinct points $t_1,\ldots ,t_m\in \IC$  and 
diagonal matrices $A,T\in \End(U)$.
This determines the  irregular type 
$Q= Az^2/2 + Tz$ at $z=\infty$.
Consider connections on $U\times \IP^1\to \IP^1$ with Fuchsian singularities at each point $t_i$ and an irregular singularity at $z=\infty$ with irregular type $Q$.
Fix adjoint orbits $\cO_1,\ldots ,\cO_m \subset \End(U)$ and 
$\breve \cO_H\subset \lh$ where $\lh\subset \End(U)$ is the set of matrices that commute with both $A$ and $T$.
The problem is to decide when there are stable connections of the form
$$\cA= \left(Az+B+T + \sum_1^m \frac{R_i}{z-t_i}\right)dz$$ 
with $R_i\in \cO_i$ and $B\in \im(\ad_A)\subset \End(U)$ so that 
$\cA$ is locally holomorphically isomorphic to a connection of the form 
$$\left(Az+T + \frac{\La}{z} + holomorphic \right)dz$$
near $z=\infty$ for some $\La\in \breve \cO_H$.

To answer this, choose markings as in \S\ref{ssn: mero conns}
so a supernova graph $\wh \cG$ and data $\bla,\bd$ are determined,  and there is an isomorphism 
$$\Conn_{st}(\breve \cO_H, \{ \cO_i \})\cong \cQ(\wh\cG,\bla,\bd)$$
as in \eqref{eq: conn to qv isom},
from the space of such stable connections to 
the corresponding twisted quiver variety.
In particular $\wh\cG$ determines a Kac--Moody root system as described in \S\ref{ssn: km}.

\begin{cor}
There are stable connections $\cA$ with the given local data as above, if and only if:

1) $\bd$ is a positive root,

2) $\bla\cdot\bd  =0$, and

3) If $\bd=\bd_1+\bd_2+\cdots$  is a nontrivial sum of positive roots such that 
$\bla\cdot\bd_1=\bla\cdot\bd_2 =\cdots =0$, 
then $\De(\bd)>\De(\bd_1)+\De(\bd_2)+\cdots$, where 
$\De(\bd) = 2-(\bd,\bd)$.
\end{cor}
\pf
Since Theorems \ref{thm: cM and Q} and \ref{thm: conn isom} show that such stable connections correspond to stable points of 
a twisted quiver variety, and Crawley--Boevey \cite{CB-mmap} Theorem 1.2 has shown that these criteria characterise quiver varieties having stable points, 
it just remains to show that our twisted quiver variety is isomorphic to an (untwisted) quiver variety (or equivalently in the language of \cite{CB-mmap}, that points of $\cQ(\wh\cG,\bla,\bd)$ correspond to simple representations of the deformed preprojective algebra). 
This boils down to proving the following lemma.
Let $\cG\subset \wh \cG$ be the core of the graph $\wh \cG$, with 
nodes $I\subset \wh I$, and choose an orientation of $\cG$. 
This gives a map $\cG\hookrightarrow\bar{\cG}$ to the set of oriented edges $\bar{\cG}$ of $\cG$ (cf. \S\ref{ssn: reps of graphs}),   
so each edge $e\in\cG$ has a 
well-defined head $h(e)\in I$ and tail $t(e)\in I$.
Let $\IM$ be the symplectic vector space of 
\S\ref{sn: sympl}, which we identify with $\Rep(\cG,V)$  as usual (see Proposition \ref{prop: M as graph reps}), with $V=\bigoplus_I V_i,  V_i=\IC^{d_i}$.
Recall from Proposition \ref{prop: mmap La_i} that
 $\wh H = \Prod_I \GL(V_i)$ acts on $\IM$ with moment map $\Ga\mapsto (\La_i).$
On the other hand in the theory of quiver varieties one identifies 
\beq\label{eq: rep as cotbdle}
\Rep(\cG,V)  = T^*\bigoplus_{e\in\cG} \Hom(V_{t(e)},V_{h(e)})
\eeq
with the cotangent bundle of the space of maps along the edges following the given orientation.
In this symplectic structure a moment map for the action of $\wh H$
has $\End(V_i)$ component taking $\rho\in \Rep(\cG,V)$ to
$$\sum_{e\in \bar{\cG}, t(e)=i} \eps(e)\rho(\bar {e}) \rho(e)
$$
where $\bar{e}$ is the edge $e$ with the opposite orientation, and 
$\eps(e)=1$ if $e\in \cG$ and $\eps(e)=-1$ if $e\in \bar{\cG}\setminus \cG$.

\begin{lem}
The space $\IM$ and the cotangent bundle \eqref{eq: rep as cotbdle}
are isomorphic as Hamiltonian $\wh H$-spaces.
\end{lem}
\pf
Recall that $\La_i=-P_iQ_i$ is the $\End(V_i)$ component of 
$\Xi\Ga\in \End(V)$, so that 
if $\rho\in \Rep(\cG,V)$ is the representation corresponding to $\Ga$ then
$$\La_i=\sum_{e\in \bar{\cG}, t(e)=i} \phi(\bar{e})\rho(\bar {e}) \rho(e)$$
where $\phi(e):=\phi_{jj'}$ if $e$ is the edge from 
$i\in I_j$ to $i'\in I_{j'}$.
Thus if we define a linear map $\IM\to \IM$ taking 
a representation $\rho$ to the representation $\rho'$ defined by
$\rho'(e) = -\phi(e)\rho(e)/\eps(e)$ if $e\in\cG$ and 
$\rho'(e) = \rho(e)$ if $e\in \bar{\cG}\setminus \cG$, then
$$\phi(\bar{e})\rho(\bar {e}) \rho(e) = \eps(e)\rho'(\bar {e}) \rho'(e)$$
for all $e\in \bar\cG$, so the moment maps are intertwined as desired.
Moreover, from \S\ref{sn: sympl}, 
the symplectic structure on $\IM$ is given by 
$$
\frac{1}{2}\tr( d\Xi\wedge d\Ga) = 
\frac{1}{2}\sum_{e\in\bar\cG} \phi(\bar{e})\Tr(d\rho(\bar{e})\wedge 
d\rho(e)) = 
\frac{1}{2}\sum_{e\in\bar\cG} \eps(e)\Tr(d\rho'(\bar{e})\wedge 
d\rho'(e))$$
$$= 
\sum_{e\in\cG} \Tr(d\rho'(\bar{e})\wedge d\rho'(e))$$
which is the standard symplectic structure on the cotangent bundle \eqref{eq: rep as cotbdle}.
\epf
Thus $\cQ(\wh\cG,\bla,\bd) =\wh \IM\spq_{\bla} \wh\IG$ is isomorphic to a  quiver variety (since $\wh \IM=\IM \times \IL$ and we already used the cotangent symplectic structure on the legs $\IL$), and so we may apply \cite{CB-mmap} Theorem 1.2.
\epf

\section{Example reduced systems}

Thus in summary we can now attach an isomonodromy system 
\beq\label{eq: imd syst}
\cM^*_{st}\times \IB\to \IB\eeq
to the data $\wh \cG, \bla,\bd,\ba$ as in 
Theorem \ref{thm: intro thm} of the introduction.
If nonempty then $\dim_{\IC}(\cM^*_{st}) = 2-(\bd,\bd)$ by Proposition
\ref{prop: qv constr}, and a precise criterion for nonemptiness is given in \S\ref{sn: +iDSP}.
Thus if one chooses some local (or birational) coordinates on 
$\cM^*_{st}$ then the isomonodromy connection on 
\eqref{eq: imd syst} amounts to a system of nonlinear equations of order $2-(\bd,\bd)$.
If the core $\cG$ of $\wh \cG$ is a complete $k$-partite graph then, by moving $\ba$ around in $\bP$ using the symplectic transformations, 
this isomonodromy system controls isomonodromic deformations of linear systems $\cA$ on vector bundles with in general $k+1$ different ranks (depending on which, if any, of the $k$ elements of $\ba$ is moved to $\infty$).
In other words there are $k+1$ ways to read the graph in terms of connections (before considering the reflections).

Since $\cM^*_{st}$ is symplectic the simplest nontrivial case is dimension $2$, when $(\bd,\bd) =0$ (i.e. $\bd$ is a null root).
These cases correspond to affine Dynkin diagrams, and in turn to the (known) second order Painlev\'e systems.
The dimension vectors $\bd$ for the supernova cases are as in the following diagram ($\bd$ is the minimal imaginary root).
(One may also consider the type $E$ affine Dynkin graphs but these are only of interest for Painlev\'e difference equations, cf. \cite{quad}; they have no nontrivial isomonodromic deformations.)

\begin{figure}[ht] 
	\centering
	\input{p456.no.pstex_t}
	\caption{} \label{fig: painleve 4-6}
\end{figure}

These affine Dynkin cases are special since $\bd$ is preserved by the reflections, and only the parameters $\bla$ are moved. 
In each case one still has the choice of the parameters $\ba$, and we will now illustrate in the case of Painlev\'e VI how to read the graphs to give various equivalent Lax pairs for these Painlev\'e systems.

\subsection{Painlev\'e VI}\label{ssn: PVI}

Here it is simplest to view the affine $D_4$ diagram as a supernova graph as in Figure \ref{fig: p6} with just four nodes in the core, say $I_1$ consists of three of the feet and $I_2$ consists of the central node.
The three ways to read this  graph are then as follows.

\begin{figure}[ht] 
	\centering
	\input{p6.systems.pstex_t}
	\caption{} \label{fig: p6}
\end{figure}

1) $a_1 = \infty, a_2=0$: This is the standard Lax pair as Fuchsian systems 
$$\cA = \left(\sum_1^3 \frac{R_i}{z-t_i}\right)dz$$
on rank two bundles with four poles on the Riemann sphere.
(The three finite poles correspond to $\abs{I_1}=3$ and the rank two is the dimension on the central node.)

2) $a_1 = 0, a_2=\infty$: This is Harnad's dual Lax pair \cite{Harn94}; systems
$$\cA = \left(T_0 + \frac{R}{z}\right)dz$$
on rank three bundles with a Fuchsian singularity and an irregular singularity on the Riemann sphere.
(The single finite Fuchsian singularity corresponds to $\abs{I_2}=1$ and the rank three is the sum of the dimensions of the nodes in $I_1$.)

3) $a_1 = 0, a_2=1$: This is the generic reading of this graph; as systems of the form
$$\cA = \left(Az + B + T \right)dz$$
on rank five bundles with just an irregular singularity at $\infty$ 
on the Riemann sphere (and no others poles).
(The rank five is the sum of the dimensions of the nodes in 
$I_1\union I_2$.)
Explicitly the form of $\cA$ specified by the graph is as follows:
$
A = \diag(0,0,0,1,1), 
T = \diag(t_1,t_2,t_3,t_4,t_4)$
with $\abs{\{t_1,t_2,t_3\}}=3$ and 
$ B=  \bsmx
 0 &  B_{12}  \\ 
 B_{21} & 0 
\esmx\in \End(\IC^3\oplus \IC^2)$. 
Fixing the orbit of the residue $\La$ of the normal form at $z=\infty$
then corresponds to fixing the orbit of $B_{21}B_{12}\in \End(\IC^2)$ and the diagonal part of $B_{12}B_{21}\in \End(\IC^3)$.
The nonlinear equations are as in \eqref{eq: master}
with $ X=  \bsmx
 0 &  -B_{12}  \\ 
 B_{21} & 0 
\esmx.$

\begin{figure}[ht] 
	\centering
	\input{p6.extsystems.pstex_t}
	\caption{} \label{fig: ep6}
\end{figure}

If this graph is embedded in a larger graph, then performing reflections will give infinitely many other Lax pairs. 
For example the graph in Figure \ref{fig: p6} is the same as the graph on the left in Figure \ref{fig: ep6} (obtained by adding a node with dimension zero). 
Then if we perform three reflections working up the long leg from the toe, the diagram on the right of Figure  \ref{fig: ep6} is obtained. 
This graph may be read (as in 1) above) as Fuchsian systems on rank three bundles with four poles on $\IP^1$ such that the three residues at finite points are rank one matrices.
This is the alternative Fuchsian Lax pair for Painlev\'e VI 
used in \cite{k2p} to obtain new algebraic solutions of Painlev\'e VI
(generalising  that, for special values of the parameters, in \cite{Dub95long} Remark 3.9 and Appendix E).

\subsection{Painlev\'e IV} Similarly the triangle $\cG(1,1,1)$ may be read in various ways. If the dimension vector is $(1,1,1)$  one reading is as connections with one pole of order $3$ and one pole of order $1$ on a rank two bundle (this is the standard Lax pair for Painlev\'e IV). 
The generic reading of this graph is as a space of connections on a rank three bundle with a pole of order $3$ and no others. 
(This generic reading appears in an explicit form in \cite{JKT07}.)

\subsection{An infinite Weyl group orbit} 
Consider the case $A_2^{++}$ of the triangle with a leg of length one attached, as on the right of Figure \ref{fig: higherpainleve 6-4} (but we will use different dimension vectors here). Then \cite{ogg84, FKN08} the index two `rotation subgroup' of the corresponding 
Kac--Moody Weyl group is $\PSL_2(E)$ where $E=\IZ[\omega]$ is the ring of Eisenstein 
integers (where $\omega^3=1$).
For example if we label the nodes of $A_2^{++}$ as $1,2,3,4$ (with the $2$ in the middle, the $1$ at the foot and $2,3,4$ on the triangle)
and take dimension vector $\bd=(1,2,2,1)$ and generic parameters, then the corresponding moduli space $\cM_{st}^*(\bla,\bd)$ has complex dimension $2$. 
Indeed performing the reflections $s_1s_2s_3$ yields dimension vector $(0,1,1,1)$ so the variety is isomorphic to the space appearing in the case of Painlev\'e IV (which in one reading is thus also isomorphic to a space of connections on a rank $3$ bundle with $2$ poles of order $1$ and $3$).
On the other hand the Weyl group element (see  \cite{FKN08} 4.20): 
$$w=s_1s_4s_1s_2s_4s_1s_3s_1$$
has infinite order, realising the same space as a space of connections on bundles of arbitrarily high rank---indeed for $n\ge 1$ the space with dimension vector $w^n(1,2,2,1)$ may be read as a space of connections as above on bundles of rank $n^2+(n-1)+(n-2)^2$.
This thus gives an infinite number of Lax pairs for Painlev\'e IV.
(In general we view each such realisation as a ``representation'' of the abstract Painlev\'e system, or more generally of the corresponding nonabelian Hodge structure.) 

\subsection{Higher Painlev\'e systems}\label{sn: hps}

Recall that 
the next simplest class of Kac--Moody algebras after the affine ones are the hyperbolic Kac--Moody algebras, which by definition, have Dynkin diagrams such that any proper subgraph is either a finite or affine Dynkin diagram.
Since we know the lowest dimensional moduli spaces (those of dimension two) are related to affine Dynkin diagrams it is natural to look amongst the hyperbolic graphs to find the next simplest cases.
This was taken up in \cite{rsode} p.11-12, where a computer search was done for dimension vectors for hyperbolic supernova graphs such that $2-(\bd,\bd)=4$, i.e., so that the corresponding moduli space has dimension $4$ (or that the corresponding isomonodromy system has fourth order). 
By looking at these examples one sees there are several natural families of examples, which arise by ``doubly extending'' a finite Dynkin diagram (i.e. by adding another node next to the extending node of an affine diagram), and that these examples generalise to give examples with dimension $2-(\bd,\bd)=2n$ for any $n$. Within the class of (simply-laced) supernova graphs this includes the following three cases.

\begin{figure}[ht] 	\centering
	\input{higherpainleve2.pstex_t}	\caption{Dynkin diagrams for higher Painlev\'e systems \hPVIn, \hPVn, \hPIVn.} \label{fig: higherpainleve 6-4}
\end{figure}

\begin{prop}
Let $\wh\cG$ be one of the  supernova graphs of Figure 
\ref{fig: higherpainleve 6-4} with the given dimension vector $\bd$, and let $\bla$ be some generic parameters. Then  
$\cM^*_{st}(\wh\cG,\bla,\bd)$ has dimension $2n$. 
\end{prop}
\pf
The dimension vector has the form $\bd=n\de + e$ where $\de$ is the minimal imaginary root for the affine subdiagram (removing the extending node with dimension $1$), and $e$ is supported on the extending node, and so it is easy to compute $(\bd,\bd) = 2-2n$.
\epf

Note that for $n=1$ these are in general isomorphic to 
the usual Painlev\'e systems, since we may perform the reflection at the foot (the extending node) 
so the resulting dimension vector is supported on the affine subdiagram.

We will now write down explicitly the simplest reading of each of these graphs (for generic parameters $\bla$). 
This enables us to spot the pattern and thereby describe Lax pairs for other higher Painlev\'e systems
(we will ignore the ramified cases for simplicity, even though they present little difficulty except in notation). 
Recall that the moduli spaces are determined by fixing the local data at each pole, consisting of the irregular type and residue orbit.

Let $U=\IC^{2n}$, and $G=\GL(U)$ with $\lt\subset \g=\Lie(G)$ the diagonal matrices and fix $\Si=\IP^1$.

\subsection{Higher Painlev\'e VI}
Here $\Si$ has four marked points and all irregular types zero.
The local data consists of four semisimple orbits   
$\cO_i\subset \GL(U)$ for $i=1,2,3,4$, such that 
$\cO_1$ has three eigenvalues with multiplicities $n,n-1,1$ and each of  the other classes $\cO_2,\cO_3,\cO_4$ has just two eigenvalues, each of multiplicity $n$. 
The resulting space  $(\cO_1\times\cdots\times\cO_4)\spq G$ has dimension  $(2n^2+2n-2)+3(2n^2) - 2(4n^2-1)=2n$.

\subsection{Higher Painlev\'e V}
Here $\Si$ has three marked points. 
Two of the marked points are Fuchsian with local data 
given by semisimple orbits 
$\cO_i\subset \GL(U)$ for $i=1,2$, such that 
$\cO_1$ has three eigenvalues with multiplicities $n,n-1,1$ and $\cO_2$ has just two eigenvalues, each of multiplicity $n$. 
At the third singularity the irregular type is $Q=A_1/z$ with $A_1\in \lt$ having two eigenvalues each of multiplicity $n$.
The residue orbit at this irregular singularity 
is specified by two orbits in $\gl_n(\IC)$; we take them both to be scalar orbits (i.e. zero-dimensional).
Thus the resulting space is of the form
$$(\cO_1\times \cO_2 )\spqa{\breve\cO_H} \GL_n(\IC)\times \GL_n(\IC)$$

\noindent
where $\breve\cO_H$ is zero-dimensional and so, 
since $\cO_1$ has dimension $2n^2+2n-2$ and $\cO_2$ has 
dimension $2n^2$, the resulting space has dimension 
$4n^2+2n-2-2(2n^2-1) = 2n$.

\subsection{Higher Painlev\'e IV}
Here $\Si$ has two marked points. 
One of the marked points is Fuchsian with
local data given by a semisimple orbit  
$\cO\subset \GL(U)$ with three eigenvalues with multiplicities $n,n-1,1$.
The other singularity has irregular type 
$Q=A_2/z^2 +A_1/z$ with $A_2\in \lt$ having two eigenvalues each of multiplicity $n$ (and $A_1\in \lt$ is any element whose centraliser in $G$ contains that of $A_2$ so that the centraliser $H$ of $Q$ is 
$\GL_n(\IC)\times \GL_n(\IC)$).
The orbit of the residue of the normal form at this irregular singularity is specified by two orbits in $\gl_n(\IC)$; we take them both to be scalar orbits (i.e. zero-dimensional).

\begin{rmk}
Note that none of the readings of \hPIV\  fall within the scope of the JMU system for $n\ge 2$, whereas they all fall within the scope of this article (this was one of our original motivations). 
On the other hand, using the above readings,  the \hPVI\  system 
is a special case of the Schlesinger system 
and \hPV\   is a special case of the JMMS system.
(Note that the fourth order member of the \hPVI\ family has recently been written in explicit coordinates by Sakai \cite{sakai-4th} p.20. The corresponding graph appears in both \cite{rsode} p.12 and \cite{oshima} p.21)
\end{rmk}

\subsection{Higher Painlev\'e III}
Here $\Si$ has two marked points, say at $z=0,\infty$. 
Both marked points have irregular type of the same form  
$Q_0=A_0/z, Q_\infty=A_\infty z$ 
where $A_0,A_\infty\in \lt$ each have two eigenvalues of multiplicity $n$.
At each pole the residue orbit 
is specified by two orbits in $\gl_n(\IC)$; 
we take them both to be scalar orbits (i.e. zero-dimensional) at $0$, and at $\infty$ we take one scalar orbit plus a semisimple orbit with 
two eigenvalues of multiplicities $1,n-1$ respectively.

\subsection{Higher Painlev\'e II}
Here $\Si$ has just one marked point. 
The irregular type $Q$ is of the  form  
$A_3/z^3 +A_2/z^2+A_1/z$ where $A_3\in \lt$ has two eigenvalues each of multiplicity 
$n$
(and $A_1,A_2\in \lt$ are any elements whose centraliser in $G$ contains that of $A_3$, so that the centraliser $H$ of $Q$ is 
$\GL_n(\IC)\times \GL_n(\IC)$).
The residue orbit is specified by orbits $\breve\cO_1,\breve\cO_2\subset \gl_n(\IC)$; we take $\breve\cO_1$ to be a scalar orbit and $\breve\cO_2$ to be a semisimple orbit with two eigenvalues of multiplicities $1,n-1$ respectively.

(The actual values of all the residue eigenvalues are chosen generically, 
subject to the constraint that the sum of the traces is zero).
Similarly there are higher versions of the Painlev\'e systems with ramified normal forms (this will be discussed in detail elsewhere).

\begin{rmk}
Note that in all the above cases, except that of Painlev\'e III, the underlying two-dimensional Painlev\'e moduli space $\cM^*_{st}$ is isomorphic to an ALE hyperk\"ahler four manifold (\cite{Kron.ale}) and it turns out that the (fibre of the) 
corresponding higher Painlev\'e system is diffeomorphic to the Hilbert scheme of $n$-points on the corresponding ALE space (\cite{nakaj-email08, nakaj-sALE})\footnote{thus the ``h'' in \hPVI\  etc. might also stand for Hilbert, as well as higher or hyperbolic.}.
Presumably this also holds in the case of Painlev\'e III 
(in this case the underlying complex surface is a $D_2$ ALF space).
Further, presumably 
it should hold also for the full moduli spaces $\cM$ (as studied in \cite{wnabh}) 
and not just their approximations $\cM^*$ studied here.
Said differently (changing complex structure) this suggests the following conjecture: If $\cM_H$ is a two-dimensional (meromorphic) Hitchin system then the Hilbert scheme of $n$-points on $\cM_H$ is again diffeomorphic to a meromorphic Hitchin system, at least if the parameters involved are sufficiently generic.
(Here we mean diffeomorphic on the nose, not just modulo a birational 
map.)\footnote{
Note that one can consider the cotangent bundle of an elliptic curve
to be a trivial two-dimensional Hitchin system (the hyperk\"ahler metric is flat).
In this case analogous higher Hitchin systems are described  in \cite{donagi-witten},\cite{Donagi} p.24 in relation to Seiberg--Witten theory.}
Note that, using the graphs and then extrapolating as above, we are thus  able to predict exactly which higher dimensional Hitchin systems to look at.  
\end{rmk}

\appendix

 \section{Loop algebras and Adler--Kostant--Symes} \label{apx: AKS}

Let $z$ be a standard coordinate on $\IP^1$ and 
let $D=\{\infty\}\union\{t_i\st i\in I_\infty\}\subset \IP^1$.
Write $z_i = z-t_i$ and $z_\infty  = 1/z$, so we have a preferred local coordinate vanishing at each point of $D$.

Let $\g= \gl_n(\IC)$ where $n=\dim(U_\infty)$ and consider the Lie algebra
$\g(*D)$ of rational maps $\IP^1\to \g$ with poles of arbitrary order on $D$ (and nowhere else).
Thus, subtracting the principal parts at each $t_i$ to leave a $\g$-valued polynomial, yields a vector space isomorphism
$$\g(*D) \cong \cL^-:=\cL^-_\infty \oplus \bigoplus \cL^-_i$$
where $\cL^-_i = z_i^{-1}\g[z_i^{-1}]$ and 
$\cL^-_\infty = \g[z_\infty^{-1}]=\g[z]$.
Of course $\cL^-$ has its own (product) Lie algebra structure, and as such it is ``half'' of the larger Lie algebra:
$$\cL = \cL_\infty \oplus \bigoplus \cL_i$$
where $\cL_i = \g(\!(z_i)\!)$ and $\cL_\infty = \g(\!(z_\infty)\!)$, i.e. there is a vector space decomposition $\cL = \cL^+\oplus \cL^-$
into subalgebras where $\cL^+ = \cL^+_\infty \oplus \bigoplus \cL^+_i$
with $\cL^+_i = \g\flb z_i\frb $ and 
$\cL^+_\infty = z_\infty\g\flb z_\infty\frb$.
(Beware that the convention for $+/-$ is not uniform in the literature.)
Now each $\cL_i$ (including $i=\infty$) 
has a nondegenerate invariant bilinear form given by 
$$(X,Y) = \res_i\tr(XY dz)$$
and together these determine a bilinear form on $\cL$.
This identifies $\cL^-$ with the dual of the Lie algebra $\cL^+$ and so $\cL^-\cong \g(*D)$ inherits a Poisson structure.
The symplectic leaves are finite dimensional and are the coadjoint orbits of the group $G^+=B_\infty\times\prod G\flb z_i \frb$ corresponding to $\cL^+$, where $B_\infty\subset G\flb z_\infty \frb$ is the kernel of the map evaluating at $z_\infty=0$. 

\begin{lem}\label{lem: mmap}
The map 
$$\IM \to \g(*D);\qquad \Ga\mapsto \cA = Az+B+T+Q(z-C)^{-1}P$$
is a Poisson map, indeed it is the moment map for an action of $G^+$ on $\IM$.
\end{lem}
\begin{lem}\label{lem: AKS}(Adler--Kostant--Symes, see e.g. 
\cite{adler-vm-euc} Theorem 3.1.)
The restriction to  $\g(*D)$ of any pair of $\Ad$-invariant functions on $\cL$ are Poisson commuting (and hence so is their pull-back to $\IM$).
\end{lem}

\begin{lem}\label{lem: hamvfs}
Fix $\cA\in \g(*D)$ and choose an element $\cB\in \cL_i$ for some $i$ (or $i=\infty$).
Suppose that $[\cB, \cA_-]_-=0\in \cL_i^-$ (where $\cA_-\in \cL_i^-$ is the projection of $\cA$).
Then the Hamiltonian vector field on $\g(*D)$ at $\cA\in\g(*D)$ corresponding to the one-form $\res_i\tr (\cB\, d\cA)dz$ on $\g(*D)$ is
$$[\cA,\cB_-]_{\g(*D)}$$ 
where $\cB_-\in \cL_i^-\subset \g(*D)$ is the projection of $\cB$.
\end{lem}

\pfms (of Lemma \ref{lem: mmap}).
Let $\cO_B$ denote the coadjoint orbit of $B_\infty$ through
$(Az+T)\in \cL^-_\infty$. 
Then we have $\cO_B = \{(Az+T+B)\st B\in \Im(\ad_A)\subset \g\}$.
The tangents to $\cO_B$  are of the form $\dot B=[A,Y]$ for $Y\in \g$ and without loss of generality we may assume $Y\in \Im\ad_A$.
The (KKS) symplectic structure on $\cO_B$ is, if we have a second tangent 
$B'=[A,Z]$
$$\omega([A,Y],[A,Z]) = \res_\infty \tr(Az+T+B)[Y/z,Z/z]dz = 
-\tr A[Y,Z] = \tr Y B' $$
so that $\omega = \frac{1}{2}\tr dX\wedge dB $ where $dX = \ad_A^{-1}(dB)$, which appears in the expression \eqref{eq: omega QPXB} for the symplectic structure on $\IM$.
The other term $\Tr(dQ\wedge dP)$ in \eqref{eq: omega QPXB} is just the symplectic form on the cotangent bundle $T^*\Hom(W_\infty, U_\infty)$
which, breaking up $W_\infty$,
decomposes as $\sum_{i\in I_\infty}\Tr(dQ_i\wedge dP_i)$ on 
$\bigoplus_{i\in I_\infty}T^*\Hom(V_i,U_\infty)$.
Now $G_i^+=G\flb z_i\frb$ acts on $T^*\Hom(V_i,U_\infty)$ via the projection $G_i^+\to G$ (evaluating at $t_i$) and the action
$$g(Q_i,P_i) = (gQ_i,P_ig^{-1})$$
of $G$, which has moment map $(Q_i,P_i)\mapsto Q_iP_i/(z-t_i) \in \cL_i^-$.
Repeating for each $i$ yields the result, noting that $Q(z-C)^{-1}P =  \sum Q_iP_i/(z-t_i)$.
\epfms

\pfms (of Lemma \ref{lem: hamvfs}).
The one-form at $\cA$ on $\g(*D)\cong \cL^-$ represents an element $X$ of $\cL^+$.  
The minor subtlety here is that in the expression 
$\res_i\tr (\cB\, d\cA)dz$ we are taking the full Laurent expansion of $d\cA$ at $z_i=0$, so by definition $X\in \cL^+$ is such that 
$$\res_i\tr (\cB\, d\!\cA)dz = (X,d\cA)$$
where on the right we take the various principal parts of $d\cA$.
Thus if $\cB = \cB_+ +\cB_-$ then by the residue theorem
$X = \cB_+  - \sum_{j\ne i}\pi_j(\cB_-)\in \cL^+$ 
(where we include $\infty$ in the sum if $i\ne \infty$, and $\pi_j:\cL_i^-\to \cL^+_j$ is the map taking the Taylor expansion at $j$).
Thus using the Poisson structure on $\g(*D)$
the one-form yields the Hamiltonian vector $\ad^*_X\cA$, which is the $\cL^-$ component of $[X,\cA]_\cL\in \cL$.
This equals $[\cA,\cB_-]_{\g(*D)}$ since (as one may readily check) it has the same component in each $\cL_j^-$.
\epfms

\section{Harnad duality}\label{apx: Harnad}

By definition the Harnad dual of the rational differential operator
\beq \label{eq: harnadop1}
\frac{d}{dz}- \left(T_0 + Q(z-T_\infty)^{-1}P\right)
\eeq
is the differential operator
$$\frac{d}{dz} + \left(T_\infty + P(z-T_0)^{-1}Q\right).$$
(These formulae may be extracted from \cite{Harn94} equations 1.4, 2.23, 2.24, where we have replaced the symbols 
$F,G^T,A,Y,\lambda$ by $P,-Q,T_\infty, T_0,z$ respectively.)
This may be obtained from the Fourier--Laplace transform
as follows (this is ``well-known''\footnote{In particular thanks are due to J. Harnad for telling me about \cite{Harn94} in Luminy in 1996, whilst I was trying to understand the role of the Fourier--Laplace transform in Dubrovin's work \cite{Dub95long} on Frobenius manifolds.}), cf. \cite{BJL81, sang-wood, k2p}:
A local solution $v$ of the first operator may be written as:
$$\frac{dv}{dz} =  T_0v + Qw,\qquad (z-T_\infty)w = Pv. $$
Then replacing $\frac{d}{dz}$ by $z$ and $z$ by 
$-\frac{d}{dz}$ yields
$$zv =  T_0v + Qw,\qquad-\frac{dw}{dz}-T_\infty w = Pv $$
so $\frac{dw}{dz} +T_\infty w + P(z-T_0)^{-1}Qw$, which says $w$  is a solution of the Harnad dual operator. 
(One may interpret this in terms of presentations of Weyl algebra modules, as we do in the body of this article.)
Isomonodromic deformations of such operators are governed by the JMMS equations.
If $P,Q$ solve the JMMS equations \eqref{eq: altJMMSintro} 
on $\IB'\subset \IB$  then the connection
$$\Omega = 
 d(zT_0)  + Qd\log(z-T_\infty)P + \wt{QP}
$$
on $W_0\times (\IP^1\times\IB')$ is flat, and the $\partial/\partial z$ component of it is as in \eqref{eq: harnadop1}.
Now consider the permutation
$$
(W_0,W_\infty,P,Q,T_0,T_\infty)
\mapsto
(W_\infty,W_0,Q,-P,-T_\infty,T_0)$$
of the data obtained from the Fourier--Laplace transform.
This again constitutes  a solution of the JMMS equations (and this is one of the main points of \cite{Harn94}), and so the connection
$$\Omega' = 
 -d(zT_\infty)  - Pd\log(z-T_0)Q - \wt{PQ}
$$
on $W_\infty\times (\IP^1\times\IB')$ is also flat (and its vertical component is the Harnad dual of that above), i.e. the same nonlinear equations govern the isomonodromic deformations of two connections on different rank bundles.
Note that other generalisations  of Harnad duality (different to ours) 
are studied in 
\cite{woodh-duality, yamakawa-mc+hd} (see also \cite{Harn94} \S4).

\section{Leading term computation}\label{apx: leading term}

At $z=\infty$ the connection $\cA$  on $\IP^1$ in \eqref{eq: conn} is formally isomorphic, via $G\flb z^{-1}\frb$, to a connection of the form
\beq \label{eq: formalconn}
 \left(Az+T+\frac{\wh\La}{z}\right)dz 
\eeq
for some $\h=\Lie(H)$ valued holomorphic map $\wh\La = \La + A_1 /z +A_2/z^2 + \cdots$, where $H$ is the centraliser of $A$ and $T$.
We will say $\cA$ is {\em nonresonant} (at $\infty$) if 
$\ad_\La\bigl\vert_\lh\in \End(\lh)$ does not have any nonzero integer eigenvalues. 
If  $\cA$ is nonresonant then we may take $\wh\La=\La$ to be constant.
For example if $\lh$ is a Cartan subalgebra (as in the case considered by Jimbo--Miwa--Ueno \cite{JMU81}) 
then this condition is empty: all of their connections are nonresonant at each irregular singularity.
In the general linear case (as in the body of the text) 
$\lh = \bigoplus_{i\in I\setminus I_\infty} \End(V_i)$ 
and $\La$ has components 
$\La_i\in \End(V_i)$, and the nonresonance condition is that each $\La_i$ has no eigenvalues differing by a nonzero integer.

In this appendix we will suppose at first that $\wh\La = \La$ is constant (which is possible for example if we are in the nonresonant case). Then we will show at the end (in Corollary  \ref{cor: nononres}) how to remove this assumption.

Thus suppose we have $\wh\La = \La$.
Then note that the normal form \eqref{eq: formalconn} may be written as $d\xi$ where 
$\xi:=Az^2/2 + T z +\La\log(z)$.

The aim of this appendix is to determine explicitly the first nontrivial term of any formal isomorphism between \eqref{eq: conn} and the normal form \eqref{eq: formalconn}. 
This is a formal series $\wh g = 1 + g_1/z + g_2/z^2 +\cdots$ such that
$$\wh g[\cA] 
=\left(Az+T+\frac{\La}{z}\right)dz\qquad 
\cA=(Az+B+T)dz  + S(z)dz 
$$
where $Sdz = (QP/z + S_2/z^2+\cdots)dz$ is the Laurent expansion at $z=\infty$ of $Q(z-T_\infty)^{-1}Pdz$, and the square brackets denote the gauge action; $g[\cA] = g\cA g^{-1} + (dg)g^{-1}$. 
In general such $\wh g$ will not be unique.
Our aim is to compute $g_1$ and $\Lambda$ in terms of $A,B,\wh T, P,Q$.
The strategy is as follows.

Write $\g = \End(U_\infty)$ and let $\g^\circ = \Im(\ad_A)\subset \g$ and let 
$\h_1=\ker(\ad_A)$  so that, as vector spaces 
$$\g = \g^\circ\oplus\h_1.$$
Then write 
$\h_1^\circ= \Im(\ad_T\bigl\vert_{\h_1})\subset \h_1$ and let 
$\h=\ker(\ad_T\bigl\vert_{\h_1})$ so that 
$$\h_1 = \h_1^\circ\oplus\h,$$
and we have a nested sequence of Lie algebras $\g\supset \h_1\supset \h$.
Write $\delta:\g\to\h_1$ for the projection onto $\h_1$ along 
$\g^\circ$ and write
$\pi_\h:\h_1\to\h$ for the projection onto $\h$ along $\h_1^\circ$. 

We will find $\wh g$ as a product 
$$\wh g = \wh g_3\wh g_2 \wh g_1$$
where 
$$\wh g_3 = \left(\cdots e^{h_2/z^2}e^{h_1/z}\right),\quad
\wh g_2 = \left(\cdots e^{Y_2/z^2}e^{Y_1/z}\right),\quad
\wh g_1 = \left(\cdots e^{X_2/z^2}e^{X_1/z}\right)$$
where $X_i\in \g^\circ, Y_i\in \h_1^\circ,  h_i\in \h$ for all $i$, 
and where $\wh g_1$ moves $\cA$ into 
$\h_1$, $\wh g_2$ moves the result into $\h$ and then $\wh g_3$ removes the remaining terms which are not singular at $z=\infty$. Note that we will have 
$$g_1 = h_1 + Y_1 + X_1.$$

To simplify notation we will write $X=X_1$.
\begin{lem} \label{lem: g1 calcn}
The following formulae hold: 
$$X= \ad_A^{-1}(B), \qquad (1+\ad_\La)h_1 = \Lambdareplacement_2$$ 
$$
\La = \pi_\h(R), \quad 
Y_1 = \ad^{-1}_T(R-\La),\quad
\Lambdareplacement_2 = \pi_\h\left(R_2 + [Y_1,R]/2\right)$$
$$R=\delta\left(QP + [X,B]/2 \right)$$
$$\pi_\h(R_2) = \pi_\h\left(QT_\infty P + [X,QP] + \ad_X^2(T)/2 + \ad_X^2(B)/3\right).$$
\end{lem}

This will be proved below.
One application of these formulae is to write down the contribution at infinity to the extension to our context of the Jimbo--Miwa--Ueno one-form (used to define the $\tau$ function):
\beq\label{eq: jmu at infty}
\varpi_\infty = 
\Res_\infty \tr \left(\wh g \frac{\partial(\wh g^{-1})}{\partial z} dT z dz\right)
=-\Res_\infty\tr (g_1 dT dz/z) = \tr(g_1 dT).
\eeq
(Beware of the sign in the definition of $T^{(\infty)}$ in \cite{JMU81} (2.17), also reflected in the sign after (2.19).)
Using the formulae, this is
$$\tr(h_1 dT) = \tr(\Lambdareplacement_2 dT) = \tr(R_2 + [Y_1,R]/2)dT=$$
$$\tr(QT_\infty PdT) +\tr ([X,QP]dT) - \tr([X,T][X,dT])/2 - 
\tr([X,B][X,dT])/3+ \tr(\wt R R)/2$$
where $\wt R = \ad_T^{-1}([dT,R])$. 
If there are no simple poles, this expression remains valid for any reductive structure group 
(replacing $\tr(AB)$ by an invariant inner product $\langle A,B\rangle$ throughout). 
In the general linear case it may be rewritten as:
$$\tr(QT_\infty PdT) +\tr ([X,QP]dT) - \tr(XTXdT) + \tr(X^2TdT) - 
\tr(XBXdT)$$$$
+ \tr(\wt{QP}\delta(QP))/2
+ \tr(\wt{QP}\delta(XB))
+ \tr(\wt{XB}\delta(XB))/2.
$$
After re-ordering the terms by degree this is \eqref{eq: varpi-infty}.
The only tricky parts in deriving this are 1) to observe $\tr(\wt M N) = \tr(M\wt N)$ in general, and 2) to note that 
$$
\tr((X^2B+BX^2)dT) =
\tr((X^2[A,X]+[A,X]X^2)dT)= 
\tr((X^2AX-XAX^2)dT)$$
as the remaining terms are $[A,X^3]dT$ which is traceless.
In turn this is $\tr(X[A,X]XdT)=-\tr(XBXdT)$; this leads to the equality
$\tr[X,B][X,dT]/3 = \tr XBXdT$.

\pfms (of Lemma \ref{lem: g1 calcn}).
Suppose that 
$$(\wh g_2\wh g_1)[\cA] 
=\left(Az+T+\sum_1^\infty \Lambdareplacement_i z^{-i}\right)dz
$$
with $\Lambdareplacement_i\in \h$. (Define $\La=\Lambdareplacement_1$.)

\begin{lem}\label{lem: h1} 
$(1+\ad_\La)h_1 = \Lambdareplacement_2$ 
\end{lem}
\pf
Here and below the formula 
$\Ad_{\exp(X)}{Y} = \sum_{n=0}^\infty \ad^n_X(Y)/n!$ will be very useful. Thus 
$$e^{(h_1/z)}\left[\sum \Lambdareplacement_i z^{-i} \right] = 
\La/z +(\Lambdareplacement_2 + [h_1,\La] - h_1)/z^2 + O(1/z^3).$$
The desired formula then arises from the vanishing of the second term.
\epf
Note that under the nonresonance conditions on $\Lambda$ 
the operator $(1+\ad_\La)$ is invertible and so $h_1$ is determined by $\Lambdareplacement_2$.
Now suppose that 
$$\wh g_1[\cA] 
=\left(Az+T+\sum_1^\infty R_i z^{-i}\right)dz
$$
with $R_i\in \h_1$, and write $R=R_1$ for the residue term.
\begin{lem} The elements $Y_1,\La,\Lambdareplacement_2\in\h$ are uniquely determined by $R$ and $\pi_\h(R_2)$:
$$
\La = \pi_\h(R), \quad 
Y_1 = \ad^{-1}_T(R-\La),\quad
\Lambdareplacement_2 = \pi_\h\left(R_2 + [Y_1,R]/2\right).$$
\end{lem}
\pf
We have
$e^{Y_1/z}\left[Az + T + \sum R_i z^{-i}\right]=$
$$Az +  T + \left(R+[Y_1,T]\right)/z + 
\left(R_2 + [Y_1,R] + [Y_1,[Y_1,T]]/2 -Y_1\right)/z^2 +\cdots$$
so $\La$ is the $\h$ component of $R$ and $Y_1$ is defined so as to kill the rest of $R$ (i.e. the component in $\h_1^\circ$). 
In turn $Y_2$ will be defined to kill the $\h_1^\circ$ component of the displayed coefficient of $z^{-2}$, so that $\Lambdareplacement_2$ will be the $\h$ component. Given that $[Y_1,T] = \La-R$ and that $\pi_\h[Y_1,\La]=\pi_\h(Y_1)=0$ we obtain the stated formula for $\Lambdareplacement_2$.
\epf

Thus we finally need to compute $R\in \h_1$ and $\pi_\h(R_2)$ from $\cA$. We will write $X=X_1$.
$$e^{X/z}\left[Az + B+T + \sum S_i z^{-i}\right]=$$

$$Az + \left(B + T + [X,A]\right) +
\left(S_1 + [X,B+T] + \ad_X^2(A)/2 \right)/z+ $$ 
$$\left(S_2 + [X,S_1] + \ad_X^2(B+T)/2 + \ad_X^3(A)/6 -X\right)/z^2 +\cdots$$
Thus $$X=\ad_A^{-1}(B)$$ and $R$ will be the $\h_1$ component of the residue, which simplifies to
$$R=\delta\left(QP + [X,B]/2 \right).$$
Applying $e^{X_2/z^2}[\cdot]$ to the above expression for 
$e^{X/z}[\cA]$ just adds terms $[X_2,A]/z + [X_2,T]/z^2+\cdots$, 
 and the $\h$ component of $[X_2,T]$ is zero so we deduce that 
$$\pi_\h(R_2) = \pi_\h\left(S_2 + [X,QP] + \ad_X^2(B+T)/2 + \ad_X^3(A)/6\right)$$
$$   = \pi_\h\left(QT_\infty P + [X,QP] + \ad_X^2(T)/2 + \ad_X^2(B)/3\right)$$
since $S_2= QT_\infty P$ and $[A,X]=B$.
This completes the derivation of  the desired formulae.
\epfms

Finally we will show that in fact the nonresonance conditions are unnecessary. 

\begin{lem} \label{lem: allow resonance}
There are formal transformations $\wh g\in G\flb 1/z\frb$ (with constant term $1$) putting $\cA$ in the normal form
$$\left(Az+T+\frac{\La + A_1/z +A_2/z^2+\cdots}{z}\right)dz $$
with $\La, A_i\in \lh$ and $[\La^s,A_i] = -i A_i$, where $\La^s$ is the semisimple part of $\La$ (i.e. $A_i$ is in the generalised eigenspace of $\ad_\La$ in $\lh$ with eigenvalue $-i$).
\end{lem}

\pf
We proceed as above to put the connection in $\h$. Then it is effectively a logarithmic connection (since $A$ and $T$ commute with everything in $\h$), and so we can use the usual Gantmacher--Levelt theory (\cite{sabbah-engbook} 2.20).
\epf

Now let $\wh g=\wh g_3\wh g_2 \wh g_1$ be any of these formal isomorphisms. Lemma  \ref{lem: h1} is then modified to become
\begin{lem} \label{lem: h1 resonant}
$(1+\ad_\La)h_1 = \Lambdareplacement_2 - A_1$ \end{lem}
\pf 
As in Lemma  \ref{lem: h1} we look at the coefficient $(\Lambdareplacement_2+[h_1,\La]-h_1)$ of $1/z^2$, but in the resonant case we may not be able to choose $h_1$ so that this vanishes. Rather, the best that can be done is to decompose $\lh$ into the generalised eigenspaces of $\La$ and define $A_1$ to be the component of $\Lambdareplacement_2$ in the generalised eigenspace with eigenvalue $-1$, so that $\Lambdareplacement_2-A_1$ has no component in this subspace, and we may define
$$h_1 = (1+\ad_\La)^{-1}(\Lambdareplacement_2-A_1),$$ 
since $(1+\ad_\La)$ is invertible on the direct sum of all the other generalised eigenspaces.
\epf

However this does not affect the expression 
\eqref{eq: varpi-infty}  for the Hamiltonian one form:
\begin{cor}\label{cor: nononres}
The expression for
$\varpi_\infty = \tr(g_1dT)$ 
in terms of $B,P,Q$ etc is unchanged.
\end{cor}
\pf As before it equals $\tr(h_1dT)$, and this is still equal to $\tr(\Lambdareplacement_2 dT)$, since $\Tr(A_1 dT) = 0$. Indeed for example $A_1=[A_1,\La^s]$ and $\La^s$ commutes with $dT$. The rest of the formulae are unchanged.
\epf

Note that in the resonant case $\xi$ is not so well defined, and so we work with the expression $zdT$ directly rather than ``$d_\IB\xi$'' in the expression \eqref{eq: jmu at infty} used to define the $\tau$ function.

\section{Relating orbits} \label{apx: relating orbits}

Suppose $P:U\to V$ and $Q:V\to U$ are linear maps between  two finite dimensional complex vector spaces $U,V$
such that $P$ is surjective and $Q$ is injective.
In this appendix we will recall the exact relation between the Jordan normal forms of $QP\in \End(U)$ and $PQ\in \End(V)$.
This is often used in the relation between graphs and orbits of matrices (cf. e.g \cite{Kraft-Procesi-InvMath79, nakaj-duke94, CB-Shaw}).
Let $\cO\subset \End(U)$ be the orbit of elements conjugate to $QP$ and let  $\breve \cO\subset \End(V)$ be the orbit of $PQ$.

Recall that giving a Jordan form is equivalent to giving a partition (i.e. a Young diagram) $\pi_s$
for each complex number $s$,
so that  $\pi_s$ specifies the sizes of the Jordan blocks corresponding to the eigenvalue $s\in \IC$. 
For example the partition $\pi_0=(2,2,1)$ specifies 
the $5\times 5$ rank $2$ nilpotent matrix with three Jordan blocks of size $2,2$ and $1$ respectively
(and it corresponds to the Young diagram with three rows of lengths $2,2,1$).

\begin{prop}\label{prop: orbit relations}
Let $\{\pi_s\}$ be the  partitions giving the Jordan form of $\cO$ and let $\{\breve \pi_s\}$ be the  partitions giving the Jordan form of $\breve \cO$.
Then 
$\breve \pi_s = \pi_s \text{ if $s\neq$ 0, and}$
$\breve \pi_0$ is obtained from 
$\pi_0$ by deleting the first (i.e. longest) column of $\pi_0$.
(In other words each part of $\pi_0$ is decreased by one to obtain $\breve \pi_0$.) 
\end{prop}

\pf
Decompose $U=\bigoplus U_s, V = \bigoplus V_s$ into the generalised eigenspaces of $QP, PQ$ respectively.
Then it is easy to see $P$ maps $U_s$ to $V_s$ and $Q$ maps $V_s$ to $U_s$
and that these are isomorphisms if $s\neq 0$.
This establishes the result for $s\neq 0$ and we can reduce to the case where both $QP, PQ$ are nilpotent. 
For the nilpotent case, let $p_1 \ge p_2 \ge \cdots  \ge 0$ denote the lengths of the columns of $\pi_0$ (i.e. the sizes of the parts of the partition dual to $\pi_0$).
Then the relation between the orbit $\cO$ of $QP$ and these lengths is simply expressed as
$$\rank (QP)^j = \sum_{i>j} p_i$$
for any $j$. In particular, as is well-known, 
these ranks determine $\pi_0$.
Then to determine the partition $\breve \pi_0$ 
corresponding to $PQ$ we just note:
\begin{align*}
\rank(PQ)^j = \dim (PQ)^j(V) &= \dim Q(PQ)^j(V)\qquad &&\text{as $Q$ is injective}\\
&= \dim Q(PQ)^jP(U)\quad &&\text{as $P$ is surjective}\\
&= \dim (QP)^{j+1}(U)=\rank(QP)^{j+1}. &&\text{}
\end{align*}
Thus $\breve \pi_0$ has columns of lengths $p_2, p_3,\ldots$, as expected.
\epf

In particular, given $U,V$, the orbit $\breve \cO$ is uniquely determined by $\cO$, and
$\cO$ is uniquely determined by $\breve \cO$ (one just adds a column of length $\dim(U)-\dim(V)$ to the Young diagram of $\breve \pi_0$ to obtain 
$\pi_0$).

\section{Notation summary}

$\bP = \IC \union \{\infty\}$ ``Fourier sphere''. $J\subset  \bP$ a finite subset.

$V = \bigoplus_{j\in J} W_j,$

$\de = \sum_{j\in J} \de_j:\End(V)\to \bigoplus_{j\in J} \End(W_j)\subset \End(V)$

$U_j = V\ominus W_j := \bigoplus_{i\in J\setminus\{j\}} W_j$, so $V = W_j\oplus U_j$ for all $j\in J$.

$U_i := U_j$ for any $i\in I_j$.

$\ga  = \bmx C & P \\ Q & B+T \emx,$ 
$\wh T  = \de(\ga)=\bmx C &  \\  & T \emx\in \End(W_\infty\oplus U_\infty)$

$T_j = \de_j(\wh T)\in \End(W_j)$ semisimple, $T_\infty=C$

$W_j = \bigoplus_{i\in I_j} V_i $, eigenspaces of $T_j$, $I= \bigsqcup I_j$, so $V = \bigoplus_{i\in I} V_i$.

$\wh T = \sum_{i\in I} t_i \Id_i$ where $t_i\in \IC, \Id_i$ idempotent for $V_i\subset V$

$\Ga = \ga^\circ = \ga-\de(\ga) = \bmx 0 & P \\ Q & B \emx  \in \End(W_\infty\oplus U_\infty)$

$\Xi = \phi(\Ga)= \bmx 0 & -P \\ Q & X \emx \in \End(W_\infty\oplus U_\infty)$ where $X = \ad_A^{-1}(B)\in \End(U_\infty)$

$\IM = \bigoplus_{i\ne j\in J} \Hom(W_i,W_j)= \End(V)^\circ$ symplectic  space dependent on $\ba:J \hookrightarrow \bP$

$Q_i = \Gamma\circ \iota_i:V_i\to U_i, \ P_i = -\pi_i\circ\Xi:U_i\to V_i$

$R_i=Q_iP_i\in \End(U_i), \ \La_i:=-P_iQ_i\in \End(V_i)$

\renewcommand{\baselinestretch}{1}              %
\normalsize
\bibliographystyle{amsplain}    \label{biby}
\bibliography{../thesis/syr} 

\vspace{0.5cm}   
\'Ecole Normale Sup\'erieure et CNRS, 
45 rue d'Ulm, 
75005 Paris, 
 France

www.math.ens.fr/$\sim$boalch\qquad \qquad \qquad \qquad \qquad \quad\,

boalch@dma.ens.fr 
\end{document}